\documentclass[11pt, oneside]{amsart}

\usepackage{amsmath,amsthm,amssymb,hyperref, mdframed, eucal}
\usepackage{mathrsfs} 
\usepackage{lipsum}

\usepackage{tikz, tikz-cd}
\usepackage[all]{xy}
\usetikzlibrary{arrows,decorations.markings}
\usetikzlibrary{backgrounds,shapes}
\usetikzlibrary{patterns,hobby}
        \usepackage{pgfplots}
        \pgfplotsset{compat=1.6}
\tikzset{>=latex}
\usepackage{subfig}
\usepackage{booktabs}
\usepackage{subfiles}
\usepackage{caption}
\usepackage{makecell}
\usepackage[new]{old-arrows}

\usepackage{color}
\usepackage{hyperref}
\hypersetup{
    colorlinks=true, 
    linktoc=all,     
    linkcolor={mauve},  
    citecolor={plum},
    filecolor={plum},
    urlcolor={plum}
}

\definecolor{smoked}{RGB}{216, 212, 204}
\definecolor{mauve}{RGB}{200, 55, 171}
\definecolor{apricot}{RGB}{250, 144, 4}
\definecolor{sky}{RGB}{66, 169, 244}
\definecolor{plum}{RGB}{76, 0, 102}

\definecolor{lightmauve}{RGB}{232, 173, 220}
\definecolor{lightapricot}{RGB}{253, 211, 155}
\definecolor{lightsky}{RGB}{178, 221, 251}
\definecolor{lightplum}{RGB}{184, 153, 192}

\usepackage{geometry} 
\usepackage[parfill]{parskip}

\newcommand{\Z}{{\mathbb Z}}
\newcommand{\C}{{\mathbb C}}
\newcommand{\R}{{\mathbb R}}
\newcommand{\HH}{{\mathbb H}}

\DeclareMathOperator\Aut{Aut}
\DeclareMathOperator\Out{Out}

\DeclareMathOperator\End{End}

\renewcommand{\Im}{\mathrm{Im}}

\DeclareMathOperator{\PSL}{\mathrm{PSL}}
\DeclareMathOperator{\psl}{\PSL_2\R}
\DeclareMathOperator{\SL}{\mathrm{SL}}
\DeclareMathOperator{\GL}{\mathrm{GL}}

\DeclareMathOperator{\Tr}{\mathrm{Tr}}

\theoremstyle{plain}                    
\newtheorem{thm}{Theorem}[section]
\newtheorem{thma}{Theorem}

\newtheorem{lem}[thm]{Lemma}
\newtheorem{prop}[thm]{Proposition}
\newtheorem{cor}[thm]{Corollary}
\newtheorem{conj}[thm]{Conjecture}
\newtheorem{fact}[thm]{Fact}

\theoremstyle{definition}
\newtheorem{defn}[thm]{Definition}
\newtheorem{ex}[thm]{Example}
\newtheorem{question}[thm]{Question}
\newtheorem{problem}[thm]{Problem}

\theoremstyle{remark}
\newtheorem{rmk}[thm]{Remark}
\newtheorem{claim}[thm]{Claim}

\numberwithin{equation}{section}

\newenvironment{proofclaim}
 {\proof}
 {\endproof}

\DeclareMathOperator{\Mod}{Mod}
\DeclareMathOperator{\eu}{eu}
\DeclareMathOperator{\arccosh}{arccosh}
\DeclareMathOperator{\sech}{sech}
\DeclareMathOperator{\Tol}{Tol}
\DeclareMathOperator{\Vol}{Vol}
\DeclareMathOperator{\id}{id}
\DeclareMathOperator{\dev}{dev}
\DeclareMathOperator{\hol}{hol}
\DeclareMathOperator{\Hyp}{Hyp}
\DeclareMathOperator{\Teich}{Teich}
\DeclareMathOperator{\Sym}{Sym}
\newcommand{\quotient}[2]{{\raisebox{.2em}{$#1$}\left/\raisebox{-.2em}{$#2$}\right.}}

\title[Branched hyperbolic surfaces]{On the symplectic geometry of branched hyperbolic surfaces in genus two}

\author{Gianluca Faraco}
\address[G.~Faraco]{Dipartimento di Matematica e Applicazioni U5, Universita` degli Studi di Milano-Bicocca, Via Cozzi 55, 20125 Milano, Italy.}
\email{gianluca.faraco@unimib.it}
\email{gianluca.faraco.math@gmail.com}

\author{Arnaud Maret}
\address[A.~Maret]{Université de Strasbourg, IRMA, 7 rue Descartes, 67000 Strasbourg, France.}
\email{maret.arnaud@unistra.fr}

\date{\today}

\begin{document}

\begin{abstract} 
We construct analogues of Fenchel--Nielsen coordinates on an open and dense subset of the space of holonomies of branched hyperbolic structures on a closed genus-2 surface. We show that these coordinates satisfy an analogue of Wolpert’s magic formula, and thus provide Darboux charts for the Goldman symplectic form. To this end, we revisit the parametrization of hyperbolic structures on a one-holed torus and describe a simple polygonal model that makes both length and twist parameters transparent. Gluing two such polygons leads to the notion of bow-tie representations of a genus-2 surface group. We prove that bow-tie representations account for most holonomies of branched hyperbolic structures, though not all: for example, Le Fils’ pentagon representations form a real codimension-2 family of holonomies lying outside the bow-tie locus.
\end{abstract}

\maketitle

\maketitle

\section{Introduction}\label{sec:intro}
\subsection{A story of bow-ties and pentagons}\label{sec:motivations_results}
Given a closed and oriented surface $S$ of genus $g\geq 2$, the space of conjugacy classes of representations $\pi_1S\to\psl$ is known as the \emph{character variety} of $S$ and $\psl$ and is denoted by $\chi(S,\psl)$ (Definition~\ref{defn:character-varieties}). Goldman proved that $\chi(S,\psl)$ has $4g-3$ connected components indexed by an integer $|e|\leq 2g-2$ called \emph{Euler number} (Section~\ref{sec:character-varieties-components})~\cite{goldman-components}. Each component with $|e|\geq 1$ is a smooth symplectic manifold of real dimension $6g-6$ (Section~\ref{sec:character-variety-symplectic-nature})~\cite{goldman-symplectic}.

\begin{figure}[h]
\centering
\begin{tikzpicture}[scale=11/13]
\node[anchor = south west, inner sep=0mm] at (0,0) {\includegraphics[width=11cm]{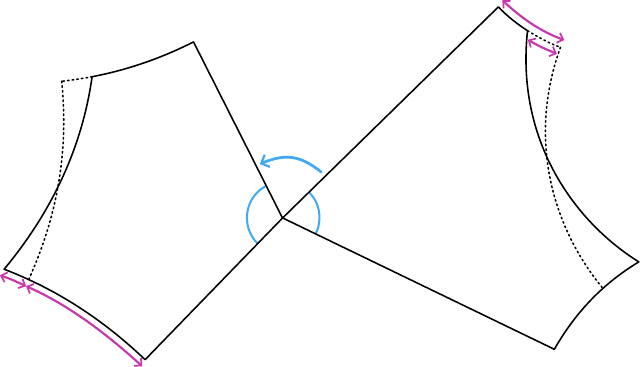}};
\node[sky] at (6,4.5) {\small $\gamma$};
\node[sky] at (7,3.15) {\small $\beta/2$};
\node[sky] at (4.05,3.15) {\small $\pi-\beta/2$};
\node[mauve] at (1.4,.6) {\small $\ell_1/2$};
\node[mauve] at (0.0,1.32){\small $\tau_1/2$};
\node[mauve] at (11.1,7.3) {\small $\ell_2/2$};
\node[mauve] at (11.92,6.37) {\small $\tau_2/2$};
\end{tikzpicture}
\caption{Geometric parametrization of bow-tie representations.}\label{fig:bow-tie-intro}
\end{figure}

This article focuses on the case where $S$ is a genus-2 surface. We introduce \emph{bow-tie} representations (Definition~\ref{defn:bow-tie-representations}) as the representations $\pi_1S\to\psl$ with Euler number $\pm 1$ that map a separating simple closed curve on $S$ to an elliptic element. Bow-tie representations admit a concrete geometric parametrization: they correspond to hyperbolic polygons obtained by attaching two pentagons along a single vertex, giving rise to the characteristic ``bow-tie'' shape, see Figure~\ref{fig:bow-tie-intro}. Those polygons have six parameters of deformations: two lengths $\ell_1,\ell_2>0$ and one angle $\beta\in (0,2\pi)$, as well as two dual twists $\tau_1,\tau_2\in \R$ and one dual angle $\gamma\in \R/2\pi\Z$ (Section~\ref{sec:parametrization-bow-tie-representations}). The parametrization depends on chosen topological data on $S$, including a separating curve. Our main result is the following.
\begin{thma}[Theorems~\ref{thm:coordinates-bow-tie-representations} \&~\ref{thm:wolpert-formula-bow-tie-representations}]\label{thm:bow-tie-representations-intro}
For every separating simple closed curve $c$ on $S$, the open set of conjugacy classes of bow-tie representations mapping $c$ to an elliptic element is diffeomorphic to two disjoint copies of $\R_{>0}\times\R\times (0,2\pi)\times \R/2\pi\Z\times \R_{>0}\times \R$, globally parametrized by the coordinates $(\ell_1,\tau_1,\beta,\gamma,\ell_2,\tau_2)$.  Moreover, the 2-form
\[
-d \ell_1\wedge d\tau_1 -d\beta\wedge d\gamma -d\ell_2\wedge d\tau_2
\]
coincides with the Goldman symplectic form.
\end{thma}

If $\mathcal{B}$ denotes the set of all conjugacy classes of bow-tie representations and $\mathcal{B}_c\subset \mathcal{B}$ the subset of those mapping $c$ to an elliptic element, then Theorem~\ref{thm:bow-tie-representations-intro} provides an atlas of Darboux charts for $\mathcal{B}$ indexed on the curves $c$ and some extra topological data (see the discussion at the beginning of Section~\ref{sec:parametrization-bow-tie-representations} for more details). Thanks to the ergodicity result of Marché--Wolff (Theorem~\ref{thm:goldman-ergocicity-genus-2}), we know that $\mathcal{B}$ is an open and dense subset of the two components of $\chi(S,\psl)$ with Euler number $\pm 1$. In other words, the Darboux charts of Theorem~\ref{thm:bow-tie-representations-intro} are almost a global parametrization of those two components.

Theorem~\ref{thm:bow-tie-representations-intro} establishes an analogue of the classical Fenchel--Nielsen coordinates for Teichmüller spaces (Section~\ref{sec:character-varieties-components}). Theorem~\ref{thm:bow-tie-representations-intro} also provides the analogue of Wolpert's magic formula~\cite{wolpert-formula} for our coordinates. Unlike Teichmüller space which is contractible, the two components of $\chi(S,\psl)$ with Euler number $\pm 1$ have the topology of a rank-2 complex vector bundle over $S$ by a result of Hitchin~\cite{hitchin-self-duality}. In particular, no single global Darboux chart for the Goldman symplectic can exist on those components. Nevertheless, Theorem~\ref{thm:bow-tie-representations-intro} provides a countable collection of Darboux charts that together cover an open and dense subset of each component.

Beyond bow–tie representations, we shall also examine a further class of representations: the so-called \textit{pentagon representations} (Definition \ref{defn:pentagon-representation}). They were introduced by Le Fils in~\cite{thomas} to show the existence of non-elementary surface group representations into \(\text{PSL}_2\,\mathbb C\) that do not admit a Schottky decomposition, that is, a pants decomposition for which the restriction of the representation to each pair of pants is an isomorphism onto a Schottky group. Our interest in pentagon representations arises from the fact that they have Euler number \(\pm 1\) and they are not bow–tie representations. In fact, they map all separating simple closed on $S$ to hyperbolic elements of $\psl$ (Lemma~\ref{lem:properties-pentagon}). Conjugacy classes of pentagons representations are precisely the fixed points of the \emph{hyperelliptic involution} (Section~\ref{sec:hyperelliptic-involution}) in the two components of $\chi (S,\psl)$ with Euler number $\pm 1$ (Lemma~\ref{lem:pentagon-iff-fixed-by-hyperelliptic-involution}).

We denote by \(\mathcal{P}\) the set of conjugacy classes of pentagon representations. It is a real codimension-$2$ symplectic subspace of the two components of $\chi(S,\psl)$ with Euler number $\pm 1$. It turns out that \(\mathcal{P}\) consists of twelve connected components (this is mostly a consequence of Theorem~\ref{thm:coordinates-pentagon-representations}), six of which have Euler number \(1\) and the other six have Euler number $-1$. The six components with the same Euler number are images of each other by Dehn twists (Table~\ref{tab:6-components-pentagons}).

\begin{figure}[h]
\centering
\begin{tikzpicture}
\node[anchor = south west, inner sep=0mm] at (0,0) {\includegraphics[width=9cm]{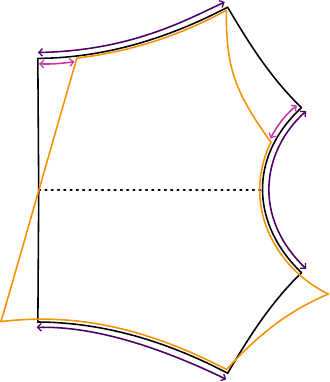}};

\node[plum] at (3.3,.75) {\small $\ell_1/2$};
\node[plum] at (3.3,9.6) {\small $\ell_1/2$};
\node[plum] at (7.85,5.2) {\small $\ell_2/2$};
\node[mauve] at (7.5,7.3) {\small $\tau_2$};
\node[mauve] at (1.5,8.4) {\small $\tau_1/2$};
\end{tikzpicture}
\caption{Geometric parametrization of pentagon representations.}
\label{fig:pentagon-intro}
\end{figure}

As for bow–tie representations, pentagon representations also admit an explicit geometric description: they correspond to (sometimes non-convex) hyperbolic hexagons\footnote{An aptonym is a rare miracle.} obtained by deforming right-angled hexagons in a particular way, as for instance the orange hexagon on Figure~\ref{fig:pentagon-intro}. Those hexagons are parametrized by two lengths \(\ell_1,\ell_2>0\) and two dual twists \(\tau_1,\tau_2\in\mathbb R\) (Section~\ref{sec:parametrization-pentagon-representations}).

\begin{thma}[Theorems~\ref{thm:coordinates-pentagon-representations} \&~\ref{thm:wolpert-formula-pentagon-representations}]\label{thm:pentagon-representations-intro}
Every connected component of \(\mathcal P\) is diffeomorphic to $\R_{>0}\times\R\times\R_{>0}\times\R$ and the parameters \((\ell_1,\tau_1,\ell_2,\tau_2)\) provide a system of global coordinates. Moreover, the 2-form
\[ -d\ell_1\wedge d\tau_1-d\ell_2\wedge d\tau_2
\] coincides with the restriction of the Goldman symplectic form.
\end{thma}

An intriguing question is whether there exist representations of a genus-\(2\) surface group into \(\psl\) with Euler number \(\pm 1\) that are neither bow--tie nor pentagon representations.

\begin{question}\label{question:bow-tie-and-pentagon-what-else}
Is every representation \(\rho\) of a genus-\(2\) surface group into \(\psl\) with Euler number \(\pm 1\) either a bow--tie or a pentagon representation?
\end{question}

Question~\ref{question:bow-tie-and-pentagon-what-else} can be reformulated as two independent questions concerning representations with Euler number \(\pm 1\): (i) Are pentagon representations characterized by the property that they map every separating simple closed curve to a hyperbolic element? (ii) If \(\rho\) maps some separating simple closed curve to a parabolic element, must it necessarily map another separating simple closed curve to an elliptic element? Viewed in this way, Question~\ref{question:bow-tie-and-pentagon-what-else} may be regarded as a refinement of Bowditch’s question, which was answered affirmatively in genus~\(2\) by Marché--Wolff (Theorem~\ref{thm:bowditch-genus-2}). We do not attempt to answer Question~\ref{question:bow-tie-and-pentagon-what-else} here, and leave it as a direction for future work.

\subsection{Some ideas about the proofs}
In order to prove that the coordinates we introduce to parametrize conjugacy classes of bow-tie and pentagon representations are Darboux coordinates for the Goldman symplectic form (Theorems~\ref{thm:bow-tie-representations-intro} \&~\ref{thm:pentagon-representations-intro}), we follow an approach similar to Wolpert's in~\cite{wolpert-formula}. After proving that the coordinates define a diffeomorphism to either $\R^5\times \R/2\pi\Z$ or $\R^4$ depending on the case, we expand the Goldman symplectic form in those coordinates. The most subtle part of the argument consists in showing that the coefficients in front of the terms $d\ell_i\wedge d\ell_j$ and $d\ell_i\wedge d\beta$ vanish. Here, $\ell_i$ and $\beta$ denote length/angle coordinates. The strategy works as follows.
\begin{enumerate}
    \item We start by showing that those coefficients are constant along the Hamiltonian flow orbits of the functions $\ell_i$ and $\beta$. This means that it is enough to prove that the coefficients vanish at particular points along those orbits. 
    \item We select those particular points as the ones with ``zero twist'' by setting $\tau_i=0$. In the case of bow-tie representations, we also set $\gamma=\beta/2$. Those points enjoy some extra symmetries in the sense that they are fixed by an anti-symplectic involution of the ambient space.
    \item We carefully compute the pullbacks of the exterior derivatives of all our coordinates by this anti-symplectic involution and the right combination of minus signs shows that the desired coefficients vanish.
\end{enumerate}

\subsection{A simple model for one-holed tori}
As a by-product of our work on genus-$2$ surface group representations into \(\psl\) with Euler number \(\pm 1\), we describe a polygonal model for visualizing and measuring Fenchel--Nielsen length and twist coordinates on hyperbolic one-holed tori (Section~\ref{sec:polygonal-model-one-holed-tori}). While the underlying ideas of this model appear to be familiar to some experts and do not seem to come as a surprise, we were unable to find an explicit presentation of this viewpoint in the existing literature. We believe that this model offers a somewhat simpler perspective on twist coordinates than the standard definitions commonly used. We briefly sketch the construction here and refer the reader to Section~\ref{sec:polygonal-model-one-holed-tori} for a detailed account.

\begin{figure}[h]
\centering
\begin{tikzpicture}[scale=1.193]
\node[anchor = south west, inner sep=0mm] at (3.15,3.15) {\includegraphics[width=8cm]{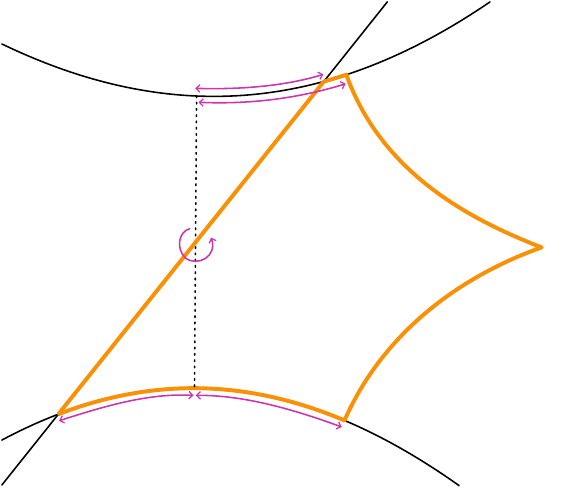}};

\node[mauve] at (4.7,4) {\small $\tau/2$};
\node[mauve] at (6.1,8) {\small $\tau/2$};
\node[mauve] at (6.05, 7.45) {\small $\ell/2$};
\node[mauve] at (6.05,3.95) {\small $\ell/2$};

\node at (3.7,4.15) {\small $X$};
\node at (5.57,4.47) {\small $H$};
\node at (7.1,4.15) {\small $Y$};
\node at (5.3,7.88) {\small $H'$};
\node at (6.9,8.13) {\small $X'$};
\node at (7.43,7.87) {\small $Y'$};
\node[apricot] at (8,4.9) {\small $s$};
\node[apricot] at (8,7.1) {\small $s'$};
\node at (9.65,6) {\small $B$};

\node at (8.1,3.73) {\tiny $\rho(a)$};
\node at (8.4,8.2) {\tiny $\rho(bab^{-1})$};
\node at (4.3,5) {\tiny $\rho(b)$};

\draw[sky, thick] (9,6.12) arc(160:200:.5);
\node[sky] at (8.35,5.95) {\small $\pi-\beta/2$};
\end{tikzpicture}
\caption{How to read Fenchel--Nielsen coordinates and visualize twist flows in a simple geometrical model.}
\label{fig:fenchel-nielsen-coordinates-one-holed-torus-intro}
\end{figure}

Let $\Sigma$ be an oriented torus with one puncture $p$. The fundamental group of $\Sigma$ is a free group on two generators $a$ and $b$ that correspond to homotopy classes of simple loops on $\Sigma$ that intersect once. We consider hyperbolic structures on $\Sigma$ with either a geodesic boundary of length $\lambda>0$, a cusp, or a cone point of angle $\beta\in (0,2\pi)$ at $p$. We let $\ell>0$ be the length of the geodesic representative of $a$ (Definition~\ref{defn:length-coordinate-one-holed-torus}) and $\tau\in \R$ be the dual twist measured using the geodesic representative of $b$ (Definition~\ref{defn:twist-coordinate-one-holed-torus}). The pair $(\ell,\tau)$ are \emph{Fenchel--Nielsen coordinates} of the hyperbolic structure. The holonomies of those hyperbolic structures are representations $\rho\colon\pi_1\Sigma\to\psl$ that map $[a,b]$ to an element of trace $<2$ (the trace of a commutator in $\psl$ is a well-defined real number). This forces $\rho(a)$ and $\rho(b)$ to be two hyperbolic elements of $\psl$ whose axes intersect, say at $X$ (Fact~\ref{fact:tr<2-implies-crrossing-axes}). The image of $X$ by $\rho(b)$ is denoted $X'$. The axes of $\rho(a)$ and $\rho(bab^{-1})$ have a common perpendicular geodesic line that intersect them at $H$ and $H'$. We let $Y$ and $Y'$ be two points on the axes of $\rho(a)$ and $\rho(bab^{-1})$ at distance $\ell/2$ from $H$ and $H'$, and we denote by $s$ and $s'$ the geodesic line through $Y$ and $Y'$ that are perpendicular to the axes of $\rho(a)$ and $\rho(bab^{-1})$. Depending on the nature of $\rho([a,b])$, the geodesic lines $s$ and $s'$ may or may not intersect. Figure~\ref{fig:fenchel-nielsen-coordinates-one-holed-torus-intro} illustrates the case where $s$ and $s'$ intersect at $B$ (see Figure~\ref{fig:polygon-construction-3-types} for the other cases).

\begin{thma}[Proposition~\ref{prop:twist-coordinate-is-geometric-distance} and the discussion beforehand]\label{thm:polygonal-model-intro}
If $\rho([a,b])$ is elliptic or parabolic with fixed point $B$, then $s$ and $s'$ intersect at $B$ making an angle $\pi-\beta/2$ (with $\beta=2\pi$ in the parabolic case). If $\rho([a,b])$ is hyperbolic, then $s$ and $s'$ are distance $\lambda/2$ of each other. Moreover, in all cases, the (signed) distance between $X$ and $H$ is equal to $\tau/2$. 
\end{thma}

\subsection{Goldman's geometrization conjecture}\label{ssec_geometrization_goldmanconj} 
For a general closed and oriented surface $S$ of genus $g\geq 2$, Goldman proved that \emph{all} the representations in the components of $\chi(S,\psl)$ with extremal Euler number $e=\pm(2g-2)$, also known as \emph{Fuchsian representations}, are holonomies of hyperbolic structures on $S$~\cite[Corollary~C]{goldman-components}. He later conjectured that most of the representations in a component with Euler number $1\leq |e|\leq 2g-3$ should be holonomies of \emph{branched hyperbolic structures} on $S$ (Definition~\ref{defn:branched-hyperbolic-structures}) with $n=2g-2-|e|$ branch points (see Conjecture~\ref{conj:geoemtrization-conjecture} for a precise statement). We will call those representations \emph{geometrizable} (Definition~\ref{defn:geometrizable-representations}). When $g=2$, Goldman's geometrization conjecture predicts that \emph{all} representations with Euler number $\pm 1$ are geometrizable. 

While Fuchsian representations are holonomies of a unique hyperbolic structure on $S$, it is possible to deform branched hyperbolic structures in a non-trivial way without changing their holonomy. More precisely, there is a real $2n$-dimensional smooth deformation space of branched hyperbolic structures on $S$ with $n$ branch points that have the same given holonomy. This space is known as an \emph{isomonodromic leaf} (see the discussion at the end of Section~\ref{sec:holonomy}).

Our next results show that pentagon representations are geometrizable and provide a new method to geometrize bow-tie representations. They also describe concrete isomonodromic deformations as a real $2$-dimensional locus of branched hyperbolic structures with prescribed bow-tie holonomy, and a $1$-dimensional such locus for each pentagon holonomy (see also Question~\ref{question:isomonodromic-deformations}).

\begin{figure}[h]
\centering
\begin{tikzpicture}[scale=9/13]
\node[anchor = south west, inner sep=0mm] at (0,0) {\includegraphics[width=9cm]{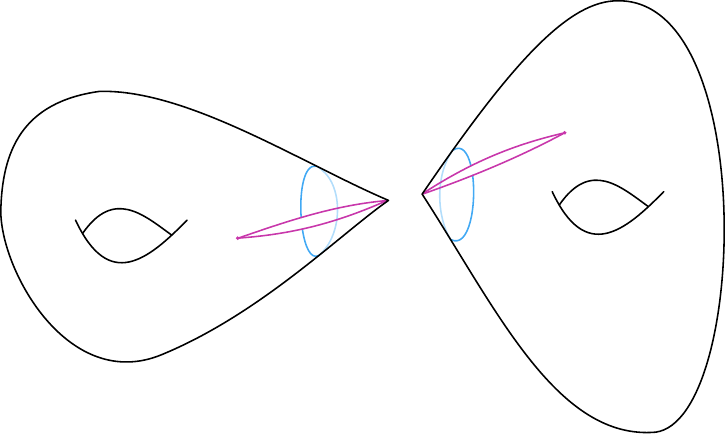}};
\node[sky] at (5.1,4.3) {$\beta_1$};
\node[sky] at (8.85,4.3) {$\beta_2$};
\draw[mauve] (5.8,4.05) edge[out=80, in=130, <->] (9.2,5.18);
\draw[mauve] (5.8,3.7) edge[out=-80, in=-70, <->] (9.3,4.91);
\end{tikzpicture}
\caption{Cut-and-paste procedure on a slit tori assembly.}
\label{fig:cut-and-paste-intro}
\end{figure}

Bow-tie representations were already known to be geometrizable since the work of Mathews~\cite[Theorem~2]{mathews-genus-2}. We geometrize them in a different way using a cut-and-paste procedure on a \emph{slit tori assembly} (Definition~\ref{defn:slit-tori-assemblies})---two hyperbolic one-holed tori with complementary cone singularity ($\beta_1+\beta_2=2\pi$ on Figure~\ref{fig:cut-and-paste-intro}), each equipped with a slit (in mauve) of equal length. Concretely, we cut open both slits and glue the two tori together by identifying the lips of the slits. The outcome is a hyperbolic genus-$2$ surface with one branch point, as desired. This construction is inspired from the samosa assembly construction from~\cite{fenyes-maret} and allows for explicit isomonodromic deformations. 
\begin{thma}[Theorem~\ref{thm:bow-tie-geometrization} and Corollary~\ref{cor:isomonodromic-deformations-bow-tie}]\label{thm:geometrization-bow-tie-intro}  
Bow-tie representations are geometrizable as holonomies of branched hyperbolic structures constructed from slit tori assemblies. Moreover, isomonodromic deformations can be achieved by (i) changing the slits' length, or (ii) by simultaneous rotations of the slits around the cone point on each one-holed torus before gluing.
\end{thma}

\begin{figure}[h]
    \centering
    \begin{tikzpicture}[scale=10/12]
    \node[anchor = south west, inner sep=0mm] at (0,-.2) {\includegraphics[width=10cm]{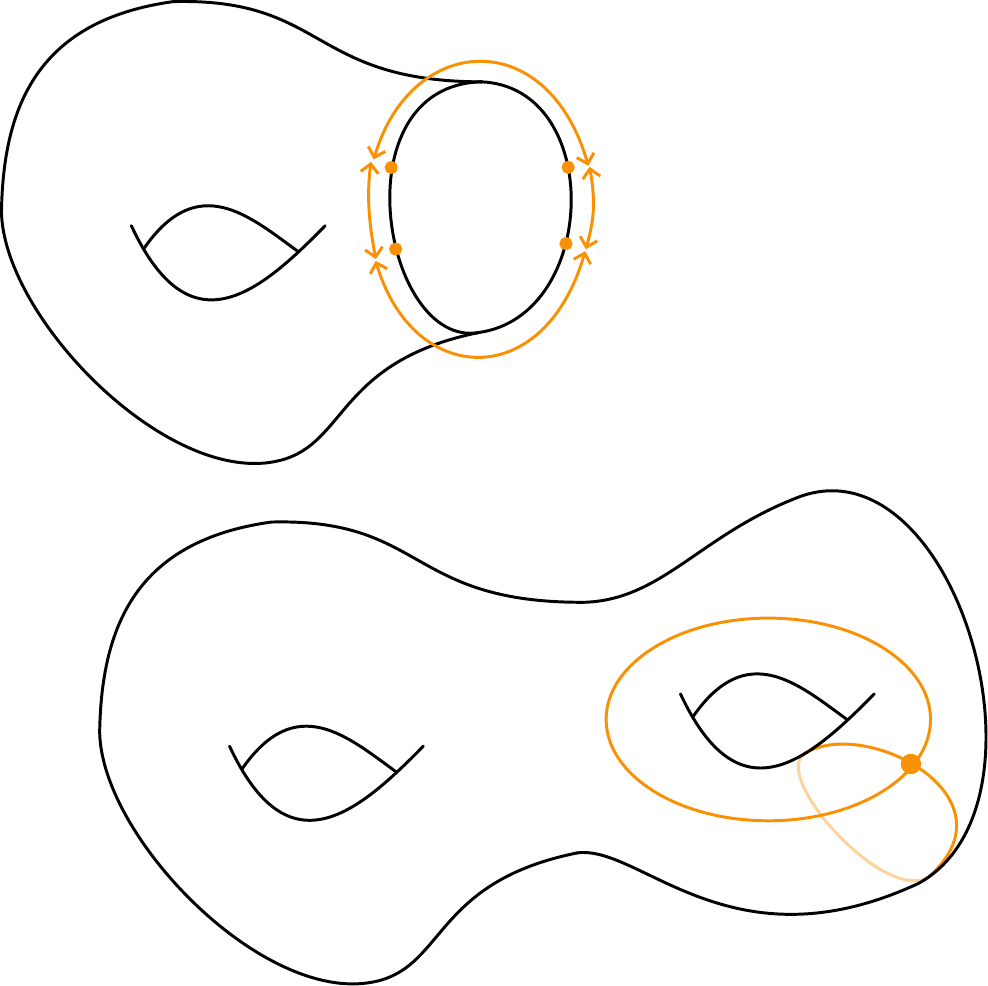}};
    \draw[mauve] (6,7.54) edge[out=80, in=280, <->] (6,10.97);
    \draw[mauve] (4.55,9.1) edge[out=-5, in=185, <->] (7.15,9.1);
    \draw (6.8,7.6) edge[out=-30, in=120, ->] (8.8, 5.5);
    \node[apricot] at (10.9, 2.85) {$4\pi$};
    \end{tikzpicture}
    \caption{Turning a one-holed torus with geodesic boundary into a branched hyperbolic surface of genus $2$.}
    \label{fig:geometrization-pentagon-intro}
\end{figure}

Le Fils showed that pentagon representations are holonomies of \emph{branched projective structures}\footnote{A \emph{branched projective structure} can be defined analogously to a branched hyperbolic structure (Definition~\ref{defn:branched-hyperbolic-structures}), except that charts take values in $\mathbb{CP}^1$ and transition maps are elements of $\mathrm{PSL}_2\C$.} on $S$ with one branch point~\cite[Theorem~1.6]{thomas}. We are interested in geometrizing pentagon representations using branched hyperbolic structures instead. We do so by starting from a hyperbolic one-holed torus with geodesic boundary, and then considering four points on its boundary that divide it into four segments with each pair of opposite segments having the same length. We then crush the boundary by identifying both segments in each pair, imitating the classical gluing of a rectangle's sides to produce a torus, see Figure~\ref{fig:geometrization-pentagon-intro}. This construction results into a branched hyperbolic surface of genus $2$ with one branch point being the image of the four points we picked on the boundary.

\begin{thma}[Theorem~\ref{thm:geometrization_pentagon_reps} and Corollary~\ref{cor:isomonodromic-deformations-pentagon}]\label{thm:geometrization-pentagon-intro}   
Every pentagon representation is equivalent, up to mapping class group action, to one that can be geometrized as holonomies of branched hyperbolic structures constructed by crushing the geodesic boundary of a hyperbolic one-holed torus as we just described. Moreover, changing the distance between the four points while keeping the length of each pair of opposite segments equal is a way of performing isomonodromic deformations.
\end{thma}

If the answer to Question~\ref{question:bow-tie-and-pentagon-what-else} is positive, then Theorem~\ref{thm:geometrization-pentagon-intro}, together with Mathews' result~\cite[Theorem~2]{mathews-genus-2}, would prove Goldman's geometrization conjecture (Conjecture~\ref{conj:geoemtrization-conjecture}) for genus-$2$ surfaces.

\subsection{Charting the next frontiers}
Before developing the results described above, we wish to present some possible further developments of the present work. A first open question has already been stated in Question~\ref{question:bow-tie-and-pentagon-what-else} above, which, we recall, asks whether all genus-$2$ surface group representations into \(\mathrm{PSL}(2,\mathbb{R})\) with Euler number \(\pm 1\) are either bow-tie or pentagon representations.

One other natural direction to pursue is to try to extend Theorems~\ref{thm:bow-tie-representations-intro} and~\ref{thm:pentagon-representations-intro} beyond genus $2$.
Our arguments strongly rely on genus-$2$ features (e.g.~the hyperelliptic involution used in Section \ref{sec:pentagon-representations}, Marché--Wolff's ergodicity results). We expect any construction of Fenchel--Nielsen type coordinates on large domains of the intermediate components of $\chi(S,\psl)$ to be inevitably related to Goldman's ergodicity conjecture (Conjecture~\ref{conj:geoemtrization-conjecture}) and Bowditch's question about existence of simple closed curves with non-hyperbolic images (see the discussion around Theorem~\ref{thm:bowditch-genus-2}).

\begin{problem}\label{problem:Darboux-coordinates}
    For every \(1\leq |e| \leq 2g-3\), describe a system of Darboux coordinates of geometric origin for the Goldman symplectic form on the connected component of \(\chi(S,\psl)\) with Euler number \(e\).
\end{problem}

Problem~\ref{problem:Darboux-coordinates} also makes sense for the component with \(e=0\), but those representations are never holonomies of branched hyperbolic structures (Section~\ref{sec:holonomy}). 

The next question concerns the mapping class group action on $\chi(S,\psl)$ (Section~\ref{sec:dynamics}). When $S$ has genus $2$, Marché--Wolff proved that the action is ergodic on each components with Euler number $\pm 1$ (Theorem~\ref{thm:goldman-ergocicity-genus-2}). This implies that almost every orbit is dense in those components, but not all of them are. Goldman posed the problem of finding sufficient conditions for an orbit to be dense~\cite[Problem~2.7]{goldman-conjectures}. Goldman's problem has been considered in different contexts including representations into $\mathrm{SU}(2)$~\cite{previte-xia-1, previte-xia-2, golsefidy-tamam}, representations into compact abelian Lie groups \cite{yohann-gianluca} and, more recently, for a special class of genus-0 surface group representations into \(\psl\), known as \emph{DT representations}~\cite{minimality}. In this latter case, crucial to the proof is the interpretation of DT representations as \emph{hyperbolic triangle chains}---a polygonal model similar to the ones introduced in this note for bow-ties and pentagon representations. This leads to the following question.

\begin{question}\label{question:MCG-orbit-closures}
    Is it possible to use of the polygonal models for bow-tie and pentagon representations to classify the mapping class group orbit closures of their conjugacy classes?
\end{question}

Note that some mapping class group orbits are discrete. For instance, conjugacy classes of representations with discrete image lie in discrete mapping class group orbits (see e.g.~\cite[Section~3.4]{finite-orbits}). Typical examples in our context include pentagon representations coming from a right-angled pentagon, see Examples~\ref{ex:pentagon-representation} and~\ref{ex:pentagon-representation-P5}. Bow-tie representations may have discrete image too, see Example~\ref{ex:discrete-bow-tie}.

The last question is closely related to geometrization by branched hyperbolic structures and concerns the study of the isomonodromic foliation given by the fibers of the holonomy map (Section~\ref{sec:holonomy}). As we already mentioned, Fuchsian representations are holonomies of a unique hyperbolic structure. This means that they encode enough geometric data to determine a unique geometric structure. Non-Fuchsian representations, on the other hand, do not encode enough geometric data to determine the structure uniquely---sometimes, they are not even geometrizable by branched hyperbolic structures (Example~\ref{ex:non-geometrizable-representations}). However, when they are geometrizable, they are holonomies of infinitely many branched hyperbolic structures which constitute its isomonodromic leaf. We have concretely described portions of those isomonodromic leaves for bow-tie representations (Theorem~\ref{thm:geometrization-bow-tie-intro}) and pentagon representations (Theorem~\ref{thm:geometrization-pentagon-intro}). Isomonodromic leaves are known to be smooth and locally connected, as infinitesimal deformations can be obtained through \textit{movements of branched points} (see \cite[Section 7]{tan} and \cite[Appendix]{CalDerFra14}) which are the geometric counterpart of \textit{Schiffer variations} \cite{NS}. However, it is not known in general whether two branched hyperbolic structures with the same holonomy can be obtained from one another via branch point moves, and the following question naturally arises.

\begin{question}\label{question:isomonodromic-deformations}
    If $\hol^{-1}([\rho])$ denotes the set of isotopy classes of branched hyperbolic structures on $S$ with holonomy $[\rho]$, then is $\hol^{-1}([\rho])$ connected? More generally, can one parametrize $\hol^{-1}([\rho])$ globally using concrete geometric deformations?
\end{question}

The counterpart of Question~\ref{question:isomonodromic-deformations} has been answered positively in several cases in the broader setting of \emph{branched projective structures} on closed surfaces. In the case of structures with quasi-Fuchsian (and, in particular, Fuchsian) holonomy, this was established by Calsamiglia--Deroin--Francaviglia~\cite{CalDerFra14}; see also Billon~\cite{billon} and Francaviglia--Ruffoni~\cite{francaviglia-ruffoni}. At the other end of the spectrum are \emph{translation surfaces }which, like branched hyperbolic structures, constitute special examples of branched projective structures and provide another situation in which the answer is affirmative. In this setting, the holonomy is reducible rather than quasi-Fuchsian, and a positive solution was obtained by Calsamiglia--Deroin--Francaviglia~\cite{calsamiglia-deroin-francaviglia_transfer} for almost every reducible holonomy representation and more recently in~\cite{far-tahar-zhang}.

\subsection{Organization of the paper}
In Section~\ref{sec:background}, we provide a brief recap on character varieties, reviewing some important results and open questions about $\psl$ character varieties of closed surfaces (Section~\ref{sec:character-varieties}). We also provide different perspectives on branched hyperbolic structures on surfaces and Goldman's geometrization conjecture, including the Hodge theoretic viewpoint (Section~\ref{sec:branched-hyperbolic-structures}). 

We revisit some classical results about hyperbolic one-holed tori (e.g.~Fenchel--Nielsen coordinates and Goldman's Hamiltonian twist flows) and we explain how to visualize twist coordinates in Section~\ref{sec:one-holed-tori} and prove Theorem~\ref{thm:polygonal-model-intro}. We also give explicit formulae for the holonomy of a hyperbolic structure on a one-holed torus with prescribed Fenchel--Nielsen coordinates (Section~\ref{sec:one-holed-torus-explicit-paramaetrization}).

Section~\ref{sec:pentagon-representations} treats pentagon representations, starting by recalling Le Fils' definition (Definition~\ref{defn:pentagon-representation}) and stating their main properties. We then go on with the description of our coordinates for pentagon representations (Section~\ref{sec:parametrization-pentagon-representations}) and prove Theorem~\ref{thm:pentagon-representations-intro} by showing that they satisfy a Wolpert-type formula (Section~\ref{sec:symplectic-structure-pentagon-representations}). Finally, we explain how to geometrize pentagon representations and prove Theorem~\ref{thm:geometrization-pentagon-intro} (Section~\ref{sec:geometrization-pentagon-representations}).

We study bow-tie representations in Section~\ref{sec:bow-tie-representations}. After describing our Fenchel--Nielsen-type coordinates for bow-tie representations (Section~\ref{sec:parametrization-bow-tie-representations}), we prove Theorem~\ref{thm:bow-tie-representations-intro}, namely that they also satisfy a Wolpert-type formula (Section~\ref{sec:symplectic-structure-bow-tie-representations}). We conclude this note by showing how to geometrize bow-tie representations using slit tori assemblies and prove Theorem~\ref{thm:geometrization-bow-tie-intro} (Section~\ref{sec:geometrization-bow-tie-representations}).

\subsection{Acknowledgments}
We express our gratitude to Thomas Le Fils for sharing enriching insights about pentagon representations and especially for explaining to us their characterization as fixed points of the hyperelliptic involution. Thank you to Nalini Anantharaman and Gustave Billon for fruitful conversations. GF would like to thank Nalini Anantharaman and the University of Strasbourg for their hospitality and support for this work. GF also acknowledges the support of INDAM–GNSAGA. AM thanks Samuele Mongodi and the University of Milano-Bicocca for their generous hospitality. AM was supported by the European Research Council (ERC) under the European Union’s Horizon 2020 research and innovation program (Grant agreement No.~101096550).
\section{Background}\label{sec:background}

\subsection{Character varieties}\label{sec:character-varieties}
We give a rapid recap on character varieties of surface groups and their fundamental attributes. We will mostly focus on $\psl$ character varieties. For further details and references, the reader may consult~\cite{maret-notes}.

\subsubsection{Definition}
Let $S$ be a closed and oriented topological surface of genus $g\geq 2$ and $G$ be a connected Lie group. Broadly speaking, the character variety of the pair $(S,G)$ is the deformation space of conjugacy classes of all representations $\pi_1S\to G$. Depending on the expected properties of that space, increasingly sophisticated definitions of character varieties are possible. We opt for the simplest possible one.
\begin{defn}\label{defn:character-varieties}
We call \emph{character variety} of the pair $(S,G)$ the topological quotient of the space of all representations $\pi_1S\to G$ (equipped with the compact-open topology) by the conjugation action of $G$. We denote it by $\chi(S,G)$.
\end{defn}

\subsubsection{Topological components}\label{sec:character-varieties-components}
The points of $\chi(S,G)$ may also be interpreted as equivalence classes of flat principal $G$-bundles over $S$. Using the theory of characteristic classes of $G$-bundles, we can define a conjugacy invariant for representations $\pi_1S\to G$ which materializes as a continuous map $\chi(S,G)\to H^2(S,\pi_1G)\cong \pi_1G$, where the isomorphism is given by the orientation of $S$. Equivalently, this invariant may be defined as an obstruction to lift a representation $\pi_1S\to G$ to a representation $\pi_1S\to \widetilde G$, where $\widetilde G$ is the universal cover of $G$ (see~\cite{goldman-components} for more details). 

In the case where $G=\psl$ is the group of orientation-preserving isometries of the upper half-plane, then this invariant is known as \emph{Euler number}. Since $\pi_1(\psl)\cong\Z$, the Euler number is an integer-valued invariant. The Euler number $\eu(\rho)$ of a representation $\rho\colon\pi_1S\to\psl$ can be computed as
\begin{equation}\label{eq:Euler-number-integral}
\eu(\rho)=\frac{1}{2\pi}\int_S f^\ast \omega_{\HH},
\end{equation}
where $f\colon \widetilde S\to \HH$ is any smooth $\rho$-equivariant map from the universal cover $\widetilde S$ of $S$ to the upper half-plane $\HH$, and $\omega_{\HH}$ is the standard hyperbolic area form\footnote{In particular, $\int_D\omega_\HH>0$ for anti-clockwise oriented disks $D$.} on $\HH$. It turns out that $\eu(\rho)$ is only constrained by the Milnor--Wood inequality which reads
\[
|\eu(\rho)|\leq 2g-2,
\]
where $g\geq 2$ is the genus of $S$. Goldman showed that the induced map 
\[
\eu_\ast\colon\pi_0\big(\chi(S,\psl)\big)\to\{2-2g,\ldots,2g-2\}
\]
is a bijection~\cite[Theorem~B]{goldman-components}. In other words, the character variety $\chi(S,\psl)$ has $4g-3$ connected components indexed by the Euler number. The two components with extremal Euler number $\pm(2g-2)$ are two copies of the \emph{Teichmüller space} of $S$---the space of conformal structures on $S$, up to isotopy---and consists of precisely all the conjugacy classes of discrete and faithful representations $\pi_1S\to \psl$~\cite[Corollary~C]{goldman-components}. All the representations $\pi_1S\to \psl$ with extremal Euler number $e=\pm (2g-2)$---also known as \emph{Fuchsian} representations---are holonomies of hyperbolic structures on $S$ (or on $S$ with opposite orientation) and therefore map every non-trivial element of $\pi_1S$ to a hyperbolic element of $\psl$. Bowditch famously asked whether Fuchsian representations are characterized by the property of mapping every \emph{simple} closed curve on $S$ to a hyperbolic element of $\psl$~\cite[Question~C]{bowditch}. Bowdtich's question admits a positive answer when $g=2$ by the work of Marché--Wolff, but remains completely open for larger genuses.

\begin{thm}[{\cite{marche-wolff, marche-wolff-2}}]\label{thm:bowditch-genus-2}
If $S$ is a closed and oriented surface of genus $2$ and $\rho\colon\pi_1S\to \psl$ sends every element of $\pi_1S$ representing a simple closed curve on $S$ to a hyperbolic element of $\psl$, then $\rho$ is a Fuchsian representation.   
\end{thm}

The two components of $\chi(S,\psl)$ with Euler number $e$ and $-e$ are isomoprhic to each other under the involution 
\begin{equation}\label{eq:involution-j}
\jmath\colon\chi(S,\psl)\to\chi(S,\psl)
\end{equation}
given by post-composition by the unique non-trivial outer automorphism of $\psl$. In concrete words, $\jmath$ is simply the conjugation by any matrix with determinant $-1$.

All the components with Euler number $e\neq 0$ are smooth manifolds of dimension $6g-6$. Hitchin actually showed that each of them is real-analytically diffeomorphic to a complex vector bundle over $\Sym^{2g-2-e}(S)$---the space of multi-sets of $2g-2-e$ elements in $S$~\cite[Theorem~10.8]{hitchin-self-duality}. For instance, the two components with Euler number $e=\pm(2g-3)$ (or co-Euler number~$1$) are diffeomorphic to complex vector bundles of complex rank $3g-4$ over $S$. 

\subsubsection{Symplectic nature}\label{sec:character-variety-symplectic-nature}
If $G$ is a \emph{quadrable} Lie group,\footnote{That is, if the Lie algebra of $G$ admits a symmetric and non-degenerate bi-linear form which is invariant under the adjoint action of $G$.} then the smooth locus of $\chi(S,G)$ carries a natural symplectic structure $\omega_\mathcal{G}$ named after Goldman~\cite{goldman-symplectic}. The quadrability assumption is satisfied by all reductive Lie groups, like $\psl$,\footnote{On to the two components of $\chi(S,\psl)$ with extremal Euler number---those isomorphic to two copies of the Teichmüller space of $S$---\(\omega_{\mathcal G}\) coincides with the Weil--Petersson symplectic form.} but not only (see e.g.~\cite[Section~1.1.3]{maret-notes}). We will mostly use this symplectic structure to associate to any real-valued smooth function $h$ defined on a smooth open subset of $\chi(S,G)$ its Hamiltonian flow. Recall that the \emph{Hamiltonian flow} of a function $h$ on a symplectic manifold $(M,\omega)$ is the flow generated by the \emph{Hamiltonian vector field} $X_h$ defined by the relation 
\begin{equation}\label{eq:Hamiltonian-vf}
    \omega(X_h,\cdot)=dh.
\end{equation}
The Goldman symplectic form on $\chi(S,G)$ has an associated Liouville measure called \emph{Goldman measure}. It also defines a Poisson bracket for smooth functions given by
\begin{equation}\label{eq:Poisson-bracket}
\{f,g\}=\omega_\mathcal{G}(X_f,X_g)=df(X_g)=-dg(X_f).
\end{equation}

When $G=\psl$, the involution $\jmath$ from~\eqref{eq:involution-j} preserves the Goldman symplectic form, therefore defining a symplectomorphism between the two components of Euler number $e$ and $-e$.

\subsubsection{Mapping class group dynamics}\label{sec:dynamics}
The \emph{mapping class group} of $S$ is the group of isotopy classes of orientation-preserving homeomorphisms of $S$. We will denote it by $\Mod(S)$. The Dehn--Nielsen--Baer theorem (see e.g.~\cite[Theorem~8.1]{mcg-primer}) identifies it with an index-2 subgroup of $\Out(\pi_1S)$---the group of outer automorphisms of $\pi_1S$. The mapping class group of $S$ therefore acts on $\chi(S,G)$ by pre-composition. This action is symplectic---it preserves the Goldman symplectic form, hence the Goldman measure too. On the other hand, the elements of $\Out(\pi_1S)$ induced by orientation-reversing homeomorphisms of $S$ act on $\chi(S,G)$ by anti-symplectic diffeomorphisms. 

In the case where $G=\psl$, the mapping class group of $S$ preserves each connect component of $\chi(S,\psl)$. The action of $\Mod(S)$ is known to be properly discontinuous on the two components with extremal Euler number, while it was conjectured by Goldman to be ergodic (with respect to the Goldman measure) on all other components~\cite[Conjecture~2.6]{goldman-conjectures}. All ergodic components of $\chi(S,\psl)$ were classified by Marché--Wolff in the case $g=2$, but Goldman's conjecture remains open in larger genuses. More precisely, Marché--Wolff obtained the following result as a consequence of Theorem~\ref{thm:bowditch-genus-2}.
\begin{thm}[{\cite{marche-wolff, marche-wolff-2}}]\label{thm:goldman-ergocicity-genus-2}
Let $S$ be a closed and oriented surface of genus $2$. The action of $\Mod(S)$ on the two components of $\chi(S,\psl)$ with Euler number $\pm 1$ is ergodic. Moreover, the action on the component of Euler number $0$ has two disjoint ergodic components. 
\end{thm}

For further considerations about mapping class groups, we refer the reader to~\cite{mcg-primer}. We also refer the reader to Appendix~\ref{apx:generators}, where an explicit generating family for the mapping class group of a closed surface of genus two is described.

\subsection{Branched hyperbolic surfaces}\label{sec:branched-hyperbolic-structures}
This section serves as a brief introduction to branched hyperbolic structures on surfaces. We recall some important results and conjectures.

\subsubsection{Definition}
Let $S$ be a closed and oriented surface of genus $g\geq 2$. A (standard) \emph{hyperbolic structure} on $S$ is given by a maximal atlas of charts $\varphi_i\colon U_i\to \HH$ with values in the hyperbolic plane and transition maps in $\psl$. The charts are by definition homeomorphisms onto their images and preserve orientation. In order to define branched hyperbolic structures, we relax the local homeomorphism condition.

\begin{defn}\label{defn:branched-hyperbolic-structures}
A \emph{branched hyperbolic chart} around $p\in S$ is a continuous, orientation-preserving map $\varphi\colon U\to\HH$ from an open neighbourhood $U$ of $p$ to the Poincaré disk model of $\HH$ such that $\varphi(p)=0$ and such that there exist an integer $k\geq 1$ and a local coordinate $z$ centred at $p$ in $U$ for which $\varphi$ behaves as $z\mapsto z^k$.\footnote{More precisely, we require a \emph{local coordinate} $z\colon U\to \C$---an orientation-preserving homeomorphism onto its image---satisfying $z(p)=0$ and $\varphi\circ z^{-1}(w)=w^k$ for every $w\in z(U)$.}
When $k=1$, the chart is a local homeomorphism; when $k\geq 2$, the point $p$ is a \emph{branch point} with a \emph{cone singularity} of angle $2\pi k$. A \emph{branched hyperbolic structure} on $S$ is a maximal atlas $\{(U_i,\varphi_i)\}$ of branched hyperbolic charts with transition maps $t_{ij}\colon \varphi_j(U_i\cap U_j)\to \varphi_i(U_i\cap U_j)$ that are restrictions of elements of $\psl$ and satisify $t_{ij}\circ \varphi_j = \varphi_i$ on $U_i\cap U_j$.
\end{defn}

When $S$ is equipped with a branched hyperbolic structure, we say that $S$ is a \emph{branched hyperbolic surface}. The set $\mathcal{P}$ of branch points associated to a branched hyperbolic surface $S$ is a discrete, hence finite, subset of $S$. Away from the branch points, the structure is a standard hyperbolic structure.  

An equivalent definition consists of saying that a branched hyperbolic surface is a compact surface obtained by gluing finitely many hyperbolic triangles along isometries between their edges such that the total angle around each vertex is an integer multiple of $2\pi$. 

Branched hyperbolic structures may also be defined as the data of a Riemannian metric $h$ of constant curvature $-1$ on $S\setminus \mathcal{P}$ whose metric completion at each $p\in \mathcal{P}$ is a cone singularity of angle $2\pi k_p$ for some $k_p\geq 2$. Those metrics are special examples of \emph{hyperbolic cone metrics} on $S\setminus \mathcal{P}$. The Gauss--Bonnet theorem expresses the volume of $h$ as
\begin{equation}\label{eq:volume-branched-hyperbolic-metric}
\Vol(h)=4\pi(g-1) + 2\pi|\mathcal{P}| - 2\pi\sum_{p\in\mathcal{P}} k_p.
\end{equation}
Since $\Vol(h)$ is positive, we deduce that $\sum_{p\in\mathcal{P}} (k_p-1) \leq 2g-3$. In other words, the maximal angle excess of $h$ is $2\pi (2g-3)$ which is distributed (not necessarily evenly, but always by units of $2\pi$) among the branch points.

For the upcoming discussion and following~\cite{goldman-conjectures}, it will be more convenient to allow $\mathcal{P}$ to be a multi-set of points on $S$ and assume that $k_p=2$ for every $p\in\mathcal{P}$. The geometric cone angle at $p\in\mathcal{P}$ is now given by $2\pi (m_p+1)$, where $m_p$ is the multiplicity of $p$ in~$\mathcal{P}$. 

Recall that $\Sym^n(S)$ denotes the space of multi-sets of $n$ points on $S$. Any choice of complex structure on $S$ turns it into a smooth manifold of real dimension $2n$. Given an integer $0\leq n\leq 2g-3$, we will denote by 
\[
\Hyp_n(S)
\]
the space of all branched hyperbolic structures on $S$, up to isotopy, with branch point set in $\Sym^n(S)$. We emphasize that isotopies are allowed to move branch points. The mapping class group of $S$ acts on $\Hyp_n(S)$. In the special case where $n=0$ and $\mathcal{P}=\emptyset$, we recover the Teichmüller space of~$S$. Note that we assumed $n\leq 2g-3$, as there are no branched hyperbolic structures on $S$ with $n> 2g-3$ branch points by the Gauss--Bonnet theorem. We also introduce the space of all branched hyperbolic structures on $S$, up to isotopy, which we denote by
\[
\Hyp(S)=\bigsqcup_{n=0}^{2g-3}\Hyp_n (S).
\]

\subsubsection{Uniformization}
Since Möbius transformations are holomorphic, a branched hyperbolic structure on $S$ with branch point set $\mathcal{P}$ also defines a \emph{complex structure} on $S\setminus \mathcal{P}$---a maximal atlas of charts to $\C$ with biholomorphic change of coordinates. Since the local model at branch points is given by $z\mapsto z^k$, the complex structure extends over branch points.\footnote{For instance, if $(\varphi_i, U_i)$ and $(\varphi_j,U_j)$ are two branch hyperbolic charts around a branch point $p$, then by definition we can find local coordinates $z_i\colon U_i\to \C$ and $z_j\colon U_j\to \C$ with $z_i(p)=z_j(p)=0$ and such that $\varphi_i\circ z_i^{-1}(w)=w^k$ and $\varphi_j\circ z_j^{-1}(w)=w^k$. These local coordinates will serve as complex charts around the branch points. If $t_{ij}\colon\varphi_j(U_i\cap U_j)\to \varphi_i(U_i\cap U_j)$ denotes the transition map between the two branch hyperbolic charts, then $(z_i\circ z_j^{-1}(w))^k=t_{ij}(w^k)$ for every $w\in z_j(U_i\cap U_j)$. Since $t_{ij}$ is a non-constant holomorphic map with $t_{ij}(0)=0$, its Taylor expansion at $0$ reads $t_{ij}(w)=w^nu(w)$ for some holomorphic function $u(w)$ with $u(0)\neq 0$ and some positive integer $n$. This means that $z_j\circ z_i^{-1}(w)=w^n u(w^k)^{1/k}$ for the right $k^{\textrm{th}}$-root of the function $u(w^k)$ (they differ by a $k^{\textrm{th}}$-root of unity), and is therefore holomorphic.} Using the equivalence between complex structures and conformal classes, we get a map
\begin{equation}\label{eq:complex-structure-of-branched-hyperbolic-structure}
\Hyp (S)\to \Teich (S),
\end{equation}
where $\Teich(S)$ denotes the Teichmüller space of $S$, seen here as the space of conformal structures on $S$, up to isotopy.

McOwen and Troyanov showed that every choice of conformal structure on the surface $S$ and marked point set $\mathcal{P}$ leads to a unique hyperbolic cone metric with prescribed angle at each point of $\mathcal{P}$, provided that the angles satisfy a Gauss--Bonnet condition.
\begin{thm}[{\cite{mcowen, troyanov}}]\label{thm:troyanov-uniformization}
For every conformal structure on the surface $S$ and every multi-set of marked points $\mathcal{P}$ on $S$ satisfying
\begin{equation}\label{eq:Gauss-Bonnet-inequality}
|\mathcal{P}| \leq 2g-3,
\end{equation}
there exists a unique branched hyperbolic structure on $S$ which induces the same conformal structure and with a cone singularity of angle $2\pi (m_p+1)$ at each $p\in \mathcal{P}$, where $m_p$ is the multiplicity of $p$ in $\mathcal{P}$.
\end{thm}

Theorem~\ref{thm:troyanov-uniformization} gives a new perspective on the space $\Hyp_n(S)$. Consider the fiber bundle over $\Teich(S)$ with total space $\mathcal{S}^n(S)$ consisting of all isotopy classes of pairs of a conformal structure on $S$ together with a multi-set $\mathcal{P}\in\Sym^n(S)$. This bundle can be given the structure of a holomorphic fiber bundle over $\Teich(S)$ which is topologically equivalent to the trivial bundle $\Teich(S)\times \Sym^n(S)$ because $\Teich(S)$ is contractible. By also recording the branch point set, we can upgrade the map~\eqref{eq:complex-structure-of-branched-hyperbolic-structure} to
\[
\Hyp_{n}(S)\stackrel{\cong}{\longrightarrow}\mathcal{S}^n(S),
\]
which is invertible by Theorem~\ref{thm:troyanov-uniformization}. This identification turns $\Hyp_{n}(S)$ into a smooth manifold of real dimension $6g-6+2n$.\footnote{It is also possible to topologize $\Hyp_{n}(S)$ with the induced $\mathcal{C}^\infty$ topology from the space of developing maps. Nguyen proved that the two topologies coincide~\cite[Theorem~3.2.8 and the discussion beforehand]{nguyen}.}

\subsubsection{Holonomy}\label{sec:holonomy}
A branched hyperbolic structure on $S$ with branch points set $\mathcal{P}$ defines a standard hyperbolic structure on $S\setminus \mathcal{P}$ and a \emph{developing map} $\dev\colon\widetilde{S\setminus \mathcal{P}}\to \HH$ which is unique up to pre-composition by isotopies of $S\setminus \mathcal{P}$ lifted to the universal cover and post-composition by $\psl$. It is a local homeomorphism that preserves orientation and is equivariant with respect to a unique representation $\rho\colon\pi_1(S\setminus \mathcal{P})\to \psl$ known as \emph{holonomy}. Since the holonomy around each branch point is trivial, $\rho$ factorizes through the projection $\pi_1(S\setminus \mathcal{P})\to \pi_1S$ to define a representation $\rho'\colon\pi_1 S\to \psl$. For the same reason, $\dev$ factorizes through the natural map $\widetilde{S\setminus \mathcal{P}}\to\widetilde{S}$, producing a map $\dev'\colon\widetilde{S}\to\HH$ which is $\rho'$-equivariant. The difference between $\dev$ and $\dev'$ is that $\dev'$ fails to be a local homeomorphism at the pre-images of branch points, but is locally a branched covering. Overall, taking the holonomy of a branched hyperbolic structure on $S$ defines a $\Mod(S)$-equivariant map
\begin{equation}\label{eq:holonomy-map}
\hol\colon\Hyp(S) \to \chi(S,\psl).
\end{equation}

If $[\rho]\in \chi(S,\psl)$ is the holonomy of a branched hyperbolic structure on $S$ belonging to $\Hyp_n(S)\subset \Hyp(S)$, then from~\eqref{eq:Euler-number-integral} and~\eqref{eq:volume-branched-hyperbolic-metric} we obtain
\[
\eu(\rho)=\frac{1}{2\pi}\int_S \dev^\ast \omega_\HH=2g-2-n.
\]
This means that $\hol(\Hyp_n(S))$ is contained in the component of $\chi(S,\psl)$ of Euler number $2g-2-n$. Since $n\leq 2g-3$, the Euler number of the holonomy of a branched hyperbolic structure on $S$ is always contained in $\{1,\ldots, 2g-2\}$ and is never $0$.\footnote{If we flip the orientation of $S$, then the same developing map $\dev$ is now orientation-reversing. It is henceforth equivariant with respect to a conjugate of $\rho$ by a determinant $-1$ matrix which changes the sign of the Euler number. So, we may think of $\{-1,\ldots, -(2g-2)\}$ as Euler numbers of branch hyperbolic structures on $S$ with opposite orientation.} 

Going back to McOwen--Troyanov uniformization theorem (Theorem~\ref{thm:troyanov-uniformization}), recall that the fibers of $\Hyp_n(S)\to\Teich(S)$ define a trivial foliation of $\Hyp_n(S)$ where each fiber is a copy of $\Sym^n(S)$. We will call it the \emph{isoconformal foliation} of $\Hyp_n(S)$. Note that the restriction of $\hol$ to each isoconformal leaf is injective. Indeed, given two different multisets of branch points $\mathcal{P}$ and $\mathcal{P}'$ on the same Riemann surface $S$, the corresponding branched hyperbolic structures are distinct and determined by their developing maps $\dev\colon\widetilde{S}\to\HH$ and $\dev'\colon\widetilde{S}\to\HH$. We may assume that $\dev$ and $\dev'$ are the unique holomorphic maps in their isotopy class, up to post-composition by a Möbius transformation. The developing maps are equivariant with respect to two holonomy representations $\rho\colon\pi_1S\to\psl$ and $\rho'\colon\pi_1S\to\psl$. If $\rho$ and $\rho'$ were conjugate, then $\dev$ and $\dev'$ would differ by post-composition by a Möbius transformation. In particular, they would have the same set of branch points, contradiction. In conclusion, every choice of conformal structure on $S$ defines a slice
\begin{equation}\label{eq:slice-Sym^n(S)}
    \Sym^n(S)\longhookrightarrow \chi(S,\psl)
\end{equation}
that lies in the component of Euler number $2g-2-n$. As we will explain at the end of Section~\ref{sec:Hodge}, all these slices are actually embedded submanifolds. 

Another important foliation of \(\Hyp_n(S)\) is its \emph{isomonodromic foliation}, defined by the fibers of the holonomy map $\Hyp_n(S) \to \chi(S,\psl)$. Each isomonodromic leaf is a \(2n\)-dimensional submanifold of \(\Hyp_n(S)\) (see the discussion before Question~\ref{question:isomonodromic-deformations}). 
The global topology of the isomonodromic foliation is currently unclear to us (see also Question~\ref{question:isomonodromic-deformations}); for instance, is it everywhere transverse to the isoconformal foliation?
Nevertheless, one of the main goals of this paper is to describe explicit isomonodromic deformations for a large class of holonomy representations (see Corollaries~\ref{cor:isomonodromic-deformations-bow-tie} and~\ref{cor:isomonodromic-deformations-pentagon}).

\subsubsection{Geometrization conjecture}
We will denote the restriction of $\hol$ to $\Hyp_n(S)$ by $\hol_e$, where $e=2g-2-n$. It takes image in the component of $\chi(S,\psl)$ of Euler number~$e$. Goldman showed that $\hol_{2g-2}$ is surjective~\cite[Corollary~C]{goldman-components}. In other words, every representation with Euler number $2g-2$ is the holonomy of a hyperbolic structure on $S$. When $e\neq 2g-2$, it is an open problem to describe the image of $\hol_e$ (see e.g.~\cite[Question~3.1.14]{tholozan-hdr} or~\cite[Questions~4 \&~5]{mathews-genus-2}). The research in the field is guided by the following conjecture of Goldman known as the \emph{geometrization conjecture}.
\begin{conj}[{\cite[Conjecture~3.9]{goldman-conjectures}}]\label{conj:geoemtrization-conjecture}
If $e=2g-3$, then $\hol_e$ is surjective; if $1\leq e\leq 2g-4$, then every representation with dense image lies in the image of $\hol_e$.
\end{conj}
\begin{defn}\label{defn:geometrizable-representations}
A representation $\rho\colon\pi_1S\to \psl$ is \emph{geometrizable} if its conjugacy class or its image by the involution $\jmath$ from~\eqref{eq:involution-j} belongs to the image of the holonomy map $\hol$.
\end{defn}

There are geometrizable representations in every component of Euler number $e\neq 0$. For instance, the union of all the slices described in~\eqref{eq:slice-Sym^n(S)} coming from the McOwen--Troyanov uniformization theorem is made of geometrizable representations. Also note that there exist geometrizable representations with discrete image and non-extremal Euler number such as every pentagon representation coming from a right-angled pentagon (Examples~\ref{ex:pentagon-representation} and~\ref{ex:pentagon-representation-P5}), as well as the bow-tie representation from Example~\ref{ex:discrete-bow-tie} for instance. So, the conclusion of Conjecture~\ref{conj:geoemtrization-conjecture} is not a complete characterization of the image of $\hol_e$. On the other hand, some representations with Euler number $1\leq \pm e\leq 2g-4$ are not geometrizable. The following examples were described by Tan.

\begin{ex}[{\cite[Example~8]{tan}}]\label{ex:non-geometrizable-representations}
Let $S$ and $S'$ be two closed and oriented surfaces. Let $j\colon\pi_1 S'\to \psl$ be a Fuchsian representation and $f\colon S\to S'$ be any continuous map preserving orientation that is not homotopic to a branched cover. For instance, Tan takes $f$ to be a handle-crushing map. Then $\rho=j\circ f_\ast\colon\pi_1 S\to \psl$ is not geometrizable. Indeed, assume for the sake of contradiction that $\rho$ was geometrizable. There would then exist an orientation-preserving developing map $\dev\colon \widetilde{S}\to \HH$ which descends to a map $\overline{\dev}\colon S\to \HH/\rho(\pi_1 S)$. Post-composing $\overline{\dev}$ with the projection $\HH/\rho(\pi_1 S)\to \HH/j(\pi_1 S')\cong S'$, would lead to a branch covering map $F\colon S\to S'$. Note that $F$ and $f$ would then induce the same morphism at the level of fundamental groups, implying that they are homotopic, a contradiction.
\end{ex}

Tan also observed that the non-geometrizable representations from Example~\ref{ex:non-geometrizable-representations} can be approximated by geometrizable ones. This led to the following question.
\begin{question}[{\cite{tan}}]\label{question:tan}
Is the set of conjugacy classes of geometrizable representations dense in every component of $\chi(S,\psl)$ of Euler number $e\in\{\pm 1,\dots, \pm(2g-3)\}$? 
\end{question}
The Ehresmann--Thurston principle implies that the set of geometrizable representations is an open subset of $\chi(S,\psl)$. It is also invariant under the action of the mapping class group of $S$. So, if Goldman's ergodicity conjecture was true, then the answer to Question~\ref{question:tan} would be positive. In particular, the answer is yes for surfaces of genus $2$ by the work of Marché--Wolff (Theorem~\ref{thm:goldman-ergocicity-genus-2}). Similarly, since the set of faithful representations is dense in every component of $\chi(S,\psl)$ by a result of DeBlois--Kent~\cite{deblois-kent} and injective representations with non-extremal Euler number $e\neq \pm(2g-2)$ have dense image, it follows that a positive answer to Question~\ref{question:tan} is also implied by Goldman's geometrization conjecture (Conjecture~\ref{conj:geoemtrization-conjecture}).

Partial progress towards Question~\ref{question:tan} and Conjecture~\ref{conj:geoemtrization-conjecture} was achieved by Mathews. For instance, he proved that if the genus of $S$ is $2$, then every representation with Euler number~$\pm 1$ that sends a separating simple closed curve on $S$ to a non-hyperbolic element of $\psl$ is geometrizable~\cite[Theorem~2]{mathews-genus-2}. He also proved that for any surface $S$ of genus $g\geq 2$, almost every representation with Euler number $\pm(2g-3)$ that send a non-separating simple closed curve on $S$ to an elliptic element of $S$ is geometrizable~\cite[Theorem~3]{mathews-genus-2}. The case of purely hyperbolic representations---representations whose image consists of hyperbolic elements of $\psl$ and the identity---was developed in~\cite{gianluca}.

\subsubsection{Hodge theoretic perspective}\label{sec:Hodge}
We already explained that the data of a branched hyperbolic structure on $S$ naturally turns $S$ into a Riemann surface. With respect to this complex structure, the developing map $\widetilde S\to \HH$ becomes a $\rho$-equivariant holomorphic map.

Conversely, given a representation $\rho\colon\pi_1S\to \psl$, if there exist a complex structure on $S$ and a $\rho$-equivariant holomorphic map $f\colon\widetilde S\to \HH$ which is not constant, then $f^\ast g_\HH$ is branched hyperbolic metric on $S$. Here $g_\HH$ denotes the standard hyperbolic metric on $\HH$. Indeed, since $f$ is not constant, it has isolated critical points. Let $\mathcal{P}$ be the finite set of all projections of the critical points onto $S$. Away from its critical points, $f$ is a local biholomorphism and so $f^\ast g_\HH$ defines a hyperbolic metric on $S\setminus \mathcal{P}$. Near a critical point $p$, in a local coordinate $z$, $f^\ast g_\HH\sim |z|^{2(k_p-1)}|dz|^2$, where $k_p\geq 2$ is the branching order of $f$ at $p$. This corresponds to a cone singularity of angle $2\pi k_p$.

In conclusion, given a representation $\rho\colon\pi_1 S\to \psl$, the following are equivalent:
\begin{enumerate}
    \item $\rho$ is geometrizable (Definition~\ref{defn:geometrizable-representations}).
    \item There exist a complex structure on $S$ and a $\rho$-equivariant holomomorphic or anti-holomorphic map $\widetilde S\to \HH$ which is non-constant.
\end{enumerate}

Given a complex structure on $S$, the non-abelian Hodge correspondence developed by Hitchin~\cite{hitchin-self-duality}, Donaldson~\cite{donaldson}, Corlette~\cite{corlette}, and Simpson~\cite{simpson} gives a real-analytic isomorphism between the subspace of $\chi(S,\psl)$ made of conjugacy classes of reductive representations and a certain moduli space of polystable Higgs bundles. Roughly speaking, a \emph{Higgs bundle} is a rank-2 holomorphic vector bundle $E$ on $S$ with a holomorphic 1-form on $S$ with values in $\End(E)$, which is known as the \emph{Higgs field}. Scaling the Higgs field by a non-zero complex number defines an action of $\C^\times$ on the moduli space of polystable Higgs bundles. If $\rho\colon\pi_1 S\to \psl$ is a reductive representation, then the following condition is equivalent to $\rho$ being geometrizable:
\begin{enumerate}\setcounter{enumi}{2}
    \item There exists a complex structure on $S$, such that the equivalence class of Higgs bundles on $S$ corresponding to $[\rho]$ under the non-abelian Hodge correspondence is a fixed point of the $\C^\times$ action. Equivalently, the Higgs bundles come from a \emph{complex variation of Hodge structures} (see~\cite{simpson-cvhs} for more details).
\end{enumerate}

For a fixed complex structure on $S$, Hitchin fully described the locus of fixed points for the $\C^\times$ action, which also corresponds to the minima of some energy functional~\cite{hitchin-self-duality}. It turns out that the locus of fixed points in the moduli space of Higgs bundles corresponding to the components of $\chi(S,\psl)$ with Euler number $e=2g-2-n$ is precisely the slice $\Sym^n(S)$ described in~\eqref{eq:slice-Sym^n(S)}. We emphasize that different choices of conformal structures on $S$ lead to different loci of fixed points and different slices in the character variety. In particular, it shows that each slice is an embedded submanifold of $\chi(S,\psl)$ isomorphic to $\Sym^n(S)$, hence of real dimension $2n$. Hitchin also explained that following the flow lines of the $\C^\times$ action defines a homotopy equivalence between the slice $\Sym^n(S)$ and the moduli component.
\section{Revisiting one-holed tori}\label{sec:one-holed-tori}
\subsection{Overview}
Hyperbolic one-holed tori are the elementary pieces of our construction. We recall some classical results (Section~\ref{sec:Fenchel-Nielsen-coordinates-one-holed-torus}) and provide a unified method for encoding their geometry into a polygon (Section~\ref{sec:polygonal-model-one-holed-tori}). This model helps visualizing Fenchel--Nielsen coordinates, in particular the twist coordinate, and the Hamiltonian twist flow of the associated length function. Finally, we give explicit formulae for the holonomy of a hyperbolic one-holed torus with prescribed Fenchel--Nielsen coordinates (Section~\ref{sec:one-holed-torus-explicit-paramaetrization}).

\subsection{Fenchel--Nielsen coordinates}\label{sec:Fenchel-Nielsen-coordinates-one-holed-torus}
\begin{figure}[h]
\centering
\begin{tikzpicture}[scale=1]
\node[anchor = south west, inner sep=0mm] at (0,0) {\includegraphics[width=3.8cm]{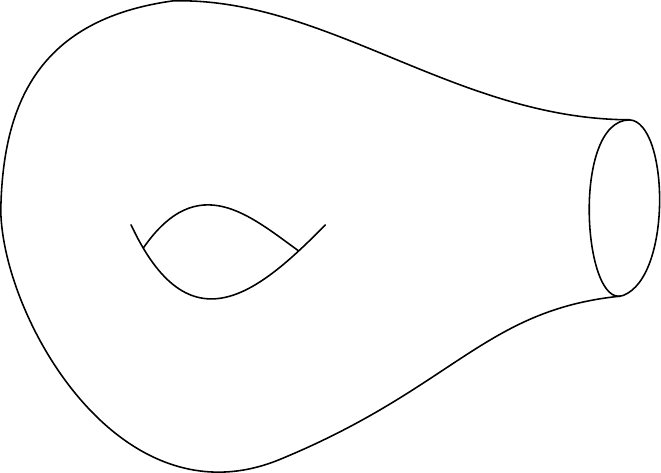}};
\node[anchor = south west, inner sep=0mm] at (4.5,0.2) {\includegraphics[width=4.2cm]{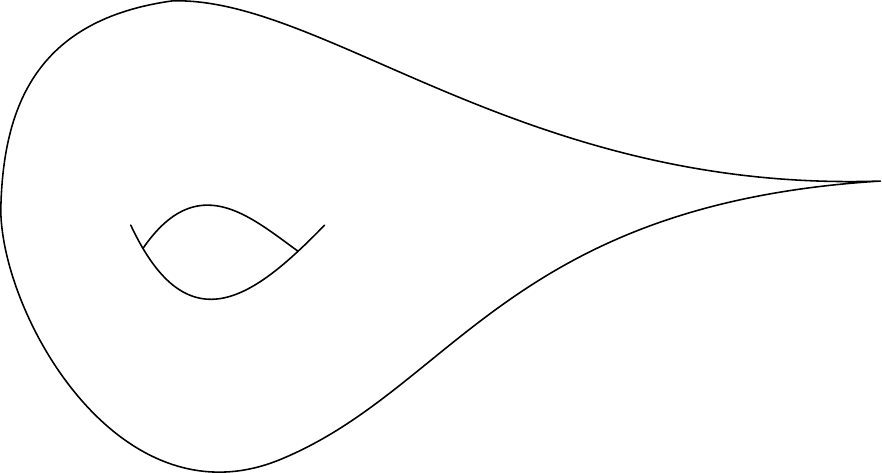}};
\node[anchor = south west, inner sep=0mm] at (9,0) {\includegraphics[width=4cm]{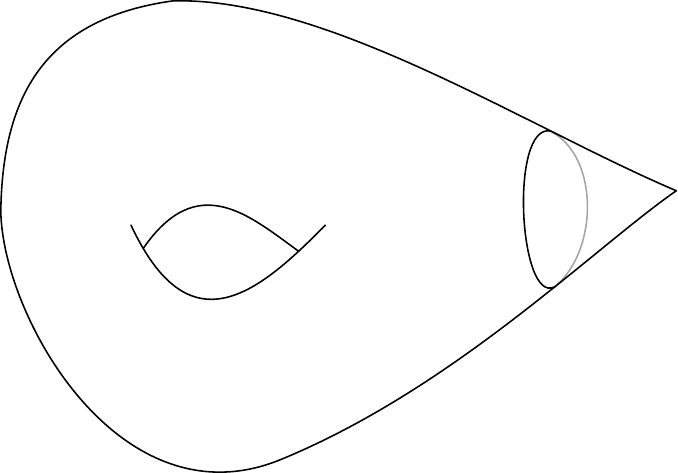}};
\end{tikzpicture}
\caption{From the left: a hyperbolic one-holed torus with a geodesic boundary, a cusp, and a cone point.}
\label{fig:hyperbolic-one-holed-tori}
\end{figure}

Let $\Sigma$ be an oriented topological torus with one puncture. The fundamental group of $\Sigma$ can be presented as
\begin{equation}\label{eq:presentation-pi_1-one-holed-torus}
\pi_1\Sigma =\langle a,b,c\, : \, [a,b]=c\rangle,
\end{equation}
where $c$ is the homotopy class of a loop circling around the puncture of $\Sigma$. Note that $\pi_1\Sigma$ is isomorphic to the free group on two generators. This means that any representation $\rho\colon\pi_1\Sigma\to\psl$ can be lifted to four representations $\pi_1\Sigma\to\SL_2\R$. The image of~$c$ is the same for every lift, which means that we can define $\Tr\rho(c)$ to the be the trace of the image of~$c$ for any lift of $\rho$. We are interested in the three regimes $\Tr\rho(c)<-2$, $\Tr\rho(c)=-2$, and $-2<\Tr\rho(c)<2$ because they correspond to the three types of hyperbolic structures illustrated on Figure~\ref{fig:hyperbolic-one-holed-tori} that~$\Sigma$ can be endowed with~\cite[Section~3]{goldman-one-holed-torus}.
\begin{itemize}
    \item[(H)] $\Tr\rho(c)<-2$: $\rho$ is the holonomy of a hyperbolic structure on $\Sigma$ with geodesic boundary of length $2\arccosh(-t/2)$, where $t=\Tr\rho(c)$.
    \item[(P)] $\Tr\rho(c)=-2$: $\rho$ is the holonomy of a hyperbolic structure on $\Sigma$ with a cusp at the puncture.
    \item[(E)] $-2<\Tr\rho(c)<2$: $\rho$ is the holonomy of a hyperbolic structure on $\Sigma$ with a cone singularity of angle $2\arccos(-t/2)$ at the puncture, where $t=\Tr\rho(c)$.
\end{itemize}

Given $t<2$, we denote by $\chi_t(\Sigma,\psl)$ the space of conjugacy classes of representations $\rho\colon\pi_1\Sigma\to\psl$ for which $\Tr\rho(c)=t$. We call it the \emph{$t$-relative character variety} of~$\Sigma$. Each relative character variety $\chi_t(\Sigma,\psl)$ is a smooth manifold diffeomorphic to two disjoint copies of $\R^2$. The count of components follows from~\cite{goldman-one-holed-torus}. The two connected components of $\chi_t(\Sigma,\psl)$ are images of each other by conjugation by a matrix of determinant $-1$. Each $t$-relative character variety can be equipped with Goldman's symplectic form to be turned into a symplectic manifold~\cite{goldman-symplectic}. The symplectic structure can be described explicitly using trace coordinates as in~\cite[Section~2.4]{goldman-one-holed-torus} or, alternatively, using Fenchel--Nielsen coordinates which we now describe. We will use the following simple, but important observation.

\begin{fact}[{e.g.~\cite[Lemma~3.4.5]{goldman-one-holed-torus}}]\label{fact:tr<2-implies-crrossing-axes}
If $\rho\colon\pi_1\Sigma\to\psl$ is a representation satisfying $\Tr \rho(c)<2$, then $\rho(a)$ and $\rho(b)$ are hyperbolic elements of $\psl$ whose axes intersect inside the hyperbolic plane. Here, $a,b,c$ refer to the generators of $\pi_1\Sigma$ from~\eqref{eq:presentation-pi_1-one-holed-torus}.
\end{fact}

Fenchel--Nielsen coordinates on $\chi_t(\Sigma,\psl)$ for $t<2$ consist of a length coordinate and a twist coordinate. Traditionally, the definition of Fenchel--Nielsen coordinates depends on a choice of pants decomposition of the surface (each pants curve gives a length coordinate) and a choice of transverse seams to define twist coordinates (see e.g.~\cite[Section~10.6.1]{mcg-primer}). The algebraic alternative consists of choosing a geometric presentation of $\pi_1\Sigma$ like~\eqref{eq:presentation-pi_1-one-holed-torus}. The simple closed curve represented by $a$ determines the pants decomposition of $\Sigma$ and the curve represented by $b$ the transverse seam.

The length coordinate is the simplest to define.
  
\begin{defn}\label{defn:length-coordinate-one-holed-torus}
The \emph{length coordinate} is defined as the function
\begin{align*}
    \ell_a\colon \chi_t(\Sigma,\psl)&\to \R_{>0}\\
    [\rho] &\mapsto 2\arccosh\big(\vert\Tr\rho(a)\vert/2\big),
\end{align*}
where $a$ is the generator of $\pi_1\Sigma$ from~\eqref{eq:presentation-pi_1-one-holed-torus}. 
\end{defn}
Note that by Fact~\ref{fact:tr<2-implies-crrossing-axes}, $\rho(a)$ is hyperbolic and thus $\vert\Tr\rho(a)\vert>2$, which makes $\ell_a$ well defined. Equivalently, we can define $\ell_a([\rho])$ as the length of the geodesic representative of $a$ for the marked hyperbolic structure on $\Sigma$ corresponding to~$[\rho]$. 

\begin{figure}[h]
\centering
\begin{tikzpicture}
\node[anchor = south west, inner sep=0mm] at (0,0) {\includegraphics[width=6.5cm]{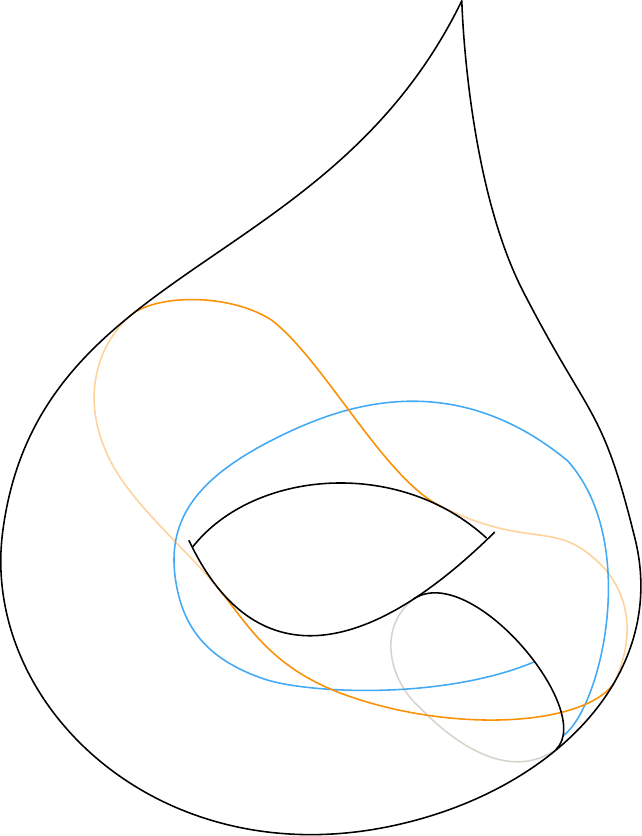}};

\node at (5.7,.6) {$\alpha$};
\node[apricot] at (2.5, 5) {$\beta$};

\node[anchor = south west, inner sep=0mm] at (7,0.2) {\includegraphics[width=6.5cm]{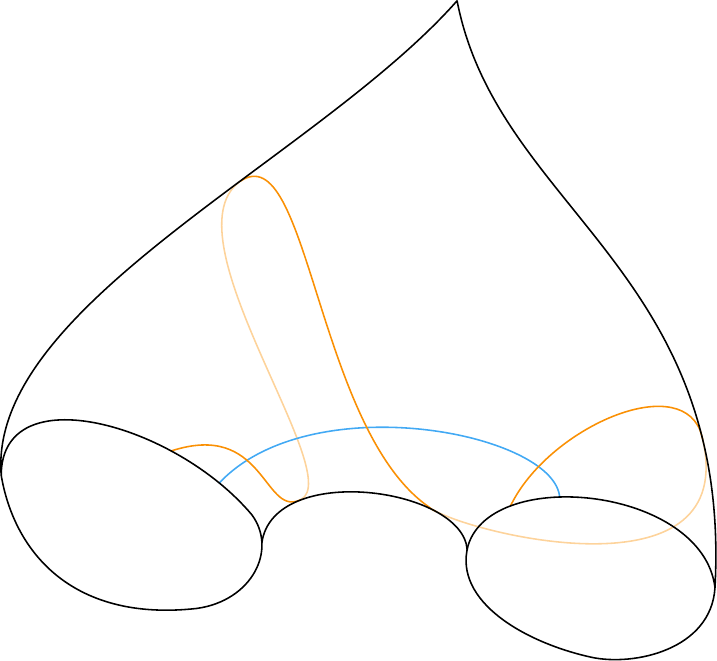}};

\node[sky] at (11,2.5) {$\gamma$};
\node at (12.5, -.1) {$\alpha_2$};
\node at (8, .4) {$\alpha_1$};
\end{tikzpicture}
\caption{Before and after cutting a hyperbolic one-holed torus along the geodesic $\alpha$.}
\label{fig:one-holed-torus-cut-and-uncut}
\end{figure}

The twist coordinate of $[\rho]\in\chi_t(\Sigma,\psl)$ is defined by first choosing a hyperbolic one-holed torus $\Upsilon$ with a marking $f\colon\Sigma\to \Upsilon$ whose holonomy is $[\rho]$. Such a choice is unique up to isotopy. We will denote by $\alpha$ and $\beta$ the geodesic representatives of $f_\ast a$ and $f_\ast b$ on $\Upsilon$. Cutting $\Upsilon$ along $\alpha$ produces a hyperbolic pair of pants (see Figure~\ref{fig:one-holed-torus-cut-and-uncut}) with two geodesic boundary curves $\alpha_1$ and $\alpha_2$ of length $\ell_a([\rho])$ (the same length as $\alpha)$, as well as a third geodesic boundary curve, a cusp, or a cone point depending on the value of $t<2$. There is a shortest geodesic arc $\gamma$ on the pair of pants joining $\alpha_1$ and $\alpha_2$. Back on $\Upsilon$, the geodesic $\beta$ is homotopic to the concatenation of $\gamma$ with an unique arc of $\alpha$ (that may go several times around $\alpha$). The signed length of that arc is denoted by $\tau$. The signed of $\tau$ is determined by the orientation of $\alpha$. 

\begin{defn}\label{defn:twist-coordinate-one-holed-torus}
The \emph{twist coordinate} is the function
\[
\tau\colon\chi_t(\Sigma,\psl)\to\R,
\]
which we have just defined.
\end{defn}
It is a classical result of Bers~\cite{bers} and Kerckhoff~\cite{kerckhoff} that the pair $(\ell_a,\tau)$ defines a smooth (even real-analytic) diffeomorphism from each of the two connected components of $\chi_t(\Sigma,\psl)$ to $\R_{>0}\times\R$ for every $t<2$, see also~\cite[Section~10.6]{mcg-primer}. Moreover, if $\omega_\mathcal{G}$ denotes the Goldman symplectic form on $\chi_t(\Sigma,\psl)$, then Wolpert\footnote{Formula~\eqref{eq:wolpert-formula-one-holed-torus} and Wolpert's original formula differ by a minus sign because we are using a different sign convention in the definition of Hamiltonian vector fields, namely $\omega(X_h,\cdot)=dh$ as in~\eqref{eq:Hamiltonian-vf}, which leads to different duality formulas between lengths and twists.} in~\cite{wolpert-duality-formula, wolpert-formula} proved that 
\begin{equation}\label{eq:wolpert-formula-one-holed-torus}
\omega_\mathcal{G}=-d\ell_a\wedge d\tau.
\end{equation}
This means that $\tau$ is dual to $\ell_a$ in the sense that $\tau$ parametrizes Hamiltonian flow orbits of the function $\ell_a\colon\chi_t(\Sigma,\psl)\to\R_{>0}$. In other words, if $\Phi_a^t$ denotes the Hamiltonian flow of the function $\ell_a$---known as the \emph{twist flow} of $a$---at time $t$, then it holds that
\begin{equation}\label{eq:Goldman-formula-one-holed-torus}
\tau\circ\Phi_a^t = \tau + t.
\end{equation}

The twist flow of $a$ has been described explicitly by Goldman~\cite{goldman-invariant-functions}. It descends from the following flow on representations. For $t\in \R$, let $X(t)$ denote the element in the Lie algebra of $\psl$ given by
\[
X(t)=\begin{pmatrix}
    t/2 & 0\\ 0 & -t/2
\end{pmatrix}.
\]
Given a representation $\rho\colon\pi_1\Sigma\to\psl$ with $\Tr\rho(c)<2$, there exists $A_\rho$ in $\psl$ such that
\[
\rho(a)=A_\rho\exp\big(X\big(\ell_{a}(\rho)\big)\big)A_\rho^{-1},
\]
where $\exp$ denotes the Lie exponential map. For simplicity, we will write
\[
\xi_\rho(t)=A_\rho\exp\big(X(t)\big)A_\rho^{-1}.
\]
Note that $\xi_\rho(t)$ is a hyperbolic translation of length $t$ along the invariant axis of $\rho(a)$. In other words, $\{\xi_\rho(t):t\in\R\}$ is the centralizer of $\rho(a)$. The flow orbits $\Phi_{a}^t\big([\rho])$ inside the relative character variety can be lifted to the space of representations as
\begin{equation}\label{eq:twist-flow-one-holed-torus}
\overline{\Phi}_{a_1}^t(\rho):\begin{cases}
    a\mapsto \rho(a)\\
    b\mapsto \rho(b)\xi_\rho(t).
\end{cases}
\end{equation}
Observe that $\{\rho(b)\xi_\rho(t):t\in\R\}$ is the set of all hyperbolic transformations mapping the axis of $\rho(a)$ to its image by $\rho(b)$. 

\subsection{Polygonal model}\label{sec:polygonal-model-one-holed-tori}
A common way to construct hyperbolic surfaces is to glue hyperbolic polygons along edges of the same length. Many constructions can be found in the literature. The one we present below was inspired by a construction of Parker--Parkkonen~\cite[Section~1]{parker-parkonnen}. 

\begin{figure}[h]
\centering
\begin{tikzpicture}
\node[anchor = south west, inner sep=0mm] at (0,0) {\includegraphics[width=11cm]{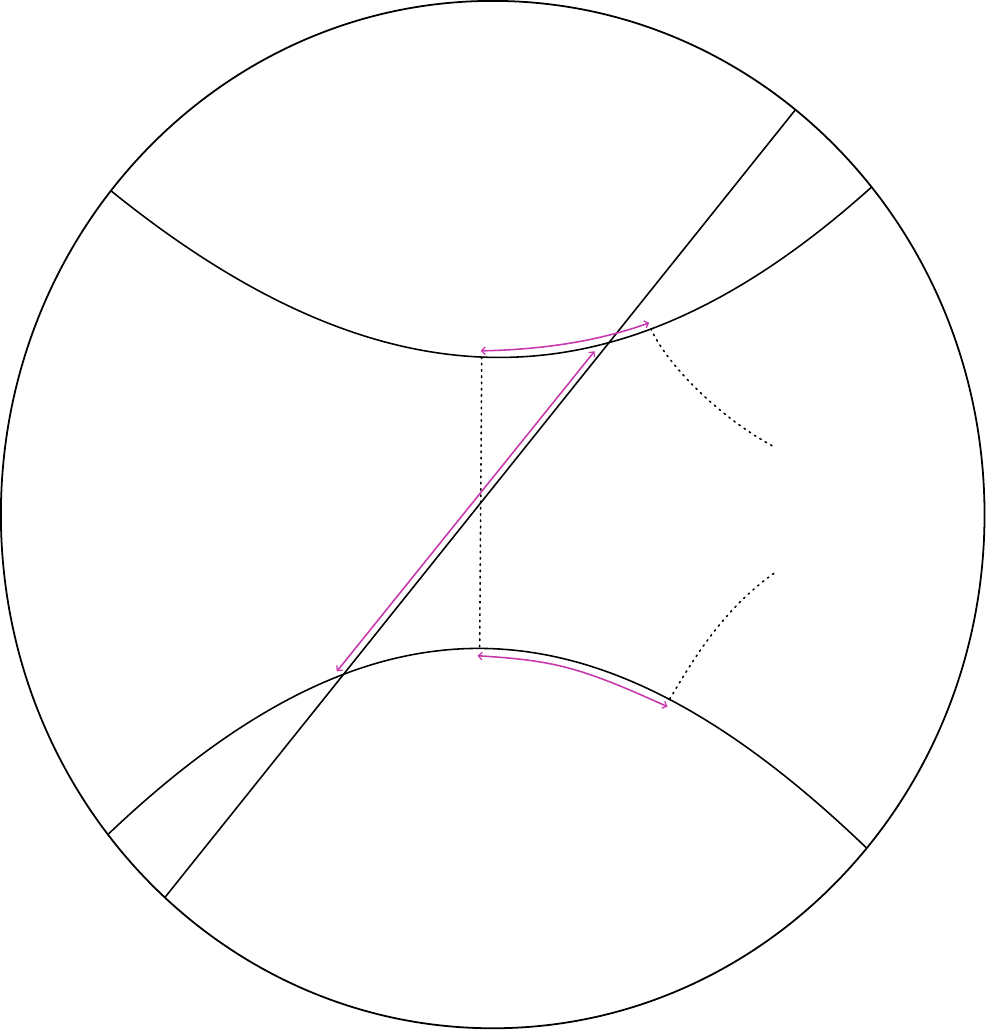}};

\node at (5.1,4.05) {$H$};
\node at (5.25,7.75) {$H'$};
\node at (3.95,3.8) {$X$};
\node at (6.82,7.99) {$X'$};
\node at (7.4,3.4) {$Y$};
\node at (7.45,8.07) {$Y'$};
\node at (5.1,6.05) {$M$};

\node at (9,2.35) {$g_a$};
\node at (9,9.25) {$g_a'$};
\node at (2.4,1.8) {$g_b$};
\node at (8.25,4.5) {$s$};
\node at (8.25,7.1) {$s'$};

\node[mauve] at (6.3,3.7) {$\ell_a/2$};
\node[mauve] at (6.1,7.87) {$\ell_a/2$};
\node[mauve] at (4.25,5) {$\ell_b$};
\end{tikzpicture}
\caption{Construction of our polygonal model.}
\label{fig:polygon-construction}
\end{figure}

Let $[\rho]\in\chi_t(\Sigma,\psl)$ for some $t<2$. The reader may want to have Figure~\ref{fig:polygon-construction} at hand to follow the construction. Fact~\ref{fact:tr<2-implies-crrossing-axes} tells us that $\rho(a)$ and $\rho(b)$, where $a$ and $b$ are the generators of $\pi_1\Sigma$ from~\eqref{eq:presentation-pi_1-one-holed-torus}, are both hyperbolic transformations of the hyperbolic plane whose axes $g_a$ and $g_b$ intersect. We label their intersection point $X$. The image by $\rho(b)$ of $g_a$ is denoted $g_a'$. It is the invariant axis of the hyperbolic transformation $\rho(ba^{-1}b^{-1})$ and contains the image $X'$ of $X$ by $\rho(b)$. Note that $X'$ is the intersection point of $g_b$ and $g_a'$. So,
\[
d(X,X')=\ell_b([\rho])=2\arccosh\big(\vert\Tr\rho(b)\vert/2\big).
\]
The geodesics $g_a$ and $g_a'$ have a common perpendicular geodesic line which intersects $g_a$ at $H$ and $g_a'$ at $H'$. It is possible that $X=H$, in which case $X'=H'$. If $X\neq H$, then we let $M$ be the intersection points of the hyperbolic segments $[XX']$ and $[HH']$.

\begin{claim}
The point $M$ is the midpoint of both segments $[XX']$ and $[HH']$.
\end{claim}
\begin{proofclaim}
The hyperbolic triangles $XHM$ and $X'H'M$ have the same interior angles by construction, so they are isometric. This implies $XM=X'M'$ and $HM=H'M$.
\end{proofclaim}

We continue our construction by letting $Y$ be the point on $g_a$ at distance $\ell_a([\rho])/2$ from~$H$ in the direction of the translation $\rho(a)$. In other words, $Y$ is the image of $H$ by the hyperbolic translation that squares to $\rho(a)$. Similarly, we let $Y'$ be the point on $g_a'$ at distance $\ell_a([\rho])/2$ from $H'$ and on the same side of the geodesic line $(HH')$ as $Y$. Now, consider $s$ the perpendicular line to $g_a$ through $Y$, and $s'$ the perpendicular line to $g_a'$ through $Y'$. The reason for introducing $s$ and $s'$ is that we can use them to recover $\rho(c)$. Let $\sigma$ and $\sigma'$ denote the orientation-reversing reflections through the geodesic lines~$s$ and~$s'$.

\begin{claim}\label{claim:rho(c)=sigma-circ-sigma'}
It holds that $\rho(c)=\sigma\sigma'$.
\end{claim}
\begin{proofclaim}
Let $\varsigma$ be the reflection through the geodesic line $(HH')$. Since $d(H,Y)=\ell_a([\rho])/2$, we have $\rho(a)=\sigma\varsigma$. Moreover, $\rho(ba^{-1}b)$ is the hyperbolic translation along $g_a'$ by $\ell_a([\rho])$ that maps $Y'$ beyond $H'$. In other words, $\rho(ba^{-1}b)=\varsigma\sigma'$. So, $\rho(c)=\rho(aba^{-1}b^{-1})=\sigma\sigma'$.
\end{proofclaim}

\begin{figure}[h]
\centering
\begin{tikzpicture}
\node[anchor = south west, inner sep=0mm] at (0,0) {\includegraphics[width=6.5cm]{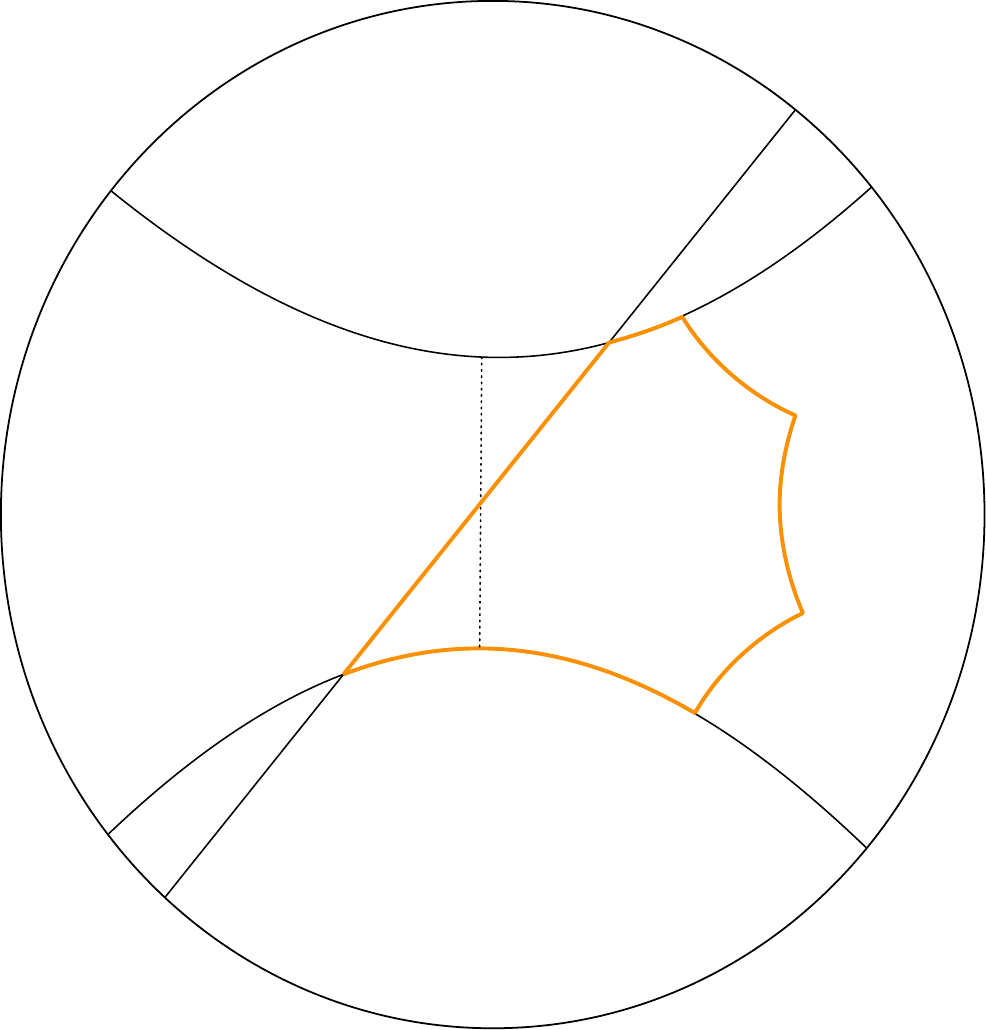}};

\node[apricot] at (5.5,4.1) {$B'$};
\node[apricot] at (5.48,2.76) {$B$};
\node[apricot] at (2.37,2.16) {$X$};
\node[apricot] at (4.5,1.9) {$Y$};
\node[apricot] at (3.9,4.75) {$X'$};
\node[apricot] at (4.68,4.95) {$Y'$};

\node[anchor = south west, inner sep=0mm] at (7.2,0) {\includegraphics[width=6.5cm]{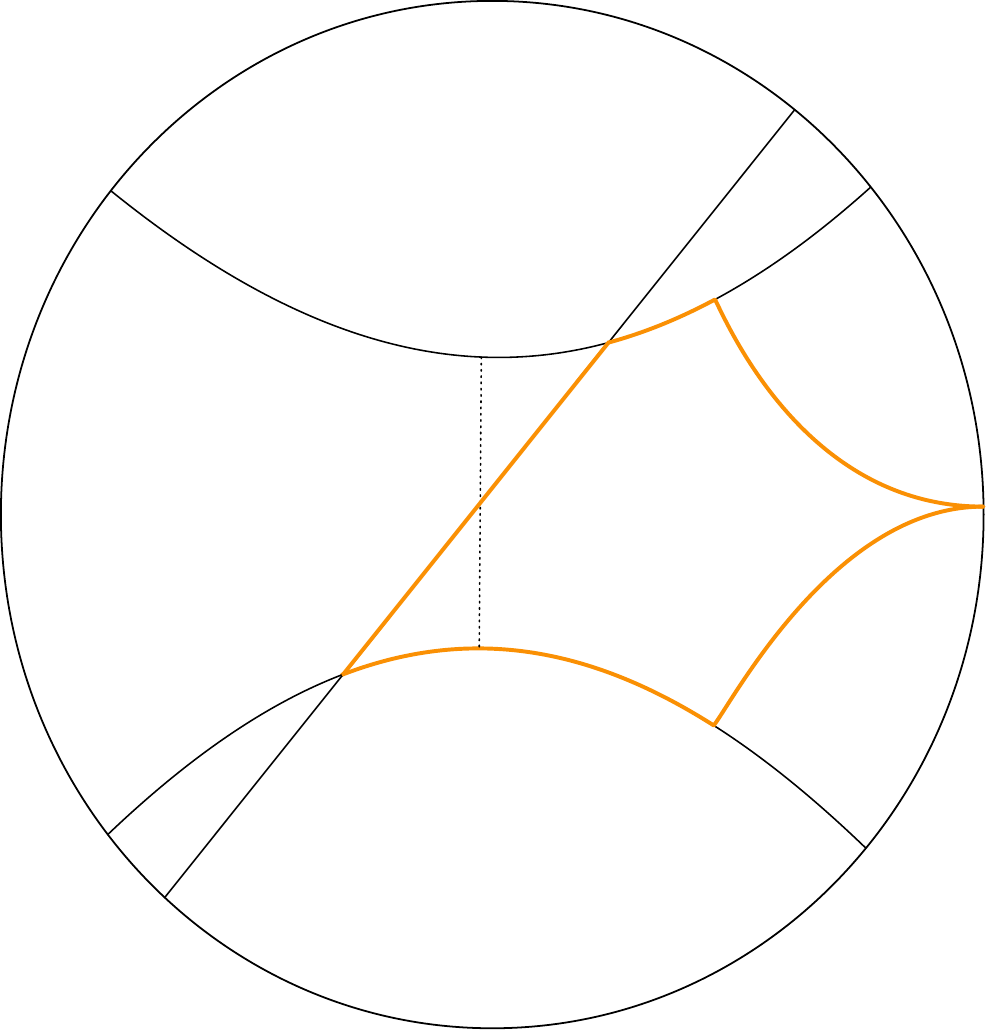}};

\node[apricot] at (9.57,2.16) {$X$};
\node[apricot] at (11.8,1.8) {$Y$};
\node[apricot] at (11.1,4.75) {$X'$};
\node[apricot] at (11.88,4.95) {$Y'$};
\node[apricot] at (13.9, 3.45) {$B$};
\end{tikzpicture}
\begin{tikzpicture}
\node[anchor = south west, inner sep=0mm] at (0,0) {\includegraphics[width=6.5cm]{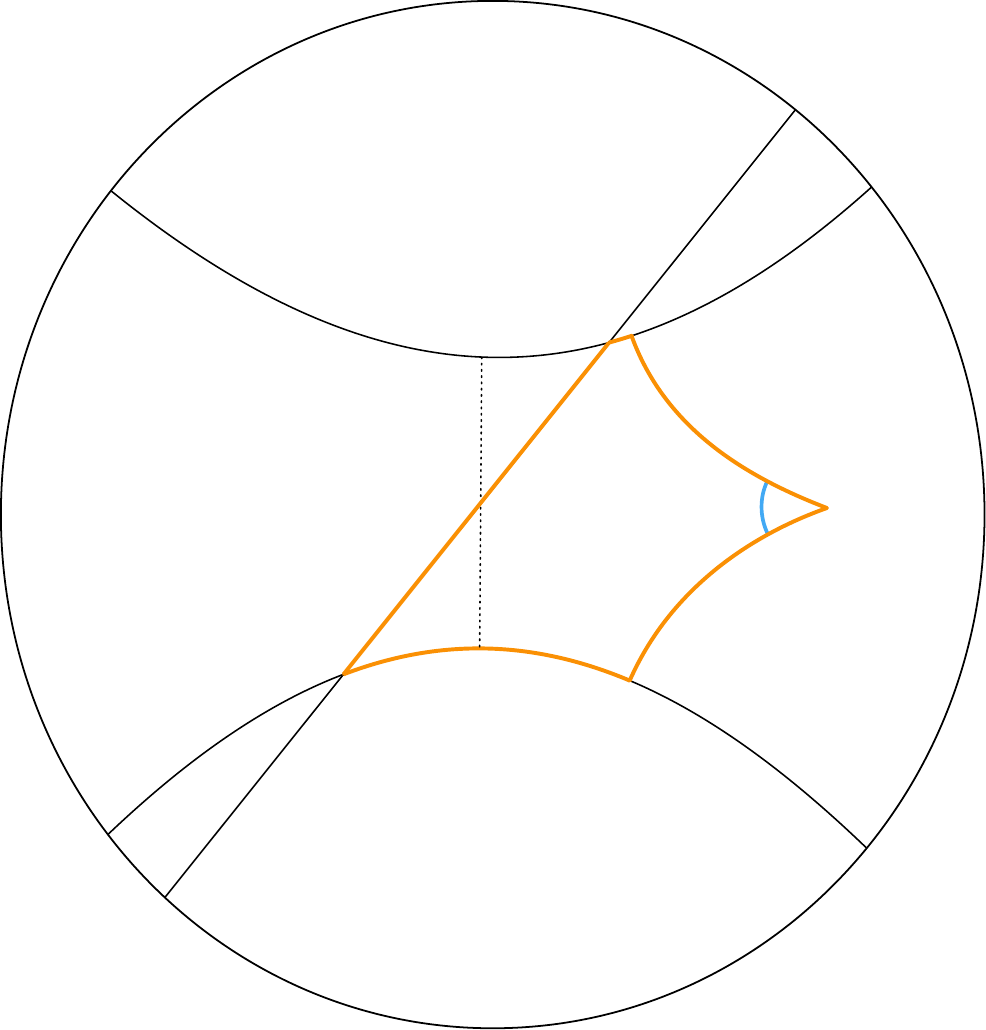}};

\node[apricot] at (2.37,2.16) {$X$};
\node[apricot] at (4.15,2.1) {$Y$};
\node[apricot] at (3.9,4.75) {$X'$};
\node[apricot] at (4.52,4.5) {$Y'$};
\node[apricot] at (5.6,3.45) {$B$};
\node[sky] at (4.4,3.45) {\footnotesize $\pi-\beta/2$};
\end{tikzpicture}
\caption{On top, from the left: our polygonal model in the hyperbolic and parabolic regimes. On the bottom: the elliptic case.}
\label{fig:polygon-construction-3-types}
\end{figure}

Depending on whether the geodesic lines $s$ and $s'$ intersect, the three configurations illustrated on Figure~\ref{fig:polygon-construction-3-types} may occur.
\begin{enumerate}
    \item[(H)] The geodesic lines $s$ and $s'$ do not intersect in the closure of the hyperbolic plane. In that case, they have a common perpendicular line which intersects $s$ at $B$ and $s'$ at $B'$. Furthermore, $\rho(c)=\sigma\sigma'$ is a hyperbolic translation along $(B'B)$ of length $2d(B',B)$. This case thus occurs when $t<-2$; more precisely, when $t=-2\cosh(d(B',B))$.
    \item[(P)] The geodesic lines $s$ and $s'$ intersect at a point $B$ in the boundary of the hyperbolic plane. In that case, $\rho(c)=\sigma\sigma'$ is a parabolic transformation that fixes $B$ and $t=-2$.
    \item[(E)] The geodesic lines $s$ and $s'$ intersect at a point $B$ inside the hyperbolic plane. This case occurs when $-2<t<2$. If $t=2\cos(\beta/2)$ with $\beta\in (0,2\pi)$, then the angle at $B$ between the lines $s$ and $s'$ is equal to $\pi-\beta/2$. So, when $\beta$ approaches $2\pi$, the angle tends to $0$ and $t$ tends to~$-2$, and we recover the parabolic regime. In this case, $\rho(c)=\sigma\sigma'$ is a rotation of angle $\beta$ around $B$. Whether the rotation is clockwise or anti-clockwise depends on the connected component of $\chi_t(\Sigma,\psl)$ (see Lemma~\ref{lem:direction-of-rotation-one-holed-torus}).
\end{enumerate}

We just explained how to associate to every point of $\chi_t(\Sigma,\psl)$ a hexagon $XYBB'Y'X'$ when $t<-2$ and a pentagon $XYBY'X'$ when $-2\leq t< 2$. Strictly speaking, a polygon is associated to a representation and an \emph{isometry class} of polygons is associated to a point of $\chi_t(\Sigma,\psl)$, but we will often neglect this technicality. The polygons are illustrated on Figure~\ref{fig:polygon-construction-3-types}.
\begin{defn}\label{defn:polygonal-model}
This association is what we refer to as the \emph{polygonal model} of $\chi_t(\Sigma,\psl)$ for $t<2$. 
\end{defn}
One can think of them as a ``half fundamental domains'' that encodes all the important geometric quantities associated with the hyperbolic torus. They actually distinguish the points of $\chi_t(\Sigma,\psl)$ for $t<2$, so that no two different conjugacy classes of representations are associated to isometric polygons (Corollary~\ref{cor:polygons-determine-representations-one-holed-tori}). 

\begin{figure}[h]
\centering
\begin{tikzpicture}
\node[anchor = south west, inner sep=0mm] at (0,0) {\includegraphics[width=6cm]{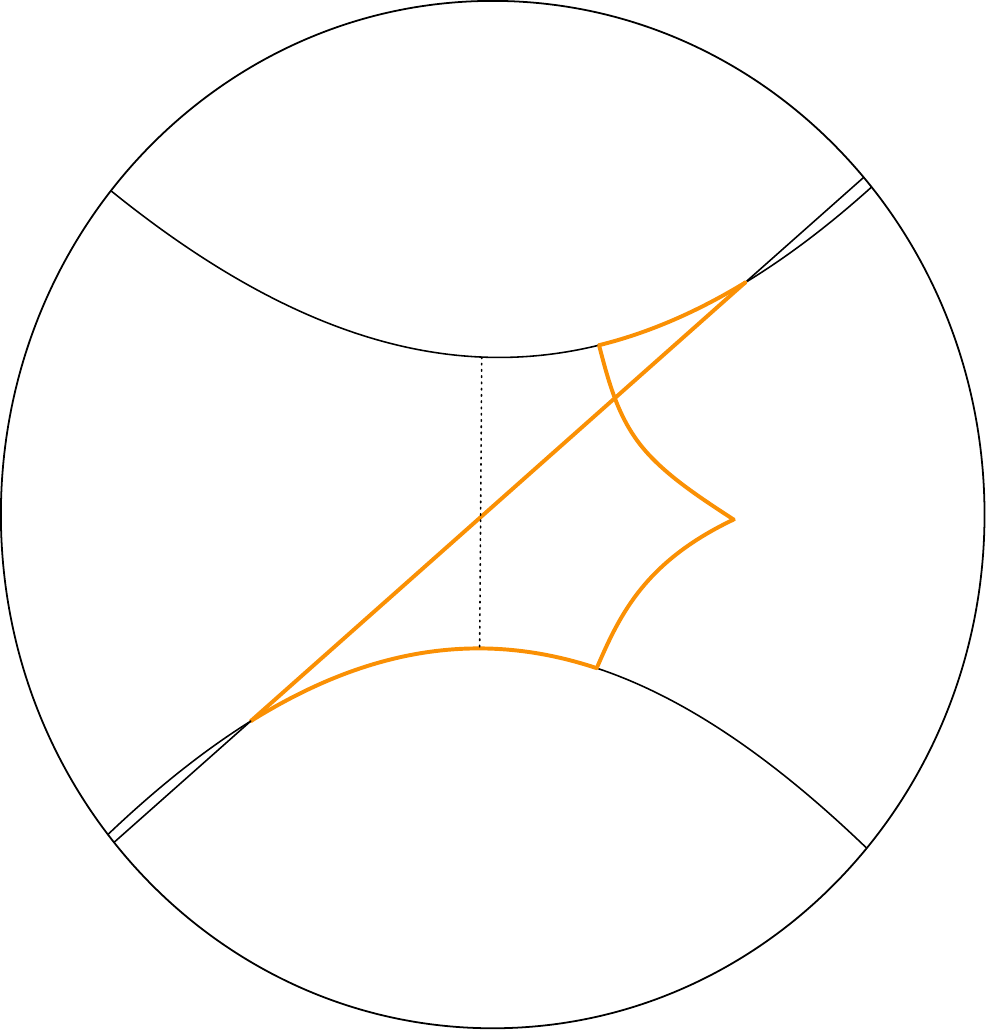}};
\end{tikzpicture}
\caption{Example of a non-convex polygon.}
\label{fig:polygon-non-convex}
\end{figure}

There are two important observations to keep in mind. First, those polygons may not be convex (see Figure~\ref{fig:polygon-non-convex}). This typically occurs when the angle between $g_a$ and $g_b$ is small, or, as we will see in Proposition~\ref{prop:twist-coordinate-is-geometric-distance}, when the twist coordinate is large.

Second, the polygons are oriented and each of the two components of $\chi_t(\Sigma,\psl)$ gives rise to polygons with opposite orientation. The polygons have a signed area. For instance, in the upper half-plane, the area may be computed by integration against the symplectic form $(dx\wedge dy)/y^2$, which ensures that anti-clockwise oriented triangles have positive area. As we will see in Lemma~\ref{lem:toledo-and-l_a-one-holed-torus}, the signed area of the polygon associated to $[\rho]\in\chi_t(\Sigma,\psl)$ is directly related to the \emph{Toledo number} of $[\rho]$ defined by Burger--Iozzi--Wienhard~\cite{BIW}. Recall that for $[\rho]\in \chi_t(\Sigma,\psl)$, we have
\begin{equation}\label{eq:Toldeo-number-one-holed-torus}
\Tol([\rho])=\begin{cases}
    \pm 1 & \text{if } t\leq -2\\
    \pm \frac{\beta}{2\pi} & \text{if } t=2\cos(\beta/2) \in(-2,2)\\
    0 & \text{if } t\ge2.
\end{cases}
\end{equation}

The two connected components of $\chi_t(\Sigma,\psl)$ are distinguished by the sign of the Toledo number when $t<2$. It can be determined from the oriented angle $\theta$ between the axes of $\rho(a)$ and $\rho(b)$. More precisely, if we switch to the Poincaré disk model and align the axis of $\rho(a)$ with the real axis so that the attracting fixed point of $\rho(a)$ is $1$ and $X$ is now the origin, then the attracting fixed point of $\rho(b)$ becomes $e^{i\theta}$ for some $\theta \in (0,\pi)\cup (\pi,2\pi)$. If $\theta\in(0,\pi)$, then $\Tol([\rho])< 0$, while $\Tol([\rho])> 0$ when $\theta\in(\pi,2\pi)$. For instance, all three polygons illustrated on Figures~\ref{fig:polygon-construction-3-types} and~\ref{fig:polygon-non-convex} have $\theta\in (0,\pi)$.

\begin{lem}\label{lem:direction-of-rotation-one-holed-torus}
Let $t=2\cos(\beta/2)$ for some $\beta\in (0,2\pi)$ and take $[\rho]\in\chi_t(\Sigma,\psl)$. 
\begin{itemize}
    \item If $\Tol([\rho])< 0$, then $\rho(c)$ is a \emph{clockwise} rotation of angle $\beta$ around $B$.
    \item If $\Tol([\rho])> 0$, then $\rho(c)$ is an \emph{anti-clockwise} rotation of angle $\beta$ around $B$.
\end{itemize}
\end{lem}
\begin{proof}
Recall from the discussion after Claim~\ref{claim:rho(c)=sigma-circ-sigma'} that $\rho(c)$ is always a rotation of angle $\beta$ around $B$. The case $\Tol([\rho])< 0$, or equivalently $\theta\in(0,\pi)$, is the one illustrated on Figure~\ref{fig:polygon-construction}, assuming the geodesic lines $s$ and $s'$ meet. In that case, $\rho(c)=\sigma\sigma'$ is indeed a clockwise rotation of angle $\beta$ around $B$. The other case is flipped.
\end{proof}

For future reference, we compile the conclusions of Lemma~\ref{lem:direction-of-rotation-one-holed-torus} and the properties of the polygonal model in the elliptic regime in the following table. Let $\beta\in (0,2\pi)$ and let us set $t=2\cos(\beta/2)$.
\begin{table}[h]
    \centering
    \begin{tabular}{ccc}
        \toprule
         & $\Tol> 0$ & $\Tol< 0$  \\ 
         \midrule 
         $\chi_t(\Sigma,\psl)$ & \makecell{anti-clockwise\\ pentagon angle $=\pi-\beta/2$} & \makecell{clockwise\\ pentagon angle $=\pi-\beta/2$} \\ \midrule 
         $\chi_{-t}(\Sigma,\psl)$ & \makecell{clockwise\\ pentagon angle $=\beta/2$} & \makecell{anti-clockwise\\ pentagon angle $=\beta/2$}\\
         \bottomrule
    \end{tabular}
    \caption{Each cell corresponds to one of the two connected components of the relative character varieties $\chi_t(\Sigma,\psl)$ and $\chi_{-t}(\Sigma,\psl)$. They indicate the direction of the rotation $\rho(c)$ around $B$ ($\rho(c)$ always rotates by $\beta$) and the interior angle at $B$ in the polygonal model of $[\rho]$.}
    \label{tab:different-cases-ellipitc-regime}
\end{table}

\begin{figure}[h]
\centering
\begin{tikzpicture}[scale=1.64]
\node[anchor = south west, inner sep=0mm] at (3.15,3.15) {\includegraphics[width=11cm]{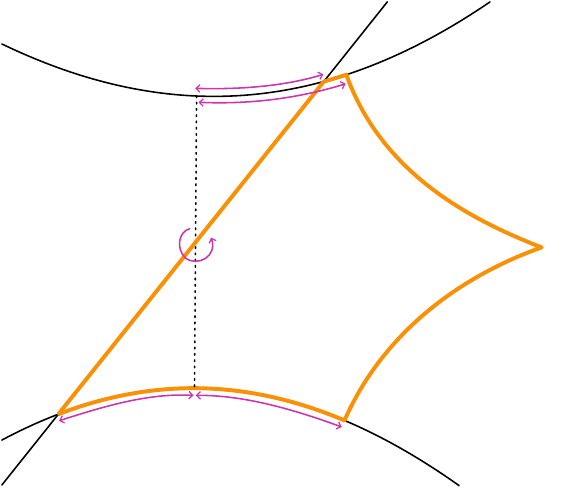}};

\node[mauve] at (4.7,4) {$\tau/2$};
\node[mauve] at (6.1,8) {$\tau/2$};
\node[mauve] at (6.05, 7.45) {$\ell_a/2$};
\node[mauve] at (6.05,3.95) {$\ell_a/2$};
\node[mauve] at (5.8,5.9) {$\Phi_a^t$};

\node at (3.7,4.15) {$X$};
\node at (5.55,4.45) {$H$};
\node at (7.1,4.15) {$Y$};
\node at (5.1,6) {$M$};
\node at (5.3,7.85) {$H'$};
\node at (6.9,8.1) {$X'$};
\node at (7.4,7.87) {$Y'$};
\end{tikzpicture}
\caption{How to read Fenchel--Nielsen coordinates and visualize the twist flow of $a$ in our polygonal model.}
\label{fig:fenchel-nielsen-coordinates-one-holed-torus}
\end{figure}

\begin{lem}\label{lem:toledo-and-l_a-one-holed-torus}
The Toldeo number of $[\rho]\in\chi_t(\Sigma,\psl)$ is equal to the signed area of the polygon corresponding to $[\rho]$ times $-1/\pi$. Furthermore,
\[
\ell_a([\rho])=\varepsilon d(X,Y)+\varepsilon' d(X',Y'),
\]
where $\varepsilon = 1$ if $X$ belongs to the geodesic ray $[YH)$ and $\varepsilon = -1$ otherwise, and similarly $\varepsilon'=1$ if $X'$ belongs to the geodesic ray $[Y'H')$ and $\varepsilon'=-1$ otherwise.
\end{lem}
\begin{proof}
Recall that the positive area of a convex hyperbolic $n$-gon in the hyperbolic plane is equal to $(n-2)\pi$ minus the sum of interior angles. Its signed area is positive if the vertices of the polygon are anti-clockwise oriented, and negative otherwise. In our case, the polygons $XYBB'Y'X'$ and $XYBY'X'$ may not be convex, so their signed area is given by the sum of the signed area of each convex piece. It is not hard to see that when the oriented angle $\theta$ between the axes of $\rho(a)$ and $\rho(b)$ is contained in $(0,\pi)$, then the signed area of $XYBB'Y'X'$ is $\pi$, while the signed area of $XYBY'X'$ is $\beta/2$, respectively $\pi$ in the case where $B$ lies on the boundary of the hyperbolic plane. Multiplying these numbers by $-1/\pi$ recovers the expression of the Toledo number of $[\rho]$ provided in~\eqref{eq:Toldeo-number-one-holed-torus} and the sign discussion thereafter. The case $\theta\in(\pi,2\pi)$ is analogous.

The statement about the length $\ell_a([\rho])$ can be seen by recalling from the construction that $\ell_a([\rho])=d(H,Y)+d(H',Y')$ and working out the different cases according to the relative positions of $X$ and $X'$, see Figure~\ref{fig:fenchel-nielsen-coordinates-one-holed-torus}.
\end{proof}

More interestingly, the polygonal model introduced in Definition~\ref{defn:polygonal-model} provides a simple geometric interpretation of the twist coordinate and the twist flow of $a$, as illustrated on Figure~\ref{fig:fenchel-nielsen-coordinates-one-holed-torus}.

\begin{prop}\label{prop:twist-coordinate-is-geometric-distance}
The twist coordinate $\tau$ of $[\rho]\in\chi_t(\Sigma,\psl)$ is equal to $2\varepsilon d(X,H)$, where $\varepsilon=-1$ if $X$ belongs to the geodesic ray $[HY)$ and $\varepsilon=1$ otherwise. Moreover, the action of the twist flow defined in~\eqref{eq:twist-flow-one-holed-torus} can be visualized on the polygonal model of $[\rho]$ by a rotation of the axis of $\rho(b)$ around $M$. In other words, it moves the points $X$ and $X'$, but preserves all the other points.
\end{prop}
\begin{proof}
The data of $[\rho]$ corresponds to a marked hyperbolic structure $f\colon\Sigma\to\Upsilon$, uniquely defined up to isotopy. The universal cover $\widetilde\Upsilon$ is isometric to the upper half-plane $\HH$ (in a unique way up to Möbius transformations). The marking $f$ lifts to universal covers as a $\rho$-equivariant homeomorphism $\widetilde f\colon\widetilde\Sigma\to\HH$. Among all the possible geodesic lines in $\widetilde\Upsilon\cong\HH$ that lie above the geodesic representatives $\alpha$ and $\beta$ of $f_\ast a$ and $f_\ast b$ on $\Upsilon$, two are the axes of $\rho(a)$ and $\rho(b)$. We label them $\widetilde\alpha$ and $\widetilde\beta$. They are naturally oriented.

Now, when we cut $\Upsilon$ along $\alpha$, we obtain a hyperbolic pair of pants with two geodesic boundary curves $\alpha_1$ and $\alpha_2$, and a third boundary which is either a geodesic boundary curve, a cusp, or a cone point. If we further cut along the shortest geodesic arcs connecting the third boundary to each of the two geodesic boundary curves $\alpha_1$ and $\alpha_2$, we obtain a hyperbolic polygon $\mathcal{P}$ like on Figure~\ref{fig:hexagon-universal-cover}. Two opposite sides of $\mathcal{P}$ are of length $\ell_a([\rho])$ and correspond to $\alpha_1$ and $\alpha_2$. They are joint by a hyperbolic segment that is perpendicular to both. Back on the hyperbolic pair of pants, this segment was the shortest geodesic arc $\gamma$ joining $\alpha_1$ and $\alpha_2$ used in the definition of the twist coordinate (Definition~\ref{defn:twist-coordinate-one-holed-torus}). The geodesic $\beta$ may not necessarily be a connected line in $\mathcal{P}$. It may instead be made of several disjoint line segments, one of which, call it $t$, passes through the centre of $\mathcal{P}$. Note that $t$ inherits an orientation from $\beta$. There is a unique way of embedding $\mathcal{P}$ in $\HH$ if we make sure that~$t$ lies on $\widetilde\beta$ with matching orientations, and the two sides of length $\ell_a([\rho])$ lie on $\widetilde\alpha$ and $\widetilde\alpha'=\rho(b)\widetilde\alpha$, again with matching orientations. 

\begin{figure}[h]
\centering
\begin{tikzpicture}
\node[anchor = south west, inner sep=0mm] at (0,0) {\includegraphics[width=11cm]{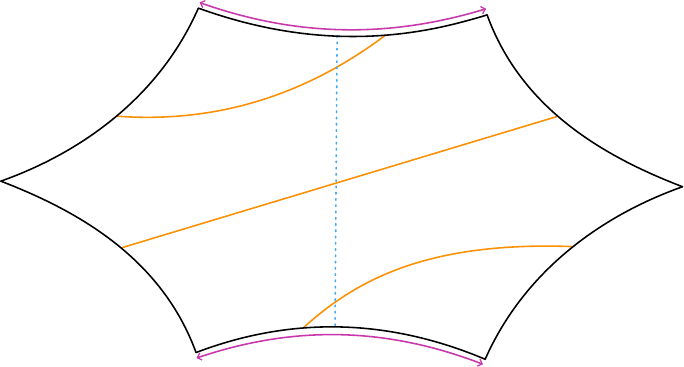}};

\node[apricot] at (7.5,3.2) {$t$};
\node[mauve] at (5.5,5.75) {$\ell_a([\rho])$};
\node[mauve] at (5.5,0.2) {$\ell_a([\rho])$};
\end{tikzpicture}
\caption{An example of a polygon $\mathcal{P}$ (this one comes the hyperbolic torus of Figure~\ref{fig:one-holed-torus-cut-and-uncut}).}
\label{fig:hexagon-universal-cover}
\end{figure}

The polygon $\mathcal{P}$ is almost our polygonal model for $[\rho]$. We recover it by first cutting $\mathcal{P}$ along $t$ and keeping the piece to the right of $t$ with respect to its orientation. If $t$ does not intersect $\widetilde\alpha$ and $\widetilde\alpha'$, then we extend $t$ on both sides until it does. The meeting points are $X$ and $X'$, respectively. The resulting polygon is our polygonal model for $[\rho]$.

The geodesic $\beta$ on $\Upsilon$ is homotopic to the concatenation of $\gamma$ and a unique arc going along $\alpha$ whose signed length is the twist coordinate of $[\rho]$ (Definition~\ref{defn:twist-coordinate-one-holed-torus}). Let $H$ and $H'$ denote the intersection points of the common perpendicular to $\widetilde\alpha$ and $\widetilde\alpha'$ with respectively $\widetilde\alpha$ and $\widetilde\alpha'$. The hyperbolic segment $[HH']$ is a lift of $\gamma$, while the hyperbolic segment $[XX']$ is a lift of $\beta$. This means that the desired arc going along $\alpha$ is the projection of $[XH]\cup [H'X']$. We conclude that $\tau([\rho])=2d(X,H)$ when $X$ lies to the left of $H$ with respect to the orientation of $\widetilde\alpha$ or equivalently when $X\notin [HY)$, and $\tau([\rho])=-2d(X,H)$ otherwise.

It remains to explain how the twist flow of $a$ can be visualized in our polygonal model. Using the explicit form~\eqref{eq:twist-flow-one-holed-torus}, we see that only the axis of the image of $b$ is affected by the flow twist of $a$. Combining Equation~\eqref{eq:Goldman-formula-one-holed-torus} with the first statement of Proposition~\ref{prop:twist-coordinate-is-geometric-distance} we just proved, we obtain that the axis of $\overline{\Phi}^t_a(\rho)(b)$ is the geodesic line through the points $X_t$ on $\widetilde\alpha$ and $X_t'$ on $\widetilde\alpha'$ which are at respective signed distance $t/2$ of $X$ and $X'$ with respect to the orientations of $\widetilde\alpha$ and $\widetilde\alpha'$. In other words, the axis of $\overline{\Phi}^t_a(\rho)(b)$ is obtained by rotating the axis of $\rho(b)$ around the midpoint $M$ of $[HH']$ by a suitable angle.
\end{proof}

\begin{cor}\label{cor:polygons-determine-representations-one-holed-tori}
Our polygonal model that associates to every point in $\chi_t(\Sigma,\psl)$ an orientation-preserving isometry class of polygons in the hyperbolic plane is injective. In other words, polygons fully determine representations.
\end{cor}
\begin{proof}
The value of $t$ is determined by the angle at $B$ or the length $d(B,B')$. Fenchel--Nielsen coordinates can be measured on polygons using Lemma~\ref{lem:toledo-and-l_a-one-holed-torus} and Proposition~\ref{prop:twist-coordinate-is-geometric-distance}. Finally, the orientation of the vertices (or the signed area of the polygons) determines the Toledo number by Proposition~\ref{prop:twist-coordinate-is-geometric-distance}, hence also determines a connected component of $\chi_t(\Sigma,\psl)$.
\end{proof}

\subsection{Explicit parametrization}\label{sec:one-holed-torus-explicit-paramaetrization}
\begin{figure}[h]
\centering
\begin{tikzpicture}
\node[anchor = south west, inner sep=0mm] at (0,0) {\includegraphics[width=14cm]{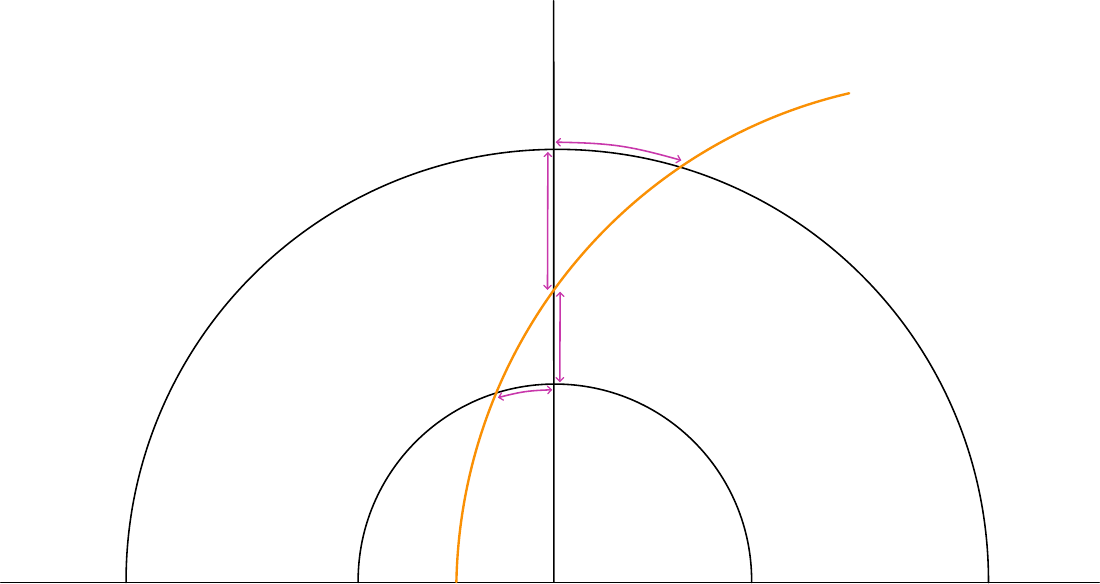}};

\node[mauve] at (7.5,3.15) {\small $\eta/2$};
\node[mauve] at (6.6,4.7) {\small $\eta/2$};
\node[mauve] at (7.85,5.8) {\small $\tau/2$};
\node[mauve] at (6.65,2.2) {\small $\tau/2$};

\node at (4.4,-.3) {\small $-\coth(\lambda/2)$};
\node at (1.8,-.3) {\small $-e^\eta\coth(\lambda/2)$};
\node at (9.7,-.3) {\small $\coth(\lambda/2)$};
\node at (12.7,-.3) {\small $e^\eta\coth(\lambda/2)$};

\node at (7.25, 2.35) {\small $H$};
\node at (6.8, 5.75) {\small $H'$};
\node at (6.15, 2.6) {\small $X$};
\node at (8.7, 5.05) {\small $X'$};
\node at (6.7, 3.8) {\small $M$};
\end{tikzpicture}
\caption{Explicit parametrization in the upper half-plane.}
\label{fig:polygon-upper-half-plane}
\end{figure}

It is possible to write explicit formulas for the images of the generators $a$ and $b$ under a representation $\rho$ representing a point of $\chi_t(\Sigma,\psl)$ with Fenchel--Nielsen coordinates $(\ell,\tau)$ when $t<2$. We include them as we think they may be useful for experiments and computer simulations. We derived the formulas below from those of Parker--Parkkonen in the parabolic case $\Tr\rho(c)=-2$~\cite[Section~1]{parker-parkonnen}. The reader may want to have Figure~\ref{fig:polygon-upper-half-plane} at hand to follow the construction.

We start by introducing two auxiliary parameters. Let $\lambda=\ell/2$ and $\eta>0$ be given by
\[
\cosh(\eta)=\frac{\cosh(\lambda)^2-t/2}{\sinh(\lambda)^2}.
\]
We let
\begin{equation}\label{eq:rho(a)-one-holed-torus}
\rho_{\ell,\tau}(a)=\pm\begin{pmatrix}
    \cosh(\lambda) & \cosh(\lambda)+1\\
    \cosh(\lambda)-1 & \cosh(\lambda)
\end{pmatrix}.
\end{equation}
Note that $\rho_{\ell,\tau}(a)$ is independent of $\tau$. The axis of $\rho_{\ell,\tau}(a)$ is the geodesic line in the upper half-plane with endpoints $\pm\coth(\lambda/2)$ on the real line (see Figure~\ref{fig:polygon-upper-half-plane}). The attracting fixed point of $\rho_{\ell,\tau}(a)$ is $\coth(\lambda/2)$ and the repelling one is $-\coth(\lambda/2)$. We may assume that $\rho_{\ell,\tau}(b)$ is a translation along the imaginary axis when $\tau=0$. Solving the equation $\Tr(\rho_{\ell,0}([a,b]))=t$ to determine the diagonal entries of $\rho_{\ell,0}(b)$ leads to
\[
\rho_{\ell,0}(b)=\pm\begin{pmatrix}
    e^{\eta/2} & 0\\
    0 & e^{-\eta/2}
\end{pmatrix}.
\]
Observe that $\rho_{\ell,0}(b)$ maps the axis of $\rho_{\ell,\tau}(a)$ to the geodesic line with endpoints $\pm e^\eta\coth(\lambda/2)$. The common perpendicular segment to these two geodesic lines has endpoints $H=i\coth(\lambda/2)$ and $H'=ie^{\eta}\coth(\lambda/2)$. The midpoint of $[HH']$ is $M=ie^{\eta/2}\coth(\lambda/2)$. In order to determine $\rho_{\ell,\tau}(b)$ for an arbitrary $\tau$, we apply the twist flow associated to the curve $a$ to $\rho_{\ell,0}(b)$ using Equation~\eqref{eq:twist-flow-one-holed-torus}. This leads to
\begin{equation}\label{eq:rho(b)-one-holed-torus}
    \rho_{\ell,\tau}(b)=\pm\begin{pmatrix}
        \cosh(\tau/2)e^{\eta/2} & \sinh(\tau/2)e^{\eta/2}\coth(\lambda/2)\\
        \sinh(\tau/2)e^{-\eta/2}\tanh(\lambda/2) & \cosh(\tau/2)e^{-\eta/2}
    \end{pmatrix}.
\end{equation}
The point 
\[
X=i\coth(\lambda/2)(\sech(\tau/2)+i\tanh(\tau/2))
\]
is the intersection points of the axes of $\rho_{\ell,\tau}(a)$ and $\rho_{\ell,\tau}(b)$. Its signed distance to $H$ with respect to the orientation of the axis of $\rho_{\ell,\tau}(a)$ is $\tau/2$. When $\tau>0$, the real part of $X$ is negative and the signed distance from $X$ to $H$ is therefore positive, which is coherent with the conclusion of Proposition~\ref{prop:twist-coordinate-is-geometric-distance}. The image of $X$ by $\rho_{\ell,\tau}(b)$ is the point 
\[
X'=ie^\eta\coth(\lambda/2)(\sech(\tau/2)-i\tanh(\tau/2)).
\]
The geodesic line $(XX')$ is the axis of $\rho_{\ell,\tau}(b)$.

\section{Pentagon representations}\label{sec:pentagon-representations}
\subsection{Overview}
This section is dedicated to a special kind of representations with Euler number $\pm 1$ known as pentagon representations. We start by recalling some generalities about the hyperelliptic involution in genus $2$ (Section~\ref{sec:hyperelliptic-involution}), before defining pentagon representations and stating their main properties (Section~\ref{sec:definition-pentagon-representations}). We then provide a parametrization of the space of pentagon representation with Fenchel--Nielsen-type coordinates (Section~\ref{sec:parametrization-pentagon-representations}) and show that the coordinates are symplectic (Section~\ref{sec:symplectic-structure-pentagon-representations}). Finally, we explain how to geometrize pentagon representations using branched hyperbolic structures (Section~\ref{sec:geometrization-pentagon-representations}).

\subsection{Hyperelliptic involution}\label{sec:hyperelliptic-involution}
\begin{figure}[h]
\centering
\begin{tikzpicture}
\node[anchor = south west, inner sep=0mm] at (0,0) {\includegraphics[width=6.5cm]{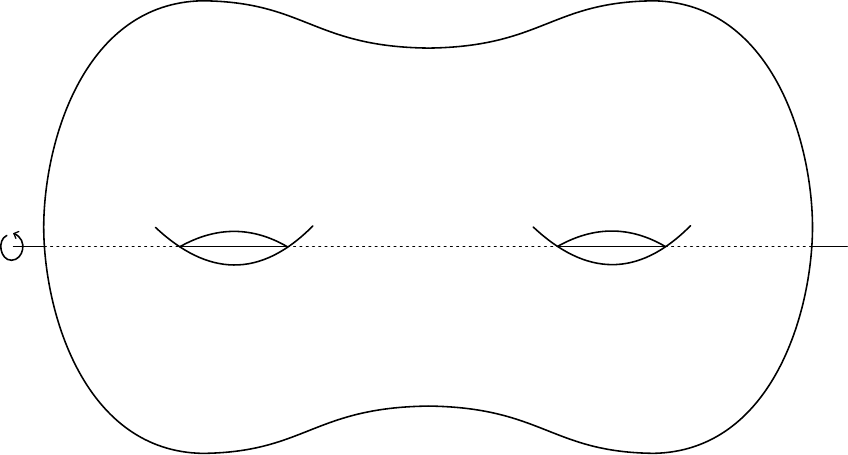}};

\node at (-.25,1.6) {$\varphi$};
\node[anchor = south west, inner sep=0mm] at (7,0) {\includegraphics[width=6.5cm]{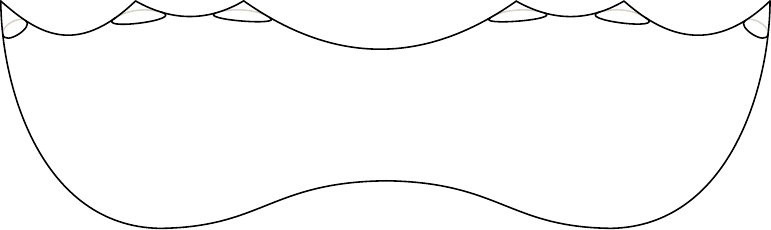}};
\end{tikzpicture}
\caption{Hyperelliptic involution in genus $2$ and the corresponding quotient orbifold.}
\label{fig:hyperelliptic-involution}
\end{figure}

The mapping class group $\Mod(S)$ of a closed and oriented surface $S$ of genus $g=2$ contains a remarkable element known as \emph{hyperelliptic involution}. It is an order $2$ element that can be defined as the isotopy class of a rigid rotation $\varphi$ of $S$ by an angle $\pi$ around the axis depicted on Figure~\ref{fig:hyperelliptic-involution}. The hyperelliptic involution is the only non-trivial central element of $\Mod(S)$ and it fixes every free isotopy class of unoriented simple closed curves on $S$ (Fact~\ref{fact:hyperelliptic-involution-fixes-scc}).

\begin{rmk}
In general, for surfaces of genus $g\geq 2$, hyperelliptic involutions are defined as mapping classes of order $2$ that act by minus the identity on the first homology of the surface with integer coefficients. Equivalently, they are isotopy classes of surface homeomorphisms of order $2$ with $2g+2$ fixed points. When $g\geq 3$, there are infinitely many hyperelliptic involutions, all conjugate to each other. They only fix some isotopy classes of simple closed curves, but not all of them. When $g=2$, the hyperelliptic involution is unique. We refer the reader to~\cite[Proposition~7.5]{mcg-primer} and the discussion thereafter for more details.
\end{rmk}

The involution $\varphi$ fixes six points on $S$ and the quotient $S/\varphi$ is a sphere with six orbifold points of degree $2$. The quotient map realizes $\pi_1S$ as an index $2$ subgroup of the orbifold fundamental group of $S/\varphi$, which we will denote by $\Gamma$. The group $\Gamma$ can be presented as follows:
\begin{equation}\label{eq:Gamma}
\Gamma = \langle q_1,\ldots, q_6: q_i^2=1,\, q_1\cdots q_6=1\rangle.
\end{equation}
The group $\pi_1S$ is classically presented as
\begin{equation}\label{eq:presentation-pi1S}
\pi_1S=\langle a_1,b_1,a_2,b_2 \, : \, [a_1,b_1][a_2,b_2]=1\rangle.
\end{equation}
An explicit injective morphism $\pi_1S\hookrightarrow \Gamma$ is given by
\begin{equation}\label{eq:injection-pi1S-Gamma}
\begin{cases} 
a_1\mapsto q_3q_2\\
b_1\mapsto q_1q_2\\
a_2\mapsto q_6q_5\\
b_2\mapsto q_4q_5.
\end{cases}
\end{equation}
In particular, $[a_1,b_1]=a_1b_1a_1^{-1}b_1^{-1}\mapsto (q_3q_2q_1)^2=(q_4q_5q_6)^2$. 

The hyperelliptic involution can be lifted as an automorphism of $\pi_1S$. A possible choice is the involution $\overline{\varphi}\in\Aut(\pi_1S)$ given by
\begin{equation}\label{eq:morphism-hyperelliptic-involution-pi1S}
\overline{\varphi}\colon\begin{cases}
a_1\mapsto a_1^{-1}\\
b_1\mapsto b_1^{-1}\\
a_2\mapsto\gamma b_2a_2^{-1}b_2^{-1}\gamma^{-1}\\
b_2\mapsto\gamma b_2^{-1}\gamma^{-1},
\end{cases}
\end{equation}
where $\gamma=a_1^{-1}b_1^{-1}a_2$ (see~\cite[Section~3]{marche-wolff}). Given a representation $\rho\colon\pi_1S\to\psl$, the image of its conjugacy class $[\rho]$ by the hyperelliptic involution is $[\rho\circ \overline{\varphi}]$.
\begin{fact}[\cite{hyperelliptic-involution}]\label{fact:hyperelliptic-involution-fixes-scc}
The hyperelliptic involution in genus two preserves the isotopy class of every unoriented simple closed curve on $S$. More precisely, it preserves the orientation of separating closed curves and reverses the orientation of non-separating simple closed curves.
\end{fact}

Depending on the Euler number of a representation $\rho\colon\pi_1S\to\psl$, the following may happen:
\begin{itemize}
    \item If $\eu(\rho)=0$, then $\rho$ lifts to a representation $\overline{\rho}\colon\pi_1S\to\SL_2\R$.  By Fact~\ref{fact:hyperelliptic-involution-fixes-scc} and since matrices in $\SL_2\R$ have the same trace as their inverses, the traces of the image by $\overline{\rho}$ and $\overline{\rho}\circ\overline{\varphi}$ of every simple closed curve on $S$ coincide. So, if $\rho$ has Zariski dense image, then $\overline{\rho}$ and $\overline{\rho}\circ\overline{\varphi}$ are conjugated by either an element of $\SL_2\R$, or by an element of $\GL_2\R$ with determinant $-1$. Marché--Wolff proved that these two eventualities lead to the only two ergodic components for the action of $\Mod(S)$ on the Euler number $0$ component of $\chi(S,\psl)$ (Theorem~\ref{thm:goldman-ergocicity-genus-2}).
    \item If $\eu(\rho)=\pm 2$, then $\rho$ is a Fuchsian representation. Since Fuchsian representations have Zariski dense images and $\Mod(S)$ preserves each connected component of $\chi(S,\psl)$, the previous argument shows that $[\rho]=[\rho\circ\overline{\varphi}]$ must hold for every Fuchsian representation $\rho$.
    \item If $\eu(\rho)=\pm 1$, then the fixed points locus of the hyperelliptic involution (which is a symplectic involution) is a codimension $2$ symplectic submanifold of the two components of $\chi(S,\psl)$ with Euler number $\pm 1$. We will describe the symplectic structure explicitly in Theorem~\ref{thm:wolpert-formula-pentagon-representations}. The fixed points coincide with the set of \emph{pentagon representations} introduced by Le Fils~\cite{thomas} (see Definition~\ref{defn:pentagon-representation} and Lemma~\ref{lem:pentagon-iff-fixed-by-hyperelliptic-involution}). In each component, we will see that the fixed point locus is made of six connected components, each of them symplectomorphic to the standard $\R^4$.
\end{itemize}

The automorphism $\overline{\varphi}\colon\pi_1S\to\pi_1S$ from~\eqref{eq:morphism-hyperelliptic-involution-pi1S} is the restriction of the inner automorphism $\widehat{\varphi}\colon\Gamma\to\Gamma$ defined by
\begin{equation}\label{eq:morphism-hyperelliptic-involution-Gamma}
\widehat{\varphi}(q_i)=q_2q_iq_2.
\end{equation}
In other words, $\widehat{\varphi}$ is the conjugation by $q_2$. This can be checked by direct computations.

\subsection{Definition and properties}\label{sec:definition-pentagon-representations}
Following~\cite[Definition~1.2]{thomas}, we define pentagon representations as follows.
\begin{defn}\label{defn:pentagon-representation}
A representation $\rho\colon\pi_1S\to\psl$ is a \emph{pentagon representation} if $\rho$ is the restriction of a representation $\overline{\rho}\colon\Gamma\to \psl$ whose kernel contains exactly one of the generators $q_i\in \Gamma$ introduced in~\eqref{eq:Gamma}. We will denote by $\mathcal{P}\subset \chi(S,\psl)$ the subset of conjugacy classes of pentagon representations.
\end{defn}

\begin{ex}\label{ex:pentagon-representation}
The first examples of pentagon representations (which justify the name) are constructed from right-angled hyperbolic pentagons. There are several ways to do so, we present one here and another one in Example~\ref{ex:pentagon-representation-P5}. Let $(C_1,\ldots,C_5)$ be the vertices of a right-angled pentagon and denote by $s_i$ the unique elliptic element of $\psl$ that fixes $C_i$ and rotates by $\pi$. The map
\begin{align*}
   \overline{\rho}\colon \Gamma &\to \psl\\
    q_i &\mapsto s_i, \quad i=1,\ldots,5\\
    q_6 &\mapsto 1
\end{align*}
is a group morphism because each $s_i$ has order $2$ and can be decomposed as the product of the two reflections through the adjacent sides of the pentagon meeting at $C_i$. The induced morphism $\rho\colon\pi_1S\to\psl$ given by~\eqref{eq:injection-pi1S-Gamma} is thus a pentagon representation according to Definition~\ref{defn:pentagon-representation}. Note that $\rho(a_1)=s_3s_2$ is a hyperbolic translation along the geodesic line $(C_2C_3)$\footnote{When we say hyperbolic translation along $(C_2C_3)$, we mean that $C_2$ is moved towards $C_3$.} of length $2d(C_2,C_3)$. Similarly, as illustrated by Figure~\ref{fig:pentagon-representation-P6}, $\rho(b_1)$ and $\rho(b_2)$ are also hyperbolic translations along a side of the pentagon of length twice the corresponding side length, while $\rho(a_2)$ is an elliptic rotation of angle $\pi$ around $C_5$.
\begin{figure}[h]
\centering
\begin{tikzpicture}
\node[anchor = south west, inner sep=0mm] at (0,0) {\includegraphics[width=7cm]{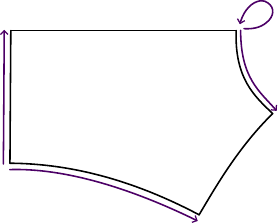}};

\node[plum] at (-.45,3.3) {$\rho(b_1)$};
\node[plum] at (2.6,.6) {$\rho(a_1)$};
\node[plum] at (6.85,4) {$\rho(b_2)$};
\node[plum] at (7.45,5.4) {$\rho(a_2)$};

\node at (.1,5.1) {$C_1$};
\node at (0.05,1.25) {$C_2$};
\node at (5.25,0) {$C_3$};
\node at (7.2,2.7) {$C_4$};
\node at (5.75,5.1) {$C_5$};
\end{tikzpicture}
\caption{First example of pentagon representation associated to a right-angled pentagon.}
\label{fig:pentagon-representation-P6}
\end{figure}
\end{ex}

We compile some important properties of pentagon representations in Lemma~\ref{lem:properties-pentagon}. They can either be found as such in Le Fils' work~\cite{thomas}, or are direct consequences of results therein.

\begin{lem}\label{lem:properties-pentagon}
Pentagon representations satisfy the following:
\begin{enumerate}
    \item Pentagon representations have Euler number $\pm1$.\label{cond:eu=pm1}
    \item The set of pentagon representations $\mathcal{P}$ is invariant under the action of $\Mod(S)$.\label{cond:mod-invariance}
    \item Pentagon representations map simple closed curves to either hyperbolic elements, or to elliptic elements of rotation angle $\pi$. Furthermore, all separating simple closed curves are mapped to hyperbolic elements.\label{cond:image-of-scc}
\end{enumerate}
\end{lem}
\begin{proof}
Condition~\eqref{cond:eu=pm1} is equivalent to saying that pentagon representations do not lift to representations into $\SL_2\R$. A proof is given in~\cite[Proposition~3.1]{thomas}. 

Mapping class group invariance of $\mathcal{P}$ (Condition~\eqref{cond:mod-invariance}) can be proved using a specific generating family of $\Mod(S)$ consisting of five Dehn twists such as those described in Appendix~\ref{apx:generators}. For each Dehn twist, Le Fils gives an explicit automorphism of $\Gamma$ which restricts to an automorphism of $\pi_1S$ that represents the Dehn twist (see~\eqref{eq:generators-Mos(S)}) and concludes mapping class group invariance~\cite[Proposition~3.2]{thomas}.

Since there are only one separating and one non-separating simple closed curve on $S$ up to mapping class group action, in order to prove Condition~\eqref{cond:image-of-scc}, it is enough to prove that if $\rho$ is a pentagon representation, then $\rho(a_1)$ is hyperbolic or elliptic of angle $\pi$, and $\rho([a_1,b_1])$ is hyperbolic. Using that $\rho$ is the restriction of some $\overline{\rho}\colon\Gamma\to\psl$ and recalling~\eqref{eq:injection-pi1S-Gamma}, we have $\rho(a_1)=\overline{\rho}(q_3)\overline{\rho}(q_2)$. Both $\overline{\rho}(q_3)$ and $\overline{\rho}(q_2)$ are order $2$ elements of $\psl$ by the definition of $\Gamma$, see~\eqref{eq:Gamma}. Only one of them can be trivial according to Definition~\ref{defn:pentagon-representation}. If they are both non-trivial, then they are elliptic rotations of angle $\pi$. The product of two elliptic rotations of angle $\pi$ is trivial only when they have the same fixed point. In that case, we would have $\rho(a_1)=1$, which contradicts~\eqref{cond:eu=pm1}. If the fixed points are distinct, then the product is a hyperbolic translation along the axis passing through the two fixed points. This shows that $\rho(a_1)$ is hyperbolic or elliptic of angle $\pi$, as desired. Similarly, we recall that $\rho([a_1,b_1])=\overline{\rho}(q_3q_2q_1)^2=\overline{\rho}(q_4q_5q_6)^2$. We know that $\overline{\rho}(q_i)=1$ for a unique index $i$. Let us assume that $i=1$ (the other cases can be treated in the same way). This means that $\rho([a_1,b_1])=(\overline{\rho}(q_3)\overline{\rho}(q_2))^2$. The above discussion shows that $\overline{\rho}(q_3)$ and $\overline{\rho}(q_2)$ are elliptic rotations of angle $\pi$ around distinct fixed points, implying that $\rho([a_1,b_1])$ is hyperbolic.
\end{proof}

The next result relates pentagon representations with the fixed points of the hyperelliptic involution. We learned about it in a discussion with Thomas Le Fils.

\begin{lem}\label{lem:pentagon-iff-fixed-by-hyperelliptic-involution}
Let $\rho\colon\pi_1S\to\psl$ be a \emph{regular} representation, which means that the centralizer of $\rho(\pi_1S)$ in $\psl$ is trivial. Then the following are equivalent.
\begin{enumerate}
    \item The conjugacy class $[\rho]$ is fixed by the hyperelliptic involution.\label{cond:rho-fixed-point}
    \item The representation $\rho$ is the restriction of a representation $\overline{\rho}\colon\Gamma\to\psl$.\label{cond:rho-hyperelliptic}
\end{enumerate}
Moreover, if~\eqref{cond:rho-fixed-point} and~\eqref{cond:rho-hyperelliptic} hold and $n$ denotes the number of generators $q_i$ in the kernel of $\overline{\rho}$, then one of the following holds.
\begin{itemize}
    \item $n=0$ and $\rho$ is either Fuchsian, or $\eu(\rho)=0$. 
    \item $n=1$ and $\rho$ is a pentagon representation.
    \item $n\geq 2$ and $\rho(\pi_1S)$ stabilizes a geodesic line in the hyperbolic plane (equivalently, after conjugation, $\rho(\pi_1S)$ is contained in the \emph{infinite dihedral subgroup} of $\psl$---the subgroup of diagonal and anti-diagonal elements). In particular, $\eu(\rho)=0$.
\end{itemize}
\end{lem}
\begin{proof}
We will use the automorphisms $\varphi\colon\pi_1S\to\pi_1S$ and $\widehat{\varphi}\colon\Gamma\to\Gamma$ representing the hyperelliptic involution which we introduced in~\eqref{eq:morphism-hyperelliptic-involution-pi1S} and~\eqref{eq:morphism-hyperelliptic-involution-Gamma}

To prove the implication $\eqref{cond:rho-hyperelliptic}\Rightarrow \eqref{cond:rho-fixed-point}$, let $g=\overline{\rho}(q_2)$ and note that, for every $x\in\pi_1S$, 
\[
\rho\circ\overline{\varphi}(x)=\overline{\rho}\circ\widehat{\varphi}(x)=\overline{\rho}(q_2xq_2)=g\overline{\rho}(x)g^{-1}=g\rho(x)g^{-1}.
\]
This implication does not require the regularity assumption on $\rho$.

Conversely, to prove $\eqref{cond:rho-fixed-point}\Rightarrow \eqref{cond:rho-hyperelliptic}$, assume that $\rho\circ\overline{\varphi}=g\rho g^{-1}$ for some $g\in\psl$. We will prove that a possible extension $\overline{\rho}\colon\Gamma\to\psl$ of $\rho$ is given by
\begin{align*}
    &\overline{\rho}(q_1)=g\rho(b_1^{-1}),\quad \overline{\rho}(q_2)=g,\quad \overline{\rho}(q_3)=\rho(a_1)g,\\
    &\overline{\rho}(q_4)=\rho(a_2^{-1}b_2^{-1})g\rho(a_1^{-1}b_1^{-1}b_2^{-1}),\quad \overline{\rho}(q_5)=\overline{\rho}(q_4)\rho(b_2),\quad \overline{\rho}(q_6)=\rho(a_2)\overline{\rho}(q_5).
\end{align*}
First, we verify that $\overline{\rho}$ is a morphism. Note that $\overline{\rho}(q_1^2)=g\rho(b_1^{-1})g\cdot\rho(b_1^{-1})=\rho(\varphi(b_1^{-1}))\rho(b_1)=\rho(b_1)\rho(b_1^{-1})=1$. Similarly, $\overline{\rho}(q_2^2)=g^2$. Since $\overline{\varphi}$ is an involution, $g^2$ lies in the centralizer of the image of $\rho$ and is therefore trivial because $\rho$ is regular. Thus, $\overline{\rho}(q_2^2)=1$. A slightly more involved computation shows that
\begin{align*}
\overline{\rho}(q_4^2)&=\rho(a_2^{-1}b_2^{-1})g\rho(a_1^{-1}b_1^{-1}b_2^{-1}a_2^{-1}b_2^{-1})g\rho(a_1^{-1}b_1^{-1}b_2^{-1})\\
&=\rho(a_2^{-1}b_2^{-1}a_1b_1\cdot\gamma\cdot b_2b_2a_2\underbrace{b_2^{-1}b_2}_{=1}\cdot\gamma^{-1}\cdot a_1^{-1}b_1^{-1}b_2^{-1})\\
&=\rho(a_2^{-1}b_2^{-1}a_1b_1\cdot a_1^{-1}b_1^{-1}a_2\cdot b_2\underbrace{b_2a_2\cdot a_2^{-1}b_1a_1\cdot a_1^{-1}b_1^{-1}b_2^{-1}}_{=1})\\
&=\rho(a_2^{-1}b_2^{-1}a_1b_1\cdot a_1^{-1}b_1^{-1}a_2\cdot b_2).
\end{align*}
The last term vanishes because of the group relation for $\pi_1S$, proving that $\overline{\rho}(q_4)^2=1$. Analogous computations show that $\overline{\rho}(q_i^2)=1$ for every other $i$. We also observe that
\begin{align*}
\overline{\rho}(q_1\cdots q_6)&=\underbrace{g\rho(b_1^{-1})\cdot g}_{=\rho(b_1)}\cdot \rho(a_1)g\cdot \underbrace{\overline{\rho}(q_4)\cdot\overline{\rho}(q_4)}_{=1}\rho(b_2)\cdot \rho(a_2)\underbrace{\overline{\rho}(q_5)}_{\overline{\rho}(q_4)\rho(b_2)}\\
&=\rho(b_1a_1)g\rho(b_2a_2)\overline{\rho}(q_4)\rho(b_2)\\
&=1,
\end{align*}
from which we conclude that $\overline{\rho}$ is indeed a morphism. It is immediate that $\overline{\rho}$ restricts to $\rho$ on $\pi_1S$, which concludes the proof of $\eqref{cond:rho-fixed-point}\Rightarrow \eqref{cond:rho-hyperelliptic}$.

Now, if $n=0$, then $\rho$ lifts to a representation into $\SL_2\R$ and therefore is either Fuchsian or has $\eu(\rho)=0$. Indeed, start by choosing an arbitrary lift $Q_i\in\SL_2\R$ of every $\overline{\rho}(q_i)\in\psl$. Note that $Q_i^{-1}=-Q_i$. The map $\pi_1S\to\SL_2\R$ defined as in~\eqref{eq:injection-pi1S-Gamma} by mapping $a_1$ to $Q_3Q_2$ and so on is a group homomorphism that lifts $\rho$ because it maps $[a_1,b_1][a_2,b_2]$ to $(Q_3Q_2Q_1)(Q_3Q_2Q_1Q_6Q_5Q_4)(Q_6Q_5Q_4)$ and $Q_3Q_2Q_1Q_6Q_5Q_4=\pm\id$ as it projects to the identity in $\psl$. If $n=1$, then $\rho$ is a pentagon representation by Definition~\ref{defn:pentagon-representation}. If $n\geq 2$, then the non-trivial images of the generators $q_i$ are elliptic rotations of angle $\pi$ whose product is trivial. Since there are at most four such rotations, the centers must lie on a common geodesic line in order for the product to be trivial. This means that $\overline{\rho}(\Gamma)$ stabilizes a geodesic line in the hyperbolic plane and thus $\eu(\rho)=0$.
\end{proof}

\subsection{Coordinates}\label{sec:parametrization-pentagon-representations}
Recall from Definition~\ref{defn:pentagon-representation} that we denoted by $\mathcal{P}$ the subset of the character variety $\chi(S,\psl)$ consisting of all pentagon representations. Lemma~\ref{lem:properties-pentagon} says that $\mathcal{P}$ lies in the two connected components of $\chi(S,\psl)$ corresponding to representations with Euler number $\pm 1$. Moreover, by Lemma~\ref{lem:pentagon-iff-fixed-by-hyperelliptic-involution}, $\mathcal{P}$ corresponds to the fixed point set of the hyperelliptic involution, and is therefore a closed embedded symplectic submanifold of those two components. We will see in Theorem~\ref{thm:coordinates-pentagon-representations} that $\mathcal{P}$ has codimension $2$. We will write 
\[
\mathcal{P}=\mathcal{P}^+\sqcup \mathcal{P}^-,
\]
where $\mathcal{P}^+$ and $\mathcal{P}^-$ denotes the elements of $\mathcal{P}$ with Euler number $+1$, respectively $-1$. Note that $\mathcal{P}^+$ and $\mathcal{P}^-$ are images of each other under the symplectic involution of $\chi(S,\psl)$ given by conjugation by any matrix of determinant $-1$. Lemma~\ref{lem:properties-pentagon} also implies that $\mathcal{P}$ decomposes as the disjoint union of six subsets 
\[
\mathcal{P}=\mathcal{P}_1\sqcup\cdots\sqcup \mathcal{P}_6,
\]
where each $\mathcal{P}_i$ consists of the conjugacy classes of the restrictions to $\pi_1S$ of all the representations $\Gamma\to\psl$ with $q_i\mapsto 1$. For instance, the conjugacy class of the pentagon representation from Example~\ref{ex:pentagon-representation} belongs to $\mathcal{P}_6$. Note that $\mathcal{P}_i$ and $\mathcal{P}_j$ are indeed disjoint when $i\neq j$, as there is one curve among $\{a_1,b_1,a_2,b_2\}$ whose image has different types in $\mathcal{P}_i$ and $\mathcal{P}_j$ (hyperbolic versus elliptic, see Table~\ref{tab:6-components-pentagons}) and so there cannot be a continuous path of pentagon representations connecting $\mathcal{P}_i$ and $\mathcal{P}_j$. We will prove below (Theorem~\ref{thm:coordinates-pentagon-representations}) that the intersections $\mathcal{P}^+_i=\mathcal{P}^+\cap \mathcal{P}_i$ and $\mathcal{P}^-_i=\mathcal{P}^-\cap \mathcal{P}_i$ are the connected components of the space $\mathcal{P}$.
\begin{table}[h]
    \centering
    \begin{tabular}{cccccccc}
        \toprule
        $[\rho]$ & $\mathcal{P}_1$ & $\mathcal{P}_2$ & $\mathcal{P}_3$ & $\mathcal{P}_4$ & $\mathcal{P}_5$ & $\mathcal{P}_6$  \\ \midrule 
        $a_1$ & hyperbolic & elliptic & elliptic & hyperbolic & hyperbolic & hyperbolic\\ 
        $b_1$ & elliptic & elliptic & hyperbolic & hyperbolic & hyperbolic & hyperbolic\\ 
        $a_2$ & hyperbolic & hyperbolic & hyperbolic & hyperbolic & elliptic & elliptic\\ 
        $b_2$ & hyperbolic & hyperbolic & hyperbolic & elliptic & elliptic & hyperbolic\\
        \bottomrule
    \end{tabular}
    \caption{The six different types of pentagon representations according to~\eqref{eq:injection-pi1S-Gamma} and the image of the generators of $\pi_1S$ from~\eqref{eq:presentation-pi1S} in every case.}
    \label{tab:6-components-pentagons}
\end{table}

Observe that the Dehn twist $\tau_{b_1}$ along the curve $b_1$ interchanges $\mathcal{P}_1$ and $\mathcal{P}_2$, by which we mean that $\tau_{b_1}(\mathcal{P}_1)=\mathcal{P}_2$ and $\tau_{b_1}(\mathcal{P}_2)=\mathcal{P}_1$, and also $\tau_{b_1}^2$ is trivial on both $\mathcal{P}_1$ and $\mathcal{P}_2$. This can be seen using the explicit formulae~\eqref{eq:generators-Mos(S)} from Appendix~\ref{apx:generators}. Similarly, the Dehn twist along the curve $a_1$ interchanges $\mathcal{P}_2$ and $\mathcal{P}_3$, and the Dehn twist along the curve $b_2$ interchanges $\mathcal{P}_4$ and $\mathcal{P}_5$, while the Dehn twist along the curve $a_2$ interchanges $\mathcal{P}_5$ and $\mathcal{P}_6$.

To fully understand the symplectic geometry of $\mathcal{P}$, it is therefore sufficient to understand the individual pieces $\mathcal{P}_i$, because they are symplectomorphic to each other. Our parametrization of $\mathcal{P}$ consists of two lengths and two twists, analogously to Fenchel--Nielsen coordinates. Note that our coordinates depend on the choice of geometric presentation of $\pi_1S$ made in~\eqref{eq:presentation-pi1S}. This choice plays the role of the topological data that traditional Fenchel--Nielsen coordinates depend on (see e.g.~\cite[Section 10.6.1]{mcg-primer}). 

We shall focus on $\mathcal{P}_5$ and start by revisiting the construction of Example~\ref{ex:pentagon-representation}. 
\begin{ex}\label{ex:pentagon-representation-P5}
Given a right-angled pentagon $ABCDE$, we denote by $b/2$ the distance $d(A,B)$, by $c/2$ the distance $d(B,C)$, and so on. The pentagon representation $\rho$ from Example~\ref{ex:pentagon-representation} was defined as follows:\footnote{Recall that the hyperbolic translation along $(BC)$ means that $B$ is moved towards $C$.}
\begin{align*}
    \rho(a_1)&=\text{hyperbolic translation along $(BC)$
    of length $c$}\\
    \rho(b_1)&=\text{hyperbolic translation along $(BA)$ of length $b$}\\
    \rho(a_2)&=\text{elliptic rotation of angle $\pi$ around $E$}\\
    \rho(b_2)&=\text{hyperbolic translation along $(ED)$ of length $e$}.
\end{align*}
Its conjugacy class belongs to $\mathcal{P}_6$. If we apply the Dehn twist $\tau_{a_2}$ along the curve $a_2$ to $[\rho]\in\mathcal{P}$, we get a point in $\mathcal{P}_5$ which is the conjugacy class of the representation $\rho'$ illustrated on Figure~\ref{fig:pentagon-representation-P5} and given by
\begin{align}
    \rho'(a_1)&=\rho(a_1)=\text{hyperbolic translation along $(BC)$ of length $c$}\label{eq:representation-right-angled-pentagon}\\
    \rho'(b_1)&=\rho(b_1)=\text{hyperbolic translation along $(BA)$ of length $b$}\nonumber \\
    \rho'(a_2)&=\rho(a_2)=\text{elliptic rotation of angle $\pi$ around $E$}\nonumber \\
    \rho'(b_2)&=\rho(b_2a_2^{-1})=\text{elliptic rotation of angle $\pi$ around $D$}.\nonumber
\end{align}
Note that $\rho'([a_1,b_1])=\rho'([b_2,a_2])=\rho(b_2)^2$ is a hyperbolic translation along $(ED)$ of length $2e$ (see Fact~\ref{fact_virtuallyabelian_char} for details).
\begin{figure}[h]
\centering
\begin{tikzpicture}
\node[anchor = south west, inner sep=0mm] at (0,0) {\includegraphics[width=8cm]{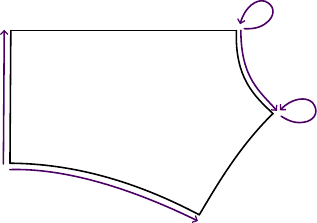}};

\node[plum] at (-.45,3.3) {$\rho'(b_1)$};
\node[plum] at (2.6,.6) {$\rho'(a_1)$};
\node[plum] at (7.25,4) {$\rho'([a_1,b_1])$};
\node[plum] at (7.45,5.4) {$\rho'(a_2)$};
\node[plum] at (8.5,2.9) {$\rho'(b_2)$};

\node at (.1,5.1) {$A$};
\node at (0.05,1.25) {$B$};
\node at (5.25,0) {$C$};
\node at (6.95,2.5) {$D$};
\node at (5.75,5.1) {$E$};

\node at (.65,3.3) {$b/2$};
\node at (2.9,4.55) {$a/2$};
\node at (5.7,3.8) {$e/2$};
\node at (5.55,1.7) {$d/2$};
\node at (2.8,1.45) {$c/2$};
\end{tikzpicture}
\caption{The pentagon representation $\rho'$ described in~\eqref{eq:representation-right-angled-pentagon}.}
\label{fig:pentagon-representation-P5}
\end{figure}
\end{ex}

Example~\ref{ex:pentagon-representation-P5} shows that to every right-angled pentagon corresponds a point in $\mathcal{P}_5$. This gives us a first $2$-dimensional slice in $\mathcal{P}_5$, because the moduli space of right-angled pentagons in the hyperbolic plane has dimension $2$. More precisely, the following holds.
\begin{fact}[{e.g.~\cite[Theorem~7.18.1]{beardon}}]\label{fact:right-angled-pentagons}
For every $c>0$ and $e>0$, there is a unique right-angled pentagon, up to isometries of the hyperbolic plane, with vertices $(A,B,C,D,E)$ and side lengths $a/2,b/2,c/2,d/2,e/2$ as on Figure~\ref{fig:pentagon-representation-P5}. Furthermore, the following formulae hold:
\[
\sinh(b/2)=\frac{\cosh(e/2)}{\sinh(c/2)}\quad\text{and}\quad \cosh(d/2)=\frac{1}{\tanh(c/2)\tanh(e/2)}.
\]
\end{fact}
\begin{rmk}\label{rem:euler-number-pentagon-representation}
The uniqueness statement in Fact~\ref{fact:right-angled-pentagons} can be refined. Two lengths $c>0$ and $e>0$ determine a unique right-angled pentagon up to isometry and two pentagons up to orientation-preserving isometries (so up to $\psl$). To each of these two pentagons one can associate the conjugacy class of the pentagon representation~\eqref{eq:representation-right-angled-pentagon} from Example~\ref{ex:pentagon-representation-P5}. Denote by $[\rho_{c,e}^+]$ the one associated with the pentagon whose vertices appear in clockwise order, and by $[\rho_{c,e}^-]$ the other one. As the notation suggests,
\[
\eu(\rho_{c,e}^+)=1\quad \text{and}\quad \eu(\rho_{c,e}^-)=-1,
\]
which means $[\rho_{c,e}^{+}]\in\mathcal{P}_5^+$ and $[\rho_{c,e}^{-}]\in\mathcal{P}_5^-$.
\end{rmk}

Our plan to parametrize all of $\mathcal{P}_5$ is to describe two directions of deformations from the $2$-dimensional slice corresponding to right-angled pentagons and then show that $\mathcal{P}_5$ is indeed $4$-dimensional (Theorem~\ref{thm:coordinates-pentagon-representations}). We will do so using the Hamiltonian flows of the length functions associated to the curves $a_1$ and $[a_1,b_1]$, which we also call \emph{twist flows} as in Section~\ref{sec:Fenchel-Nielsen-coordinates-one-holed-torus}.  More precisely, consider the open subset $U$ of the two components of $\chi(S,\psl)$ with Euler number $\pm 1$ defined by
\[
U=\left\{[\rho]\in\eu^{-1}(\pm 1)\,:\, \rho(a_1)\text{ and }\rho([a_1,b_1])\text{ are hyperbolic}\right\}.
\]
The subset $U$ is open because the set of hyperbolic elements in $\psl$ is open. Moreover, $\mathcal{P}_5\subset U$ by the definition of $\mathcal{P}_5$ and by Lemma~\ref{lem:properties-pentagon}. The length functions of $a_1$ and $[a_1,b_1]$ are smooth functions defined on $U$. Namely, if $x$ denotes either $a_1$ or $[a_1,b_1]$, then the \emph{length function} of $x$ is
\begin{align*}
    \ell_{x}\colon U&\to \R_{>0}\\
    [\rho]&\mapsto \text{translation length of }\rho(x).
\end{align*}
As in Section~\ref{sec:Fenchel-Nielsen-coordinates-one-holed-torus}, we denote the Hamiltonian flow of $\ell_x$ for the Goldman symplectic form restricted to $U$ by
\[
\Phi_x^t\colon U\to U, \quad t\in\R.
\]
The two flows are the \emph{twist flows} of $a_1$ and $[a_1,b_1]$. Since the simple curves $a_1$ and $[a_1,b_1]$ have disjoint representatives, the two twist flows commute.

The twist flows descend from explicit flows on representations, as we explained before~\eqref{eq:twist-flow-one-holed-torus}. Recall that for $t\in \R$, we denoted by $X(t)$ the Lie algebra element
\[
X(t)=\begin{pmatrix}
    t/2 & 0\\ 0 & -t/2
\end{pmatrix}.
\]
Furthermore, for every $[\rho]\in U$, there exist $A_\rho$ and $B_\rho$ in $\psl$ such that
\[
\rho(a_1)=A_\rho\exp\big(X\big(\ell_{a_1}(\rho)\big)\big)A_\rho^{-1} \quad\text{and}\quad \rho([a_1,b_1])=B_\rho\exp\big( X\big(\ell_{[a_1,b_1]}(\rho)\big)\big)B_\rho^{-1}.
\]
Using the shorthand notation
\[
\xi_\rho(t)=A_\rho\exp\big(X(t)\big)A_\rho^{-1} \quad\text{and}\quad \zeta_\rho(t)=B_\rho\exp\big(X(t)\big)B_\rho^{-1},
\]
we note that $\{\xi_\rho(t):t\in\R\}$ and $\{\zeta_\rho(t):t\in\R\}$ are the respective centralizers of $\rho(a_1)$ and $\rho([a_1,b_1])$ in $\psl$. Goldman proved in~\cite{goldman-invariant-functions} that the twist flows $\Phi_{a_1}^t\big([\rho])$ and $\Phi_{[a_1,b_1]}^t\big([\rho])$ can be lifted to the following flows on the space of representations:
\begin{equation}\label{eq:twist-flows-genus-2}
\overline{\Phi}_{a_1}^t(\rho):\begin{cases}
    a_1\mapsto \rho(a_1)\\
    b_1\mapsto \rho(b_1)\xi_\rho(t)^{-1}\\
    a_2\mapsto \rho(a_2)\\
    b_2\mapsto \rho(b_2)
\end{cases}\quad\text{and}\quad 
\overline{\Phi}_{[a_1,b_1]}^t(\rho):\begin{cases}
    a_1\mapsto \rho(a_1)\\
    b_1\mapsto \rho(b_1)\\
    a_2\mapsto \zeta_\rho(t)\rho(a_2)\zeta_\rho(t)^{-1}\\
    b_2\mapsto \zeta_\rho(t)\rho(b_2)\zeta_\rho(t)^{-1}.
\end{cases}
\end{equation}

\begin{lem}\label{lem:twist-flows-preserves-P5}
The twist flows of $a_1$ and $[a_1,b_1]$ preserve $\mathcal{P}_5$. In other words, if $[\rho]\in\mathcal{P}_5$, then $\Phi_{a_1}^t\big([\rho])\in \mathcal{P}_5$ and $\Phi_{[a_1,b_1]}^t\big([\rho])\in \mathcal{P}_5$ for every $t\in\R$.
\end{lem}
\begin{proof}
Since $\mathcal{P}_5$ is disjoint from the other $\mathcal{P}_i$, it is sufficient to prove that $\Phi_{a_1}^t\big([\rho])\in \mathcal{P}$ and $\Phi_{[a_1,b_1]}^t\big([\rho])\in \mathcal{P}$ for every $t\in\R$. As we are assuming $[\rho]\in\mathcal{P}_5$, it follows from Lemma~\ref{lem:pentagon-iff-fixed-by-hyperelliptic-involution} that $[\rho]$ is fixed by the hyperelliptic involution. So, there exists some $T\in\psl$ such that $\rho=T(\rho\circ\overline{\varphi})T^{-1}$, where $\overline{\varphi}\colon\pi_1S\to\pi_1S$ is the automorphism given in~\eqref{eq:morphism-hyperelliptic-involution-pi1S} that represents the hyperelliptic involution. In particular, it holds that $\rho(a_1)=T\rho(a_1^{-1})T^{-1}$ by Fact~\ref{fact:hyperelliptic-involution-fixes-scc}. Since $[\rho]\in\mathcal{P}_5$, $\rho(a_1)$ is hyperbolic (see Table~\ref{tab:6-components-pentagons}) and $T$ must be an elliptic rotation of angle $\pi$ with center on the axis of $\rho(a_1)$. Thus, $T^{-1}=T$ and we will simply write $T$ instead of $T^{-1}$ below.

To prove $\Phi_{a_1}^t\big([\rho])\in \mathcal{P}$, we will use Lemma~\ref{lem:pentagon-iff-fixed-by-hyperelliptic-involution} and verify that $\Phi_{a_1}^t\big([\rho]\big)$ is fixed by the hyperelliptic involution. Recall that $\Phi_{a_1}^t\big([\rho]\big)$ can be lifted to $\overline{\Phi}_{a_1}^t\big(\rho\big)$ in the space of representations. We compute
\[
\left(\overline{\Phi}_{a_1}^t\big(\rho\big)\circ\overline{\varphi}\right)(a_1)=\overline{\Phi}_{a_1}^t\big(\rho\big)(a_1^{-1})=\rho(a_1^{-1})=T\rho(a_1)T.
\]
Since $\xi_\rho(t)$ is a hyperbolic translation along the same axis as $\rho(a_1)$, it commutes with $\rho(a_1)$ and it also holds that $T\xi_\rho(t)T=\xi_\rho(t)^{-1}$. Hence,
\[
T\rho(a_1)T=T\xi_\rho(t)^{-1}\cdot \rho(a_1)\cdot \xi_\rho(t)T=\xi_\rho(t)T\cdot \rho(a_1)\cdot T\xi_\rho(t)^{-1}.
\]
Similar computations give
\[
\left(\overline{\Phi}_{a_1}^t\big(\rho\big)\circ\overline{\varphi}\right)(b_1)=\overline{\Phi}_{a_1}^t\big(\rho\big)(b_1^{-1})=\xi_\rho(t)\rho(b_1^{-1})=\xi_\rho(t)T\rho(b_1)T,
\]
which further leads to
\[
\xi_\rho(t)T\rho(b_1)T=\xi_\rho(t)T\cdot \rho(b_1)\xi_\rho(t)^{-1}\cdot \xi_\rho(t)T=\xi_\rho(t)T\cdot \rho(b_1)\xi_\rho(t)^{-1}\cdot T\xi_\rho(t)^{-1}.
\]
Recalling that $\overline{\varphi}(a_2)=\gamma b_2a_2^{-1}b_2^{-1}\gamma^{-1}$ with $\gamma=a_1^{-1}b_1^{-1}a_2$, and using that $\xi_\rho(t)$ commutes with $\rho(a_1)$, we obtain
\[
\left(\overline{\Phi}_{a_1}^t\big(\rho\big)\circ\overline{\varphi}\right)(a_2)=\xi_\rho(t)(\rho\circ\overline{\varphi})(a_2)\xi_\rho(t)^{-1}=\xi_\rho(t)T\cdot \rho(a_2)\cdot T\xi_\rho(t)^{-1}.
\]
Analogous computations show that
\[
\left(\overline{\Phi}_{a_1}^t\big(\rho\big)\circ\overline{\varphi}\right)(b_2)=\xi_\rho(t)T\cdot \rho(b_2)\cdot T\xi_\rho(t)^{-1}.
\]
Consequently, $\overline{\Phi}_{a_1}^t\big(\rho\big)\circ\overline{\varphi}$ is conjugate to $\overline{\Phi}_{a_1}^t\big(\rho\big)$ by $\xi_\rho(t)T$, implying that $\Phi_{a_1}^t\big([\rho]\big)$ is fixed by the hyperelliptic involution and thus $\Phi_{a_1}^t\big([\rho]\big)\in\mathcal{P}$. Similar arguments also show that $\Phi_{[a_1,b_1]}^t\big([\rho]\big)\in\mathcal{P}$.
\end{proof}

\begin{figure}[h]
\centering
\begin{tikzpicture}
\node[anchor = south west, inner sep=0mm] at (0,0) {\includegraphics[width=10cm]{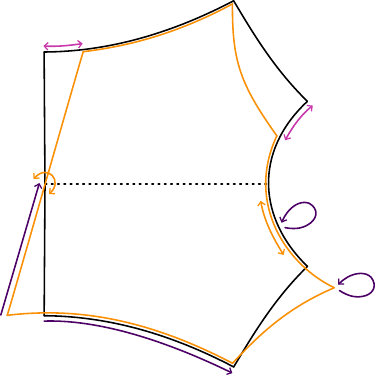}};

\node at (.6,5.1) {$A$};
\node at (0.9,1.4) {$B$};
\node at (0.95,8.75) {$B'$};
\node at (6.5,.25) {$C$};
\node at (6.5,10.1) {$C'$};
\node at (8.45,2.95) {$D$};
\node at (8.46,7.45) {$D'$};
\node at (7.4,5.1) {$E$};

\node[plum] at (-.2,3.3) {$\rho_{s,t}(b_1)$};
\node[plum] at (3.3,.7) {$\rho_{s,t}(a_1)$};
\node[plum] at (9.15,4.5) {$\rho_{s,t}(a_2)$};
\node[plum] at (9.6,3.0) {$\rho_{s,t}(b_2)$};

\node[mauve] at (8.1,6.6) {$t$};
\node[mauve] at (1.7,9.05) {$s/2$};

\node[apricot] at (1.85,5.5) {\Large $\Phi_{a_1}^t$};
\node[apricot] at (6.5,3.9) {\Large $\Phi_{[a_1,b_1]}^t$};
\end{tikzpicture}
\caption{Deformations of the right-angled hexagon in black by flow twists into the orange hexagon. The purple arrows indicate the images of the fundamental group generators under the resulting representation $\rho_{s,t}=\overline{\Phi}_{a_1}^s\circ\overline{\Phi}^t_{[a_1,b_1]}(\rho)$.}
\label{fig:pentagon-representation-reflected}
\end{figure}

The twist flows $\Phi_{a_1}^s$ and $\Phi_{[a_1,b_1]}^t$ have a particularly intuitive description when applied to a representation $[\rho]\in\mathcal{P}_5$ coming from a right-angled pentagon $ABCDE$, as in Example~\ref{ex:pentagon-representation-P5}. We first reflect the pentagon across the side $[AE]$ to obtain the right-angled hexagon $BCDD'C'B'$ illustrated on Figure~\ref{fig:pentagon-representation-reflected}, where $E$ is the midpoint of $[DD']$ and $A$ is the midpoint of $[BB']$. If we restrict $\rho$ to the free subgroup of $
\pi_1S$ generated by $a_1$ and $b_1$---the fundamental group of a one-holed torus $\Sigma$ obtained by cutting $S$ along the simple closed curve $[a_1,b_1]$---then the hexagon $BCDD'C'B'$ is the polygonal model of the restriction of $\rho$ to $\pi_1\Sigma$ constructed in Section~\ref{sec:polygonal-model-one-holed-tori}. This means that the twist flow $\overline{\Phi}_{a_1}^s$ acts on $\rho$ by rotating the axis of $\rho(b_1)$ around $A$, so that the axis point on the geodesic line $(BC)$---the axis of $\rho(a_1)$---moves by the signed distance $s/2$ (Proposition~\ref{prop:twist-coordinate-is-geometric-distance}). Similarly, the explicit form~\eqref{eq:twist-flows-genus-2} for the flow $\overline{\Phi}_{[a_1,b_1]}^t$ shows that it acts by translating the fixed points of $\rho(a_2)$ and $\rho(b_2)$ by a signed distance $t$ along the axis of $\rho([a_1,b_1])$---the geodesic line $(DE)$. Both deformations can be visualized on Figure~\ref{fig:pentagon-representation-reflected}.

This visual description of the twist flow $\Phi_{[a_1,b_1]}^t$ motivates the following notion of twist. Given $[\rho]\in\mathcal{P}_5$ (not necessarily coming from a right-angled pentagon), we consider the restriction of $[\rho]\in\mathcal{P}_5^{\pm}$ to $\pi_1\Sigma$. In our construction of a polygonal model for one-holed torus representations described in Section~\ref{sec:polygonal-model-one-holed-tori}, we observed that the axes of $\rho(a_1)$ and $\rho([a_1,b_1])$ admit a common perpendicular geodesic line which intersects the axis of $\rho([a_1,b_1])$ at a point, say $Z$ (a posteriori, $Z$ will coincide with the point $D$ of Figure~\ref{fig:pentagon-representation-reflected}). The point $Z$ will serve as natural basepoint on the axis of $\rho([a_1,b_1])$, which is oriented so that $\rho([a_1,b_1])$ translates in the positive direction. Since $[\rho]\in\mathcal{P}_5^{\pm}$, both $\rho(a_2)$ and $\rho(b_2)$ are elliptic rotations of angle $\pi$ with centres on the axis of $\rho([a_1,b_1])$. 
\begin{defn}\label{defn:twist-coordinate-pentagon}
We define the \emph{twist} $\tau_{[a_1,b_1]}$ of $[\rho]$ as the signed distance from $Z$ to the fixed point of $\rho(b_2)$ along the axis of $\rho([a_1,b_1])$. 
\end{defn}

Since $\tau_{[a_1,b_1]}$ is defined as a signed distance and since fixed points and invariant axes of images depend smoothly on $[\rho]$, the twist defines a smooth function
\[
\tau_{[a_1,b_1]}\colon\mathcal{P}_5\to\R.
\]
The above discussion shows that the analogue of Formula~\eqref{eq:Goldman-formula-one-holed-torus} holds. That is, everywhere in $\mathcal{P}_5$, we have
\begin{equation}\label{eq:Goldman-formula-pentagon-twist}
\tau_{[a_1,b_1]}\circ\overline{\Phi}_{[a_1,b_1]}^t=\tau_{[a_1,b_1]}+t. 
\end{equation}

We are now ready to state and prove our parametrization result of $\mathcal{P}_5$. To avoid confusion with the twist of Definition~\ref{defn:twist-coordinate-pentagon}, we will denote by $\tau_{a_1}\colon\mathcal{P}_5\to\R$ the twist coordinate of Definition~\ref{defn:twist-coordinate-one-holed-torus} applied to restrictions of elements in $\mathcal{P}_5$ to representations of $\pi_1\Sigma$.

\begin{thm}\label{thm:coordinates-pentagon-representations}
The map
\[
\mathfrak{T}=\left(\ell_{a_1},\tau_{a_1},\frac{\ell_{[a_1,b_1]}}{2},\tau_{[a_1,b_1]}\right)\colon\mathcal{P}_5\to \R_{>0}\times \R\times \R_{>0}\times \R,
\]
is a smooth surjective map that restricts to a smooth diffeomorphism on both $\mathcal{P}_5^+$ and $\mathcal{P}_5^-$. 
\end{thm}
\begin{proof}
The map $\mathfrak{T}$ is smooth because length functions and twist maps are smooth. To prove the diffeomorphism property, we construct an inverse map. We will use the notation introduced in Remark~\ref{rem:euler-number-pentagon-representation}; namely, $[\rho_{c,e}^+]\in\mathcal{P}_5^+$ and $[\rho_{c,e}^-]\in\mathcal{P}_5^-$ denote the pentagon representations constructed from a right-angled pentagon with two side lengths given by $c/2$ and $e/2$ as in Fact~\ref{fact:right-angled-pentagons}. Consider the two maps $\mathfrak{S}^+$ and $\mathfrak{S}^-$ defined by
\begin{align*}
    \mathfrak{S}^\pm\colon \R_{>0}\times \R\times \R_{>0}\times \R&\to \mathcal{P}_5^\pm\\
    (c,t_1,e,t_2) &\mapsto \left(\Phi_{a_1}^{t_1}\circ \Phi_{[a_1,b_1]}^{t_2}\right) \left(\left[\rho_{c,e}^\pm\right]\right).
\end{align*}
The image of $\mathfrak{S}^{\pm}$ lies in $\mathcal{P}_5^{\pm}$ because of Lemma~\ref{lem:twist-flows-preserves-P5}.

First, we prove that the two maps $\mathfrak{S}^{+}$ and $\mathfrak{S}^{-}$ are smooth functions. Since Hamiltonian diffeomorphisms depend smoothly on time, it suffices to prove that $(c,e)\mapsto [\rho_{c,e}^\pm]$ are smooth functions. The formulas in Fact~\ref{fact:right-angled-pentagons} show that the other three side lengths of a right-angled pentagon depend smoothly on $c$ and $e$. If we further fix the vertex $B$ to be the origin in the Poincaré disk and $A$ to lie below $B$ on the imaginary axis as on Figure~\ref{fig:pentagon-representation-P5}, then the coordinates of the remaining three vertices are also smooth functions of $c$ and $e$. For each such pentagon, the associated representation is defined in~\eqref{eq:representation-right-angled-pentagon}. It maps each fundamental group generator to an element of $\psl$ whose entries depend smoothly on the side lengths and the coordinates of the pentagon vertices. Overall, this shows that $[\rho_{c,e}^\pm]$ are both smooth functions of $c$ and $e$.

Now, we prove that $\mathfrak{T}\circ\mathfrak{S}^{\pm}=\id$. Note that $\ell_{a_1}\circ\mathfrak{S}^{\pm}([\rho_{c,e}^{\pm}])=c$ and $\ell_{[a_1,b_1]}\circ\mathfrak{S}^{\pm}([\rho_{c,e}^{\pm}])=2e$ are consequences of the definition of $[\rho_{c,e}^{\pm}]$ given in~\eqref{eq:representation-right-angled-pentagon}. The relation $\tau_{a_1}\circ\mathfrak{S}^{\pm}=\id$ follows from~\eqref{eq:Goldman-formula-one-holed-torus}, while $\tau_{[a_1,b_1]}\circ\mathfrak{S}^{\pm}=\id$ follows from~\eqref{eq:Goldman-formula-pentagon-twist}. Conversely, to see that $\mathfrak{S}^{\pm}\circ\mathfrak{T}=\id$, first note that $[\rho]$ and $\mathfrak{S}^{\pm}\circ\mathfrak{T}([\rho])$ have the same Fenchel--Nielsen coordinates (Section~\ref{sec:Fenchel-Nielsen-coordinates-one-holed-torus}) when restricted to $\pi_1\Sigma$. Moreover, the two restrictions belong to the same component of the relative character variety $\chi_t(\Sigma,\psl)$ with $t=-2\cosh(\ell_{[a_1,b_1]}([\rho])/2)$ and are therefore equal. Finally, using the geometric interpretation of $\Phi_{[a_1,b_1]}^t$ illustrated on Figure~\ref{fig:pentagon-representation-reflected} and Definition~\ref{defn:twist-coordinate-pentagon}, we conclude that $[\rho]=\mathfrak{S}^{\pm}\circ\mathfrak{T}^{\pm}([\rho])$. Consequently, the restrictions of $\mathfrak{T}$ to $\mathcal{P}_5^+$ and $\mathcal{P}_5^-$ are smooth diffeomorphisms with respective inverses given by $\mathfrak{S}^+$ and $\mathfrak{S}^-$.
\end{proof}

\subsection{Symplectic structure}\label{sec:symplectic-structure-pentagon-representations}
The diffeomorphisms from Theorem~\ref{thm:coordinates-pentagon-representations} turn out to be symplectomorphisms when $\mathcal{P}_5$ is equipped with the restriction of the Goldman symplectic form and $(\R_{>0}\times \R)^2$ carries its standard symplectic structure.
\begin{thm}\label{thm:wolpert-formula-pentagon-representations}
If $\omega_{\mathcal{G}}$ denotes the restriction of the Goldman symplectic form to $\mathcal{P}_5$, then the following Wolpert-like formula holds: 
\begin{equation}\label{eq:Wolpert-formula-pentagon-representation}
\omega_{\mathcal{G}}=-d\ell_{a_1}\wedge d\tau_{a_1}-d\ell_{[a_1,b_1]}\wedge d\tau_{[a_1,b_1]}.
\end{equation}
\end{thm}
\begin{proof}
The map $\mathfrak{T}$ from Theorem~\ref{thm:coordinates-pentagon-representations} is a diffeomorphism on both $\mathcal{P}_5^+$ and $\mathcal{P}_5^-$. Hence, there exist smooth functions $\alpha_{x,y},\beta_{x,y},\gamma_{x,y}$ such that
\[
\omega_{\mathcal{G}}=\sum_{x,y\in\{a_1,[a_1,b_1]\}}\alpha_{x,y}\, d\ell_x\wedge d\ell_y + \beta_{x,y}\, d\tau_x\wedge d\tau_y + \gamma_{x,y}\, d\ell_x\wedge d\tau_y. 
\]
Since the curves $a_1$ and $[a_1,b_1]$ are disjoint, we have $\{\ell_x,\ell_y\}\equiv 0$ for any $x,y$, which implies $\beta_{x,y}\equiv 0$. Here, $\{\cdot,\cdot\}$ denotes the Poissson bracket associated to the Goldman symplectic form, see~\eqref{eq:Poisson-bracket}. Similarly, using Formulas~\eqref{eq:Goldman-formula-one-holed-torus} and~\eqref{eq:Goldman-formula-pentagon-twist}, and by the geometric description of twists flows illustrated on Figure~\ref{fig:pentagon-representation-reflected}, it follows that $\{\ell_x,\tau_y\}\equiv -\delta_{xy}$, where $\delta_{xy}$ stands for the Kronecker delta. Hence, $\gamma_{x,y}\equiv -\delta_{xy}$. In other words, 
\begin{equation}\label{eq:symplectic-form-intermediate-form}
\omega_{\mathcal{G}}=\alpha_{x,y}\, d\ell_x\wedge d\ell_y - d\ell_x\wedge d\tau_x-d\ell_y\wedge d\tau_y,
\end{equation}
where now $x=a_1$ and $y=[a_1,b_1]$. It remains to prove that $\alpha_{x,y}\equiv 0$. As Wolpert explains in~\cite{wolpert-formula}, the coefficient $\alpha_{x,y}$ are independent of the twist parameters. So, it is sufficient to prove that $\alpha_{x,y}\equiv 0$ at points of $\mathcal{P}_5$ with $\tau_x=0$ and $\tau_y=0$. Those points correspond to all $[\rho]\in\mathcal{P}_5$ coming from a right-angled pentagon. What makes representations coming from right-angled pentagons special among $\mathcal{P}_5$ is the fact that they admit special symmetries. For instance, consider the involutive automorphism of $\pi_1S$ given by
\[
\imath\colon\begin{cases}
    a_1\mapsto b_1a_1b_1^{-1}\\
    b_1\mapsto b_1^{-1}\\
    a_2\mapsto a_2^{-1}\\
    b_2\mapsto a_2b_2a_2^{-1}.
\end{cases}
\]
The outer automorphism of $\pi_1S$ represented by $\imath$ is induced by an isotopy class of orientation-reversing homeomorphisms of $S$. So, it acts on the character variety $\chi(S,\psl)$ as an anti-symplectic involution. In other words, $\imath^\ast\omega_{\mathcal{G}}=-\omega_{\mathcal{G}}$. 

The involution $\imath$ should not be mistaken with the involution $\jmath\colon\chi(S,\psl)\to\chi(S,\psl)$ given by conjugation by a matrix of determinant $-1$, see~\eqref{eq:involution-j}. An important difference is that $\jmath$ is a symplectic involution, i.e.~$\jmath^\ast\omega_{\mathcal{G}}=\omega_{\mathcal{G}}$. 
\begin{claim}\label{claim:involution-representation-right-angles-pentagon}
However, if $[\rho]\in\mathcal{P}_5$ comes from a right-angled pentagon, then
\begin{equation}\label{eq:involution-representation-right-angles-pentagon}
    [\rho\circ\imath]=\jmath[\rho].
\end{equation}
\end{claim}
\begin{proofclaim}
To prove~\eqref{eq:involution-representation-right-angles-pentagon}, we will refer to Figure~\ref{fig:pentagon-representation-reflected}. Assume that $[\rho]$ comes the right-angled pentagon $ABCDE$. It follows from the definition of $\imath$ that $\rho\circ\imath(a_1)$ is a hyperbolic translation of the same length as $\rho(a_1)$ but along the axis that is the image by $\rho(b_1)$ of the axis of $\rho(a_1)$, which the geodesic line $(B'C')$ on Figure~\ref{fig:pentagon-representation-reflected}. Similarly, $\rho\circ\imath(b_1)$ is a hyperbolic translation along the same axis as $\rho(b_1)$---the geodesic line $(B'B)$---and of the same length, but in the opposite direction. Furthermore, $\rho\circ\imath(a_2)=\rho(a_2)$ and $\rho\circ\imath(b_2)$ is an elliptic rotation of angle $\pi$ around the image of $D$ by $\rho(a_2)$, which is the point $D'$. In other words, we have just proved that $\rho\circ\imath$ is the representation associated to the right-angled pentagon $AB'C'D'E$ which is obtained from $ABCDE$ by reflection across the geodesic line $(AE)$. This implies $[\rho\circ\imath]=\jmath[\rho]$ and proves the claim.
\end{proofclaim}

A consequence of Claim~\ref{claim:involution-representation-right-angles-pentagon} and the previous discussion is that for every $[\rho]\in\mathcal{P}_5$ coming from a right-angled pentagon, we have
\begin{equation}\label{eq:symplectic-form-involutions}
(\imath^\ast\omega_{\mathcal{G}})_{[\rho]}=-(\omega_{\mathcal{G}})_{[\rho\circ\imath]}=-(\omega_{\mathcal{G}})_{\jmath[\rho]}=-(\jmath^\ast\omega_{\mathcal{G}})_{[\rho]}.
\end{equation}

Now, observe that $\imath^\ast d\ell_x=d\ell_x$ and $\imath^\ast d\ell_y=d\ell_y$. This is because $\rho\circ\imath(a_1)$ is conjugate to $\rho(a_1)$ and $\rho\circ\imath([a_1,b_1])=[b_1,a_1]=[a_1,b_1]^{-1}$. Furthermore, since $\imath$ preserves the axis of $\rho(b_1)$ but flips its orientation, and since the axis of $\rho\circ\imath(a_1)$ is the translate of the axis of $\rho(a_1)$ by $\rho(b_1)$, we have $\imath^\ast d\tau_x=-d\tau_x$. Similarly, $\imath$ preserves the axis of $\rho([a_1,b_1])$ but flips its orientation, while rotating the fixed point of $\rho(b_2)$ by $\rho(a_2)$. This means that the relative position of the basepoint $Z$ used in Definition~\ref{defn:twist-coordinate-pentagon} and the fixed point of $\rho(b_2)$ are inverted by $\imath$, implying that $\imath^\ast d\tau_y=-d\tau_y$. On the other hand, $\jmath$ preserves lengths and twists. Hence, $\jmath^\ast d\ell_i=d\ell_i$ and $\jmath^\ast d\tau_i=d\tau_i$ for $i=x,y$. It is a good exercise to compare the images by $\imath$ and $\jmath$ of the orange hexagon on Figure~\ref{fig:pentagon-representation-reflected} to observe the different actions of the two involutions when the pentagon is not right-angled.

Finally, if we plug the expression~\eqref{eq:symplectic-form-intermediate-form} for $\omega_{\mathcal G}$ into the relation $(\imath^\ast\omega_{\mathcal{G}})_{[\rho]}=-(\jmath^\ast\omega_{\mathcal{G}})_{[\rho]}$ from~\eqref{eq:symplectic-form-involutions} and we use the pullback identities we just computed, we obtain $\alpha_{x,y}([\rho])=-\alpha_{x,y}([\rho])$ and therefore $\alpha_{x,y}([\rho])=0$, as desired. 
\end{proof}

\subsection{Geometrization}\label{sec:geometrization-pentagon-representations} In this section, we explain how to realize each pentagon representation as the holonomy of some branched hyperbolic structure on \(S\). The main result of this section is the following.

\begin{thm}\label{thm:geometrization_pentagon_reps}
Every pentagon representation arises as the holonomy of branched hyperbolic structures on \(S\).
\end{thm}

Before proving Theorem~\ref{thm:geometrization_pentagon_reps}, we recall a few preliminaries needed for the proof. Let \(\rho\) be a pentagon representation and assume \([\rho]\in\mathcal P^+\); if \([\rho]\in\mathcal P^-\) instead, then the same arguments apply to $S$ with opposite orientation.

We pick a standard geometric presentation of $\pi_1S$ as in~\eqref{eq:presentation-pi1S}:
\begin{equation}\label{eq:presentation-pi_1S-pentagon}
\pi_1S=\langle a_1,b_1,a_2,b_2: [a_1,b_1][a_2,b_2]=1\rangle.
\end{equation}
Up to applying a mapping class group element to $[\rho]$, we may assume that $[\rho]\in \mathcal{P}_5^+$ for the presentation of $\pi_1S$ above (see Table~\ref{tab:6-components-pentagons} and the discussion thereafter). The simple closed curve \(c\) on $S$ represented by the fundamental group element $[a_1,b_1]$ splits \(S\) into two sub-surfaces \(\Sigma_1\) and \(\Sigma_2\). By choosing basepoints on $c$, it is possible to see the fundamental groups of $\Sigma_1$ and $\Sigma_2$ as the free subgroups of $\pi_1S$ generated by $a_1,b_1$ and $a_2,b_2$. We denote by $\rho_i$ the restriction of $\rho$ to $\pi_1\Sigma_i$.

Since we are assuming that $[\rho]\in\mathcal{P}_5^+$, we have $\eu(\rho)=1$ by definition of $\mathcal{P}^+$. Furthermore, the results of~\cite{BIW} on Toledo number implies that \(\eu(\rho)=\Tol(\rho)=\Tol(\rho_1)+\Tol(\rho_2)\). 
Since $\rho$ is a pentagon representation, $\rho([a_1,b_1])$ is hyperbolic by Lemma~\ref{lem:properties-pentagon}. Therefore, using~\eqref{eq:Toldeo-number-one-holed-torus} and up permuting the indices, we may suppose that $\Tol(\rho_1)=1$ and $\Tol(\rho_2)=0$.

\begin{fact}\label{fact_pentadomain}
    The representation \(\rho_1\) is the holonomy of a complete hyperbolic structure on \(\Sigma_1\) with totally geodesic boundary \(\gamma\) (as explained in Section \ref{sec:Fenchel-Nielsen-coordinates-one-holed-torus}). Furthermore, there exists a pentagonal fundamental domain \(\mathcal D\) for action of \(\textnormal{Im}(\rho_1)\) on the hyperbolic plane. (Remark: this is not the polygonal model of $\rho_1$ that we described in Section~\ref{sec:polygonal-model-one-holed-tori}.)
\end{fact}

Although this claim is standard, we provide a proof here, not only to make our discussion as self-contained as possible, but also because the argument will be needed to show that \(\rho\) is geometrizable.

\begin{figure}[h]
    \centering
    \begin{tikzpicture}
    \node[anchor = south west, inner sep=0mm] at (0,0) {\includegraphics[width=11cm]{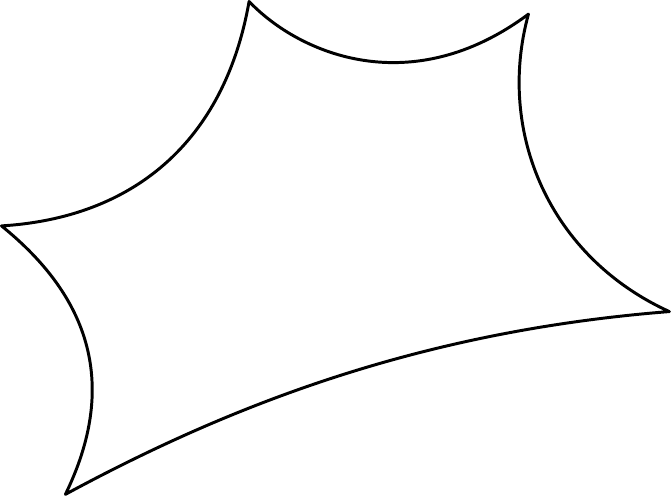}};

    \node at (4.1,8.4) {$P$};
    \node at (9,8.25) {$gP$};
    \node at (-.4,4.4) {$hP$};
    \node at (11,2.6) {$ghP=[g,h]Q$};
    \node at (.9,-.2) {$hgP=Q$};
    \draw (3.5,5.5) edge[out=340,in=200,->] (8.5,5);
    \draw (6.4,6.7) edge[out=250,in=20,->] (1.9,2.7);
    \node at (6.4, 4.3) {$g$};
    \node at (4.1,3.3) {$h$};
    \end{tikzpicture}
    \caption{The fundamental domain $\mathcal D$.}
    \label{fig:fundamental-domain-pentagon}
\end{figure}

\begin{proof}[Proof of Fact~\ref{fact_pentadomain}]
    Consider the hyperbolic structure on $\Sigma_1$ with holonomy $\rho_1$. We pick a basepoint $p$ on its boundary $\gamma$ and consider the loop $\alpha$ based at $p$ that is freely homotopic to $a_1$. When we cut $\Sigma_1$ along $\alpha$, we obtain a hyperbolic cylinder with two boundary curves that have corner points. One is a copy of $\alpha$ with corner point $p_1$ and the other one has two corner points $p_2$ and $p_3$ that partition it into a copy of $\gamma$ and a copy of $\alpha$.

    We now consider the shortest curve $\delta$ from $p_1$ to $p_2$. It is piecewise geodesic, and cannot visit $p_1$ or $p_2$ more than once by minimality. It could visit $p_3$ however. If it is the case, then we let $\beta$ be the segment from $p_1$ to $p_3$. If $\delta$ does not pass through $p_3$, then we let $\beta=\delta$. 
    
    When we cut our cylinder along $\beta$, we obtain a contractible domain that we can embed in the hyperbolic plane as a convex pentagon. To be more concrete, we embed it as follows. Let $\tilde p\in\widetilde \Sigma_1$ be a lift of $p$ and let $\dev\colon \widetilde \Sigma_1\to \HH$ be a $\rho_1$-equivariant developing map. We let $P=\dev(\tilde p)$, as well as $g=\rho_1(a_1)$ and $h=\rho_1(b_1)$. Our desired pentagon $\mathcal{D}$, illustrated on Figure~\ref{fig:fundamental-domain-pentagon} is bounded by the chain of segments
    \[ P \longmapsto gP \longmapsto ghP \longmapsto hgP \longmapsto hP \longmapsto P.
    \]
    The side $[ghP, hgP]$ is the developed image of the boundary curve $\gamma$ of $\Sigma_1$. Note that if we let $Q=hgP$, then $ghP=[g,h]Q$.
\end{proof}

\begin{rmk}\label{rmk_inner_angles}
    The sum of the interior angles in the pentagonal domain \(\mathcal{D}\) constructed in Fact~\ref{fact_pentadomain} is $\pi$. This is because the five vertices are identified with the basepoint $p$ on the boundary of $\Sigma_1$, and two pairs of sides are identified, corresponding respectively to the curves \(\alpha\) and \(\beta\) based at $p$, while the last side corresponds to the boundary of \(\Sigma_1\) given by the commutator \([\alpha, \beta]\).
\end{rmk}

\begin{rmk}\label{rmk_fundamental_domain}
    The hyperbolic plane is tessellated by isometric copies of the fundamental domain \(\mathcal D\). Notice that each tile of this tessellation can be equally entitled as fundamental domain for $\Sigma_1$ and all other tiles arises as translated copies of \(\mathcal D\) by elements of \(\textnormal{Im}(\rho_1)\). In particular, there is a (not unique) tile such the edge joining the vertices \(Q\) with \([g,h](Q)\) lies on the axis of the commutator \([g,h]\). In what follows, we shall regard these tiles as \textit{preferred domains}.
\end{rmk}

So far, we have geometrically realized only “half” of the representation \(\rho\), and we still need to treat the restriction \(\rho_2\). Naively, a fundamental domain for $\rho$ is expected to be the gluing of the fundamental domains of $\rho_1$ and $\rho_2$. However, since we are assuming that $[\rho]\in\mathcal{P}_5^+$, the image $\Im(\rho_2)$ is generated by two elliptic elements of order two with distinct fixed points. In particular, $\rho_2$ is virtually abelian, and therefore it cannot arise as the holonomy of a hyperbolic cone-manifold structure on $\Sigma_2$ with cornered geodesic boundary~\cite[Theorem~1]{mathews-torus}. This means that we cannot find a fundamental domain for $\rho_2$ as we did for $\rho_1$. Our alternative construction requires the following observation.

\begin{fact}\label{fact_virtuallyabelian_char}
    Let \(e_1,e_2\in\psl\) be two elliptic elements of order two with distinct fixed points. Let \(\ell\) be the unique geodesic line passing through their fixed points. Both transformations \(e_1e_2\) and \([e_1,e_2]\) are hyperbolic translations along \(\ell\) and every element in \(\langle e_1,e_2\rangle\) is either hyperbolic or elliptic with order two. Furthermore, the fixed points of all the elliptic elements in \(\langle e_1,e_2\rangle\) equidistribute along \(\ell\).
\end{fact}

\begin{proof}
        There is a computational argument using matrices, but the geometry here is more revealing and direct. Let \(p_1\in \ell\) and \(p_2\in\ell\) be the fixed points of \(e_1\) and \(e_2\). Denote by \(a_\ell\) the endpoint of the geodesic ray \([p_2p_1)\) and by \(r_{\ell}\) the endpoint of \([p_1p_2)\). Note that \(e_1e_2\) preserves \(\ell\) and hence is hyperbolic. Moreover, \(a_{\ell}\) and \(r_\ell\) are the attractive and repelling fixed points of \(e_1e_2\). The translation length of \(e_1e_2\) is equal to \(2d(p_1,p_2)\). Since, \(e_1^2=e_2^2=1\), it follows that \([e_1,e_2]=(e_1e_2)^2\) and hence \([e_1,e_2]\) is also hyperbolic and its translation length is equal to \(4d(p_1,p_2)\).
        
        For the second part, notice that \(\langle e_1,e_2\rangle\) is a virtually abelian group and its cyclic subgroup \(\langle e_1e_2\rangle\) has index two. Hence 
        \[\quotient{\langle e_1,e_2 \rangle}{\langle e_1e_2\rangle}\cong\langle e_1\rangle\cong\Z/2\Z.
        \] 
        Consequently, each reduced word in \(e_1\) and \(e_2\) is either hyperbolic or elliptic depending on whether its length is even or odd, respectively. In particular, all elements of the form $(e_1e_2)^ne_1$ are elliptic and their fixed point is given by
        \[
        \begin{cases}
            (e_1e_2)^kp_1 &\text{if } n=2k,\\
            (e_1e_2)^{k+1}p_2 &\text{if } n=2k+1.
        \end{cases}
        \] 
        We conclude that a point on \(\ell\) is the fixed point of an elliptic element in $\langle e_1,e_2 \rangle$ if and only if it lies at distance \(kd(p_1,p_2)\) from \(p_1\) for some integer $k$.
\end{proof}
  
In what follows we set \(e_1=\rho_2(a_2)\) and \(e_2=\rho_2(b_2)\) and we set \(\Gamma=\Im(\rho_2)=\langle e_1,e_2\rangle\) for simplicity. We are now ready to geometrize pentagon representations and prove Theorem~\ref{thm:geometrization_pentagon_reps}.

\begin{figure}[h]
    \centering
    \begin{tikzpicture}
    \node[anchor = south west, inner sep=0mm] at (0,-.2) {\includegraphics[width=12cm]{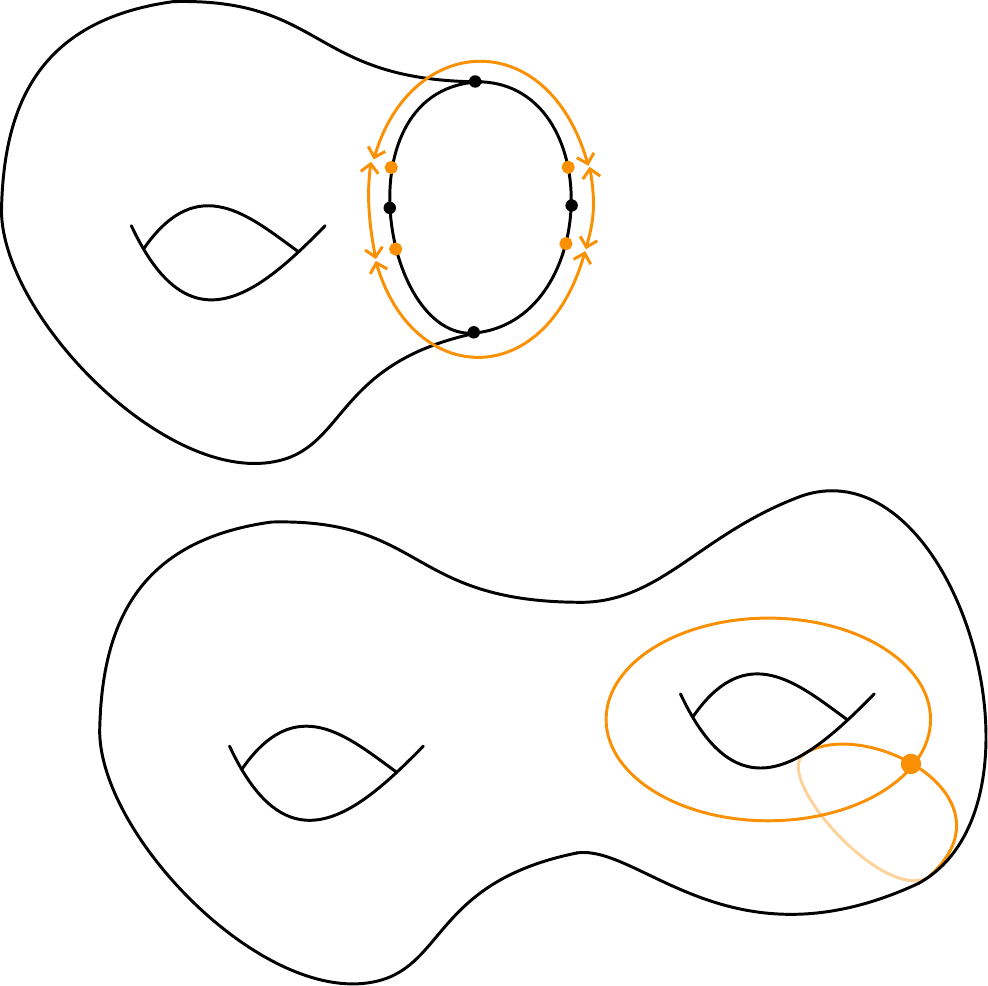}};

    \draw[mauve] (6,7.54) edge[out=80, in=280, <->] (6,10.97);
    \draw[mauve] (4.55,9.1) edge[out=-5, in=185, <->] (7.15,9.1);

    \draw (6.8,7.6) edge[out=-30, in=120, ->] (8.8, 5.5);

    \node[apricot] at (11.68, 2.35) {$\lambda_1$};
    \node[apricot] at (7.05, 3) {$\lambda_2$};
    \node[apricot] at (10.9, 2.85) {$4\pi$};
    \end{tikzpicture}
    \caption{Step-by-step gluing of the fundamental domain $\mathcal{D}'$ first into the hyperbolic torus with boundary above, and then into a branched hyperbolic surface of genus $2$ below.}
    \label{fig:pentagon-geometrization}
\end{figure}

\begin{proof}[Proof of Theorem \ref{thm:geometrization_pentagon_reps}]
    Recall that our goal is exhibit the pentagon representation $\rho$ as the holonomy of a branched hyperbolic structure on $S$. We are assuming that \([\rho]\in\mathcal P^+_5\), and we denoted by \(\rho_1\) and \(\rho_2\) the restrictions of $\rho$ to $\pi_1\Sigma_1$ and $\pi_1\Sigma_2$. We are working with the geometric presentation~\eqref{eq:presentation-pi_1S-pentagon} of $\pi_1S$ with generators $a_1,b_1,a_2,b_2$ and $\pi_1\Sigma_i$ is the free subgroup of $\pi_1S$ generated by $a_i,b_i$. 
    
    According to Fact~\ref{fact_pentadomain}, \(\rho_1\) is the holonomy of some complete hyperbolic structure on \(\Sigma_1\) with totally geodesic boundary having holonomy \([g,h]\), where $g=\rho_1(a_1)$ and $h=\rho_1(b_1)$. The invariant axis of \([g,h]\) is the geodesic line $\ell$. Moreover, according to Fact~\ref{fact_virtuallyabelian_char}, \(\Gamma=\textnormal{Im}(\rho_2)\) is virtually abelian and generated by \(e_1\) and \(e_2\) such that \([g,h]=[e_1,e_2]^{-1}\) by design. As before, let $p_1$ and $p_2$ be the fixed points of $e_1$ and $e_2$. Both $p_1$ and $p_2$ belong to $\ell$ by Fact~\ref{fact_virtuallyabelian_char}. Pick a preferred fundamental domain $\mathcal{D}$ (in the sense of Remark~\ref{rmk_fundamental_domain}) for $\rho_1$. The distance between the vertices $Q$ and $[g,h]Q$ of $\mathcal{D}$ is equal to the translation length of $[g,h]$, which is equal to $4d(p_1,p_2)$ as we have seen in the proof of Fact~\ref{fact_virtuallyabelian_char}. This means that the hyperbolic segment between $Q$ and $[g,h]Q$ contains at least four fixed points (sometimes five) of the elliptic elements in $\Gamma$ because those are equidistributed on $\ell$ with distance $d(p_1,p_2)$ between them (Fact~\ref{fact_virtuallyabelian_char}). So, up to changing the generators $a_2, b_2$ of $\pi_1S$ by applying an automorphism of $\pi_1\Sigma_2$ that preserves $[a_2,b_2]$, we may assume that $p_1$ and $p_2$ lie strictly between $Q$ and $[g,h]Q$.

    If $\lambda=d(Q,[g,h]Q)$, then we choose two positive real numbers $\lambda_1$ and $\lambda_2$ such that
    \begin{equation}\label{eq:domain-lambda}
        \lambda =2\lambda_1 +2\lambda_2.
    \end{equation}
    For each image of $p_1$ by an element of $\Gamma$ (those points always lie on $\ell$), we consider a small segment of length $\lambda_1$ centered at that point. Similarly, we consider segments of length $\lambda_2$ around each image of $p_2$ by $\Gamma$. All those segments partition the geodesic line $\ell$ in a way that a segment of length $\lambda_i$ is followed by a segment of length $\lambda_j$, and then another segment of length $\lambda_i$, and so on (see Figure~\ref{fig:fundamental-domain-pentagon-marked}). 

    \begin{figure}[h]
    \centering
    \begin{tikzpicture}
    \node[anchor = south west, inner sep=0mm] at (0,-.2) {\includegraphics[width=11cm]{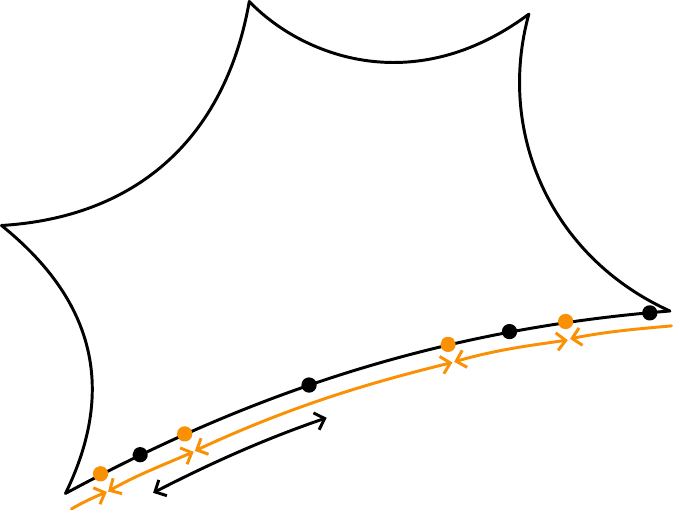}};

    \node at (4.1,8.4) {$P$};
    \node at (9,8.25) {$gP$};
    \node at (-.4,4.4) {$hP$};
    \node at (12.2,3.1) {\small $ghP=[g,h]Q$};
    \node at (.3,.1) {\small $hgP=Q$};
    \draw (3.5,5.5) edge[out=340,in=200,->] (8.5,5);
    \draw (6.4,6.7) edge[out=250,in=20,->] (1.9,2.7);
    \node at (6.4, 4.3) {$g$};
    \node at (4.1,3.3) {$h$};

    \node at (2.35, 1.1) {$p_1$};
    \node at (5, 2.2) {$p_2$};
    \node at (4, 0.4) {$\lambda/4$};
    \node[apricot] at (6.4, 1.6) {$\lambda_2$};
    \node[apricot] at (8.5, 2.15) {$\lambda_1$};

    \draw[mauve] (2.7, .4) edge[out=-40, in=240, ->] (7.9, 2.2);
    \node[mauve] at (6.6, .3) {$e_2$};
    \draw[mauve] (6.8, 1.9) edge[out=-40, in=260, ->] (9.9, 2.6);
    \node[mauve] at (8.6, 1.2) {$e_2e_1e_2$};
    \end{tikzpicture}
    \caption{The $9$-gone $\mathcal{D}'$ and the subdivision of the segment between $Q$ and $[g,h]Q$ by the segments of length $\lambda_i$ centered at $p_i$ (in orange).}
    \label{fig:fundamental-domain-pentagon-marked}
\end{figure}

    Since we are assuming that $p_1$ and $p_2$ lie strictly between $Q$ and $[g,h]Q$, those segments have four endpoints between $Q$ and $[g,h]Q$, colored in orange on Figure~\ref{fig:fundamental-domain-pentagon-marked}. By counting these four points as new vertices, our pentagon $\mathcal{D}$ is turned into a $9$-gone $\mathcal{D}'$ whose interior angles now sum to $\pi+4\pi=5\pi$. This is the desired fundamental domain for $\rho$. In fact, when we quotient $\mathcal{D}'$ by the obvious side identifications given by the isometries $g$, $h$, as well as the appropriate order-$2$ elements in $\Gamma$, the resulting space is a closed surface of genus $2$. The induced hyperbolic structure is a branched hyperbolic structure with one branch point of angle \(4\pi\) given by the identification of the endpoints of the segments of lengths $\lambda_1$ and $\lambda_2$ in the quotient, as we illustrate on Figure~\ref{fig:pentagon-geometrization}. By construction, this structure has holonomy $\rho$.    
\end{proof}

\begin{cor}\label{cor:isomonodromic-deformations-pentagon}
    All the branched hyperbolic structures constructed in the proof of Theorem~\ref{thm:geometrization_pentagon_reps} have one degree of freedom given by the choice of $\lambda_1$ and $\lambda_2$ satisfying~\eqref{eq:domain-lambda}. All these structures are distinct, and therefore define a $1$-dimensional real locus is the isomonodromic leaf over $[\rho]$.
\end{cor}
\begin{proof}
    The boundary curve of $\Sigma_1$ turns into a figure-$8$ curve on $S$ made of two geodesic loops of respective lengths $\lambda_1$ and $\lambda_2$, and meeting at the branch point (see Figure~\ref{fig:pentagon-geometrization}). Both $\lambda_1$ and $\lambda_2$ are therefore intrinsic parameters of the branched hyperbolic structure, and different choices of $\lambda_1$ and $\lambda_2$ lead to different branched hyperbolic structures.
\end{proof}
 
\section{Bow-tie representations}\label{sec:bow-tie-representations}
\subsection{Overview}
In this section, we introduce \emph{bow-tie} representations (Definition~\ref{defn:bow-tie-representations}) and state their main properties (Section~\ref{sec:definiton-bow-tie-representations}). We also describe an analogue of Fenchel--Nielsen coordinates for bow-tie representations (Section~\ref{sec:parametrization-bow-tie-representations}) and show that these coordinates satisfy a Wolpert-like formula (Section~\ref{sec:symplectic-structure-bow-tie-representations}). The section ends with a description of branched hyperbolic metrics whose holonomies realize all bow-tie representations (Section~\ref{sec:geometrization-bow-tie-representations}) and their isomondromic deformations.

\subsection{Definition and properties}\label{sec:definiton-bow-tie-representations}
As before and throughout this section, $S$ will denote a closed and oriented surface of genus $2$. We are interested in all the representations $\rho\colon\pi_1S\to\psl$ which map at least one separating simple closed curve on $S$ to an elliptic element of $\psl$. Recall that \emph{elliptic} elements of $\psl$ are those fixing a unique point in the hyperbolic plane (in particular, the identity is not elliptic). Such a representation may only have Euler number $0$ or $\pm 1$. 
\begin{defn}\label{defn:bow-tie-representations}
A representation $\rho\colon\pi_1S\to \psl$ is called a \emph{bow-tie representation} if its Euler number is $\pm 1$ and if there exists a separating simple closed curve $c$ on $S$ represented by a fundamental group element $\gamma\in\pi_1S$ such that $\rho(\gamma)$ is elliptic. The set of conjugacy classes of bow-tie representations will be denoted by $\mathcal{B}\subset \chi(S,\psl)$.
\end{defn}
The motivation behind the name ``bow-tie'' will be justified later in the discussion that leads to Figure~\ref{fig:bow-tie}. Note that $\mathcal{B}$ decomposes as 
\[
\mathcal{B}=\bigcup_{\substack{c\colon\text{separating}\\ \text{simple closed}\\ \text{curve on }S}} \mathcal{B}_c,
\]
where $\mathcal{B}_c\subset \mathcal{B}$ is the subset of all $[\rho]$ for which $\rho(\gamma)$ is elliptic for some (hence any) fundamental group representative $\gamma$ of $c$. We will see below (Theorem~\ref{thm:coordinates-bow-tie-representations}) that every point in $\mathcal{B}$ comes from the gluing of two representations of a one-holed torus with complementary elliptic boundary holonomy. For now, we will assume that $\mathcal{B}$ is nonempty. We start with some elementary properties of $\mathcal{B}$.

\begin{lem}\label{lem:properties-of-E}
Bow-tie representations satisfy the following:
\begin{enumerate}
    \item The set $\mathcal{B}$ is invariant under the action of $\Mod(S)$. \label{cond:E-mod-invriant}
    \item It is an open subset of $\chi(S,\psl)$ which is dense inside the two components of Euler number $\pm 1$.\label{cond:E-open-and-dense}
    \item The set $\mathcal{P}$ of pentagon representations (Definition~\ref{defn:pentagon-representation}) is disjoint from $\mathcal{B}$.\label{cond:E-and-P-disjoint}
\end{enumerate}
\end{lem}
\begin{proof}
Condition~\eqref{cond:E-mod-invriant} is a consequence of the fact that $\Mod(S)$ preserves separating simple closed curves. To see that $\mathcal{B}$ is open, note that each $\mathcal{B}_c$ is open because the set of elliptic elements is an open subset of $\psl$. By the ergodicity result of Marché--Wolff (Theorem~\ref{thm:goldman-ergocicity-genus-2}), since $\mathcal{B}$ is open, nonempty, and $\Mod(S)$-invariant, its complement must have zero Goldman measure in the two components of $\chi(S,\psl)$ with Euler number $\pm 1$. This implies that $\mathcal{B}$ is dense in each of those components and proves~\eqref{cond:E-open-and-dense}. Condition~\eqref{cond:E-and-P-disjoint} follows immediately from Lemma~\ref{lem:properties-pentagon}.
\end{proof}

As we did for pentagon representations, we will write
\[
\mathcal{B}=\mathcal{B}^+\sqcup \mathcal{B}^-,
\]
where $\mathcal{B}^+$ and $\mathcal{B}^-$ denotes the elements of $\mathcal{B}$ with Euler number $+1$, respectively $-1$. Similarly, we let $\mathcal{B}_c^{+}=\mathcal{B}^+\cap \mathcal{B}_c$ and $\mathcal{B}_c^{-}=\mathcal{B}^-\cap \mathcal{B}_c$. Recall that $\mathcal{B}^+$ and $\mathcal{B}^-$ are images of each other by the symplectic involution of $
\chi(S,\psl)$ given by conjugation by a determinant $-1$ matrix. Moreover, any mapping class group element that maps a separating simple closed curve $c$ to another curve $c'$ also defines a symplectomorphism between $\mathcal{B}^{\pm}_c$ and $\mathcal{B}_{c'}^{\pm}$.

\subsection{Coordinates}\label{sec:parametrization-bow-tie-representations}
In order to understand the symplectic geometry of $\mathcal{B}$, we will now describe an analogue of Fenchel--Nielsen coordinates for each $\mathcal{B}_c$ individually. Like for traditional Fenchel--Nielsen coordinates, our coordinates on $\mathcal{B}_c$ depend on two pieces of data. First, a choice of pants decomposition of $S\setminus c$, which corresponds to the choice of two simple closed curves, one on each of the two connect components of $S\setminus c$. Second, a set of seams to measure twist parameters. We will materialize these choices by a geometric presentation of $\pi_1S$ of the same form as in~\eqref{eq:presentation-pi1S}:
\begin{equation}
\pi_1S=\langle a_1,b_1,a_2,b_2\, : \, [a_1,b_1][a_2,b_2]=1\rangle.
\end{equation}
If $c$ is represented by the fundamental group element $[a_1,b_1]$, we say that the geometric presentation of $\pi_1S$ is \emph{compatible} with $c$. The simple closed curves represented by $a_1$ and $a_2$ provide the pants decomposition of $S\setminus c$, and those represented by $b_1$ and $b_2$ are the seams data. 

For any $[\rho]\in\mathcal{B}_c$, we can restrict $\rho$ to the two free subgroups of $\pi_1S$ generated by $a_1, b_1$, respectively $a_2,b_2$. These two groups are the fundamental groups of the two one-holed tori $\Sigma_1$ and $\Sigma_2$ obtained by cutting $S$ along $c$. We will denote the two restrictions by 
\[
\rho_1\colon \pi_1\Sigma_1\to\psl, \quad \rho_2\colon\pi_1\Sigma_2\to\psl.
\]
By construction, $\rho_1([a_1,b_1])=\rho_2([a_2,b_2])^{-1}$ are both elliptic fixing the same point in the upper half-plane. So, according to~\eqref{eq:Toldeo-number-one-holed-torus}, there exist $\beta_1,\beta_2\in (0,2\pi)$ such that $\Tol(\rho_1)=\pm \beta_1/2\pi$ and $\Tol(\rho_2)=\pm \beta_2/2\pi$. By the additivity of the Toledo number, it follows that $\pm \beta_1\pm\beta_2 = 2\pi$ if $[\rho]\in\mathcal{B}_c^+$ and $\pm \beta_1\pm\beta_2 = -2\pi$ if $[\rho]\in\mathcal{B}_c^-$. Since we are assuming $\beta_1,\beta_2\in (0,2\pi)$, the first equality holds if and only if $\Tol(\rho_1)=\beta_1/2\pi$ and $\Tol(\rho_2)=\beta_2/2\pi$ with $\beta_2=2\pi-\beta_1$. Similarly, the second equality holds if and only if $\Tol(\rho_1)=-\beta_1/2\pi$ and $\Tol(\rho_2)=-\beta_2/2\pi$ with $\beta_2=2\pi-\beta_1$. So, according to~\eqref{eq:Toldeo-number-one-holed-torus} again, if we write $t=2\cos(\beta_1/2)$, then in both cases we have
\begin{equation}\label{eq:relative-character-variety-restrictions-rho_1-rho_2}
[\rho_1]\in\chi_t(\Sigma_1,\psl),\quad [\rho_2]\in\chi_{-t}(\Sigma_2,\psl).
\end{equation}

\begin{rmk}\label{rmk:gluing-opertaion}
The restriction operation we just described can be inverted by a gluing operation. Two representations $\rho_1\colon\pi_1\Sigma_1\to\psl$ and $\rho_2\colon\pi_2\Sigma_2\to\psl$ can be \emph{glued} back into a representation $\rho\colon\pi_1S\to\psl$ if $\rho_1([a_1,b_1])=\rho_2([a_2,b_2])^{-1}$ because
\[
\pi_1S=\quotient{\pi_1\Sigma_1 \ast \pi_1\Sigma_2}{[a_1,b_1][a_2,b_2]}.
\]
Recall that each of the two relative character varieties $\chi_t(\Sigma_1,\psl)$ and $\chi_{-t}(\Sigma_2,\psl)$ has two connected components distinguished by the Toledo number. What the previous discussion shows (compare with Table~\ref{tab:different-cases-ellipitc-regime}) is that those four components can be paired in a unique way in order for the gluing operation along boundary curves to be possible and for the resulting representation to have Euler number $\pm 1$. Note that it is also possible to glue two representations, one from each of the two connected components of the same relative character variety, but then the resulting representation will have Euler number $0$. 
\end{rmk}

\begin{defn}\label{defn:beta-angle}
The \emph{$\beta$-angle} is the function
\begin{equation}\label{eq:angle-beta}
    \beta\colon\mathcal{B}_c\to (0,2\pi),
\end{equation}
where $\beta([\rho])$ is defined as the \emph{clockwise} rotation angle of $\rho([a_1,b_1])$.
\end{defn}

\begin{rmk}\label{rem:sign-beta-angle}
There is a subtlety in the definition of the $\beta$-angle that is worth emphasizing. Given $[\rho]\in\mathcal{B}_c$ and $\rho_1$ the restriction of $\rho$ to $\pi_1\Sigma_1$, recall from~\eqref{eq:relative-character-variety-restrictions-rho_1-rho_2} that $[\rho_1]\in\chi_t(\Sigma_1,\psl)$, where $t=2\cos(\beta_1/2)$ and $\beta_1=2\pi\vert \Tol(\rho_1)\vert$. Since $\rho([a_1,b_1])=\rho_1([a_1,b_1])$, we may be tempted to say that $\beta([\rho])$ is equal to $\beta_1$, but that is not always true. The correct relation  between the two quantities is given by Lemma~\ref{lem:direction-of-rotation-one-holed-torus} and Table~\ref{tab:different-cases-ellipitc-regime}. It reads
\[
\beta([\rho])=\begin{cases}
    \beta_1 & \text{if $\Tol(\rho_1)< 0$}\\
    2\pi-\beta_1 & \text{if $\Tol(\rho_1)> 0$.}
\end{cases}
\]
In particular, the trace of $\rho([a_1,b_1])$ is not always equal to $2\cos(\beta([\rho])/2)$.
\end{rmk}

\begin{figure}[h]
\centering
\begin{tikzpicture}
\node[anchor = south west, inner sep=0mm] at (0,0) {\includegraphics[width=13cm]{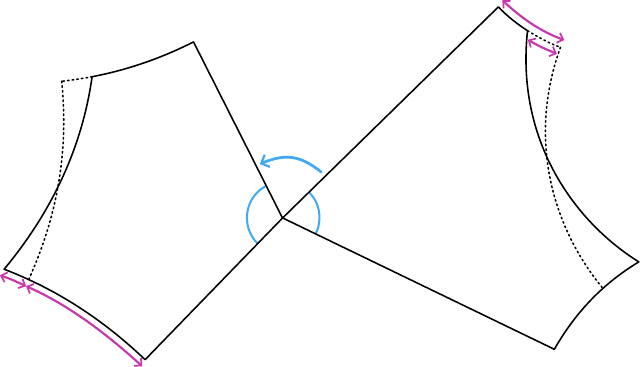}};

\node[sky] at (6,4.5) {\Large $\gamma$};
\node[sky] at (6.95,3.15) {\large $\beta/2$};
\node[sky] at (4.2,3.15) {\large $\pi-\beta/2$};
\node[mauve] at (1.45,.65) {\large $\ell_1/2$};
\node[mauve] at (0.15,1.37){$\tau_1/2$};
\node[mauve] at (11,7.3) {\large $\ell_2/2$};
\node[mauve] at (11.85,6.4) {$\tau_2/2$};

\node at (5.75, 2.75) {$B$};
\node at (-.2, 2.1) {$X_1$};
\node at (3.15, -0.05) {$Y_1$};
\node at (3.95, 6.8) {$Y_1'$};
\node at (1.85, 6.2) {$X_1'$};

\node at (11.3,.1) {$Y_2'$};
\node at (13.25,2.1) {$X_2'$};
\node at (9.9,7.5) {$Y_2$};
\node at (10.4,6.75) {$X_2$};
\end{tikzpicture}
\caption{The bow-tie of $[\rho]\in\mathcal{B}_c^-$ is made of the two pentagons $X_1Y_1BY_1'X_1'$ and $X_2Y_2BY_2'X_2'$ associated to the restrictions $[\rho_1]$ and $[\rho_2]$. The angles at $B$ are $\pi-\beta/2$ and $\beta/2$. The Fenchel--Nielsen coordinates of $[\rho_1]$ are $(\ell_1,\tau_1)$, and those of $[\rho_2]$ are $(\ell_2,\tau_2)$ (compare with Figure~\ref{fig:fenchel-nielsen-coordinates-one-holed-torus}).}
\label{fig:bow-tie}
\end{figure}

In order to parametrize $\mathcal{B}_c$, we first describe a polygonal model for bow-tie representations. One should keep in mind that polygons are associated to representations and orientation-preserving isometry classes of polygons are associated to conjugacy classes of representations, but we will often omit this difference. The polygon associated to $[\rho]\in \mathcal{B}_c$ is made of the two pentagons associated to the restrictions $[\rho_1]$ and $[\rho_2]$ that we described in Section~\ref{sec:polygonal-model-one-holed-tori}. The two pentagons are linked at the common fixed point of $\rho_1([a_1,b_1])$ and $\rho_2([a_2,b_2])$ which we label $B$. Because of~\eqref{eq:relative-character-variety-restrictions-rho_1-rho_2} and Table~\ref{tab:different-cases-ellipitc-regime}, the interior angle at $B$ in the pentagon of $[\rho_1]$ is $\pi-\beta_1/2$, while the angle at $B$ in the pentagon of $[\rho_2]$ is $\beta_1/2$. Here, $\beta_1$ refers to the angle introduced before~\eqref{eq:relative-character-variety-restrictions-rho_1-rho_2}. So, in light of Remark~\ref{rem:sign-beta-angle}, the interior angles at $B$ are $\beta/2$ and $\pi-\beta/2$ if $[\rho]\in\mathcal{B}_c^+$, respectively $\pi-\beta/2$ and $\beta/2$ if $[\rho]\in\mathcal{B}_c^-$, where $\beta=\beta([\rho])$ is the $\beta$-angle of $[\rho]$ from Definition~\ref{defn:beta-angle}. As Figure~\ref{fig:bow-tie} illustrates, the resulting polygon is shaped as a (possibly twisted) bow-tie and is our polygonal model for $[\rho]$. 

Several comments are in order.
\begin{itemize}
    \item Each of the two pentagons that make up the polygonal model of $[\rho]$ may be non-convex (see Figure~\ref{fig:polygon-non-convex} and the discussion beforehand).
    \item The two pentagons may overlap.
    \item The cyclic orientation of the vertices of each pentagon, say in the case they are both convex, is the same. In other words, the parameter $\theta$ introduced after~\eqref{eq:Toldeo-number-one-holed-torus} belongs to the same interval for both $\rho_1$ and $\rho_2$. This is because the Toledo numbers of the restrictions $\rho_1$ and $\rho_2$ have the same sign, as we observed before~\eqref{eq:relative-character-variety-restrictions-rho_1-rho_2}.
    \item Every such bow-tie determines a unique representation $\rho$ and every orientation-preserving isometry class of bow-ties determines a unique conjugacy class $[\rho]\in\mathcal{B}_c$.
\end{itemize}

The polygonal model for bow-tie representations is particularly useful as it allows us to extract a coordinate that will be the symplectic dual to the angle $\beta$.
\begin{defn}\label{defn:gamma-angle-bow-tie-representation}
The angle $\gamma\in\R/2\pi\Z$ of $[\rho]\in\mathcal{B}_c$ is the oriented angle at the vertex $B$ measured anti-clockwise from the geodesic ray $[BY_2)$ to the geodesic ray $[BY_1')$ (see Figure~\ref{fig:bow-tie}).
\end{defn}

\begin{ex}\label{ex:discrete-bow-tie}
    Bow-tie representations may have discrete image. To produce such examples, one approach consists in finding a bow-tie similar to the one of Figure~\ref{fig:bow-tie} that ``fits well'' on a triangular tessellation of the hyperbolic plane. By fitting well, we mean that the translations along the edges of the bow-tie which are the images of the four fundamental group generators by the associated representation are symmetries of the tessellation. This is the case, for instance, for the bow-tie illustrated on Figure~\ref{fig:discrete-bow-tie}. This ensures that the image of the representation is a discrete subgroup of $\psl$.
\begin{figure}[h]
\centering
\begin{tikzpicture}
\node[anchor = south west, inner sep=0mm] at (0,0) {\includegraphics[width=11cm]{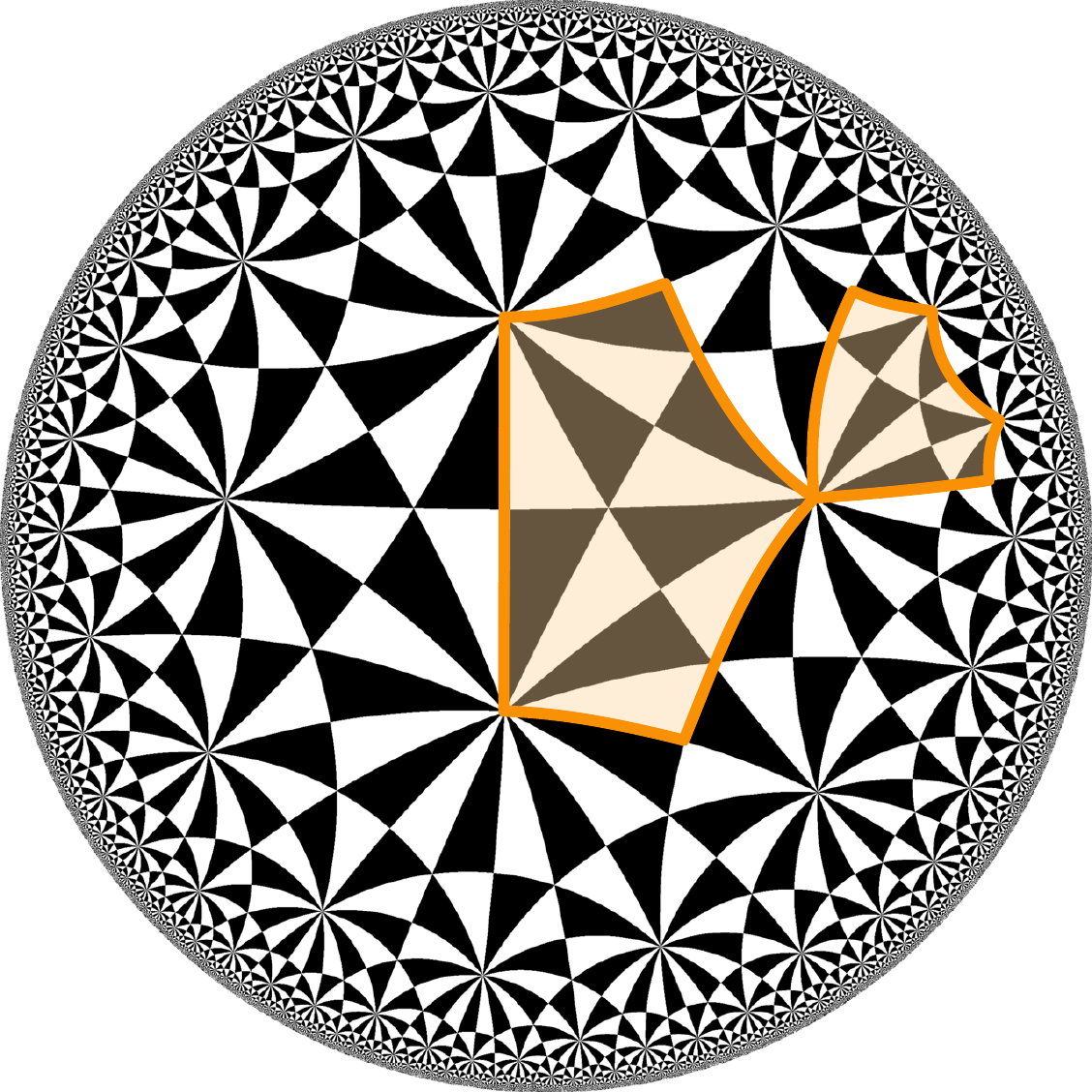}};
\end{tikzpicture}
\caption{A bow-tie fitting well on a $(2,3,8)$ tessellation of the hyperbolic plane obtained by reflecting a triangle with interior angles $(\pi/2, \pi/3, \pi/8)$ across its edges. The associated representation has $\beta=\pi$ and $\gamma=\pi/4$ (or $\gamma=3\pi/4$ depending on the parametrization), as well as $\tau_1=\tau_2=0$.}
\label{fig:discrete-bow-tie}
\end{figure}
\end{ex}

\begin{prop}\label{prop:gamma-is-dual-to-beta}
The angle $\gamma$ is the symplectic dual to the angle $\beta$ from Definition~\ref{defn:beta-angle}. In other words, if $\Phi_c^t$ denotes the Hamiltonian flow of the $\beta$-angle function $\beta\colon\mathcal{B}_c\to (0,2\pi)$, then
\[
\gamma\circ\Phi_c^t=\gamma + t.
\]
\end{prop}
\begin{proof}
The Hamiltonian flow $\Phi_c^t$ on $\mathcal{B}_c$---which we will call the \emph{twist flow} of $c$---descends from an explicit flow on representations that can be described as follows. Consider the element $Y(t)$ in the Lie algebra of $\psl$ given by
\[
Y(t)=\begin{pmatrix}
    0 & -t/2\\
    t/2 & 0
\end{pmatrix}.
\]
Its image under the Lie exponential map is the elliptic element
\[
\exp\big( Y(t)\big)=\pm\begin{pmatrix}
    \cos(t/2) & -\sin(t/2)\\
    \sin(t/2) & \cos(t/2)
\end{pmatrix}\in \psl.
\]
Geometrically, $\exp( Y(t))$ acts on the upper half-plane by a clockwise rotation of angle $t$ around the imaginary unit. In particular, for every $[\rho]\in\mathcal{B}_c$, there exists $B_\rho\in\psl$ such that
\[
\rho([a_1,b_1])=B_\rho \exp\big(Y(\beta(\rho))\big)B_\rho^{-1},
\]
because $\beta(\rho)$ was defined as the clockwise rotation angle of $\rho([a_1,b_1])$ (Definition~\ref{defn:beta-angle}). If we abbreviate $\zeta_\rho(t)=B_\rho \exp(Y(t))B_\rho^{-1}$, then $\{\zeta_\rho(t):t\in\R/2\pi\Z\}$ is the centralizer of $\rho([a_1,b_1])$ in $\psl$. Following~\cite{goldman-invariant-functions}, we define
\begin{equation}\label{eq:twist-flow-separating-curve-genus-2}
\overline{\Phi}_c^t(\rho)\colon\begin{cases}
    a_1\mapsto \rho(a_1)\\
    b_1\mapsto \rho(b_1)\\
    a_2\mapsto \zeta_\rho(t)\rho(a_2)\zeta_\rho(t)^{-1}\\
    b_2\mapsto \zeta_\rho(t)\rho(b_2)\zeta_\rho(t)^{-1}.
\end{cases}
\end{equation}
The flow $\overline{\Phi}_c^t$ is a lift at the level of representations of the twist flow $\Phi_c^t$. On the polygonal model of $\rho$, the flow $\overline{\Phi}_c^t$ does not affect the pentagon associated to the restriction $\rho_1$, but it rotates the pentagon of $\rho_2$ clockwise by an angle $t$ around their common vertex. By the definition of the angle $\gamma$, this increases its value by $t$.
\end{proof}

The angle $\gamma$ was the last missing coordinate to state our parametrization result for bow-tie representations. Like on Figure~\ref{fig:bow-tie}, we denote by $(\ell_1,\tau_1)$ and $(\ell_2,\tau_2)$ the Fenchel--Nielsen coordinates (Definitions~\ref{defn:length-coordinate-one-holed-torus} and~\ref{defn:twist-coordinate-one-holed-torus}) of the restrictions $[\rho_1]$ and $[\rho_2]$ of $[\rho]\in\mathcal{B}_c$. 

\begin{thm}\label{thm:coordinates-bow-tie-representations}
For every separating simple closed curve $c$ on $S$ and every compatible geometric presentation of $\pi_1S$, the map
\[
\mathfrak{T}\colon\mathcal{B}_c\to \R_{>0}\times \R\times (0,2\pi)\times \R/2\pi\Z \times \R_{>0}\times \R
\]
given by $\mathfrak{T}=(\ell_1,\tau_1,\beta,\gamma,\ell_2,\tau_2)$ is a smooth surjective map that restricts to a smooth diffeomorphism on both $\mathcal{B}_c^+$ and $\mathcal{B}_c^-$.
\end{thm}
\begin{proof}
The $\beta$-angle function $\beta\colon\mathcal{B}_c\to (0,2\pi)$ is smooth because it measures the rotation angle of $\rho([a_1,b_1])$. Similarly, the angle $\gamma\colon\mathcal{B}_c\to\R/2\pi\Z$ is also a smooth function, because it is the signed angle between two geodesic rays passing through points whose coordinates depend smoothly on $\rho$. So, since Fenchel--Nielsen coordinates are smooth, we conclude that $\mathfrak{T}$ is indeed a smooth map.

To prove the diffeomorphism property, we construct two inverse maps 
\[
\mathfrak{S}^{\pm}\colon \R_{>0}\times \R\times (0,2\pi)\times \R/2\pi\Z \times \R_{>0}\times \R\to\mathcal{B}_c^{\pm}.
\]
Assume that we are given a tuple $(\ell_1\tau_1,\beta,\gamma,\ell_2,\tau_2)$. First, we let $t=2\cos(\beta/2)$. Next, we consider the two points $[\rho_1^-], [\rho_1^+] \in \chi_{t}(\Sigma_1,\psl)$ with Fenchel--Nielsen coordinates $(\ell_1,\tau_1)$, and Toledo numbers $\Tol(\rho_1^-)=-\beta/2\pi$ and $\Tol(\rho_1^+)=\beta/2\pi$. 
Similarly, we consider the two points $[\rho_2^-], [\rho_2^+] \in \chi_{-t}(\Sigma_2,\psl)$ with Fenchel--Nielsen coordinates $(\ell_2,\tau_2)$, and Toledo numbers $\Tol(\rho_2^-)=\beta/2\pi-1$ and $\Tol(\rho_2^+)=1-\beta/2\pi$. These choices ensure that $\rho_1^+([a_1,b_1])$ and $\rho_2^+([a_2,b_2])$, as well as $\rho_1^-([a_1,b_1])$ and $\rho_2^-([a_2,b_2])$, are rotations of angle $\beta$ but in opposite directions.

Now, in each of the four conjugacy classes of representations we just described, there is exactly one representative $\rho_1^{\pm}$ and $\rho_2^{\pm}$ that satisfy the following conditions:
\begin{itemize}
    \item $\rho_1^{\pm}([a_1,b_1])$ and $\rho_2^{\pm}([a_2,b_2])$ all fix the imaginary unit in the upper half-plane.
    \item The vertices $Y_1'$ and $Y_2$ in the pentagon models of $\rho_1^{\pm}$ and $\rho_2^{\pm}$ both lie on the imaginary axis, above the imaginary unit.
\end{itemize}
We may now glue $\rho_1^+$ and $\rho_2^+$, as well as $\rho_1^-$ and $\rho_2^-$, along the boundary curve $[a_1,b_1]=[a_2,b_2]^{-1}$ to obtain two representations $\rho^+\colon\pi_1S\to\psl$ and $\rho^-\colon\pi_1S\to\psl$. Finally, we let 
\[
\mathfrak{S}^{\pm}(\ell_1,\tau_1,\beta,\gamma,\ell_2,\tau_2)=\Phi_c^\gamma([\rho^{\pm}]).
\]

By construction, $\eu(\rho^+)=\beta/2\pi+1-\beta/2\pi=1$ and $\eu(\rho^-)=-\beta/2\pi+\beta/2\pi-1=-1$, so $\mathfrak{S}^+$ and $\mathfrak{S}^-$ do take image in $\mathcal{B}_c^+$, respectively $\mathcal{B}_c^-$. Since the twist flow $\Phi_c^t$ depends smoothly on $t$, in order to prove that $\mathfrak{S}^{\pm}$ are smooth functions, it is sufficient to prove that $[\rho^{\pm}]$ depend smoothly on $(\ell_1,\tau_1,\beta,\ell_2,\tau_2)$. However, the maps $(\ell_i,\tau_i,\beta)\mapsto [\rho_i^{\pm}]$ are smooth because they are the inverses of the Fenchel--Nielsen coordinate maps for one-holed tori. So, it remains to prove that the gluing procedure $([\rho_1^{\pm}],[\rho_2^{\pm}])\mapsto [\rho^{\pm}]$ defined above is also smooth. This, however, is true because our choice of representatives for $[\rho_i^{\pm}]$ defines a local embedding from the components of relative character varieties containing $[\rho_i^{\pm}]$ into the space of representations $\pi_1\Sigma_i\to\psl$.

To prove that $\mathfrak{T}\circ\mathfrak{S}^{\pm}=\id$, we will denote the six components of $\mathfrak{T}$ by $\mathfrak{T}_1,\ldots,\mathfrak{T}_6$. It follows immediately from the definition of $\mathfrak{S}^{\pm}$, that $\mathfrak{T}_i\circ\mathfrak{S}^{\pm}=\id$ for all $i\neq 4$. To see that $\mathfrak{T}_4\circ\mathfrak{S}^{\pm}=\id$, note that the representations $\rho^{\pm}$ appearing in the definition of $\mathfrak{S}^{\pm}$ satisfy $\gamma(\rho^{\pm})=0$. So, $\gamma(\Phi_c^t([\rho^{\pm}])=t$ for every $t$ by Proposition~\ref{prop:gamma-is-dual-to-beta}. This proves $\mathfrak{T}_4\circ\mathfrak{S}^{\pm}=\id$.

Finally, we need to prove that $\mathfrak{S}^{\pm}\circ\mathfrak{T}=\id$. The argument is similar to the one we used in the proof of Theorem~\ref{thm:coordinates-pentagon-representations}. First, observe that $[\rho]$ and $\mathfrak{S}^{\pm}\circ\mathfrak{T}([\rho])$ have the same Fenchel--Nielsen coordinates when restricted to $\pi_1\Sigma_1$ and $\pi_1\Sigma_2$. Furthermore, the restrictions belong to the same components of the relative character varieties $\chi_{\pm t}(\Sigma_1,\psl)$ and $\chi_{\mp t}(\Sigma_2,\psl)$ for $t=2\cos(\beta([\rho])/2)$. So, the restrictions are equal. By the geometric interpretation of the function $\gamma$ illustrated on Figure~\ref{fig:bow-tie} and Proposition~\ref{prop:gamma-is-dual-to-beta}, we conclude that the restrictions glue back to the same point and thus $[\rho]=\mathfrak{S}^{\pm}\circ\mathfrak{T}([\rho])$.
\end{proof}

\subsection{Symplectic structure}\label{sec:symplectic-structure-bow-tie-representations}
The coordinates for bow-tie representations introduced in Theorem~\ref{thm:coordinates-bow-tie-representations} turn out to be Darboux coordinates for the Goldman symplectic form on $\mathcal{B}_c$, leading to a Wolpert-like formula.
\begin{thm}\label{thm:wolpert-formula-bow-tie-representations}
If $\omega_\mathcal{G}$ denotes the Goldman form on $\mathcal{B}_c$, then 
\begin{equation}\label{eq:wolpert-formula-bow-tie-representations}
\omega_{\mathcal{G}}=-d\ell_1\wedge d\tau_1 - d\beta\wedge d\gamma - d\ell_2\wedge d\tau_2.
\end{equation}
\end{thm}
\begin{proof}
As in the proof of Theorem~\ref{thm:wolpert-formula-pentagon-representations}, we will make use of the symplectic involution $\jmath\colon\chi(S,\psl)\to\chi(S,\psl)$ given by conjugation by a matrix of determinant $-1$ and introduced in~\eqref{eq:involution-j}. Recall that $\jmath$ maps $\mathcal{B}_c^-$ isomorphically onto $\mathcal{B}_c^+$. For the same reasons as in the proof of Theorem~\ref{thm:coordinates-pentagon-representations}, it holds that $\jmath^\ast d\ell_i=d\ell_i$ and $\jmath^\ast d\tau_i=d\tau_i$. Furthermore, since $\jmath$ maps clockwise rotations of angle $\beta$ in $\psl$ to anti-clockwise rotations of angle $\beta$, we have $\jmath^\ast d\beta=-d\beta$ by definition of $\beta$ (Definition~\ref{defn:beta-angle}). Finally, since $\jmath$ can be realized on polygonal models by a reflection through a geodesic line, it follows that $\jmath^\ast d\gamma=-d\gamma$ by definition of $\gamma$ (Definition~\ref{defn:gamma-angle-bow-tie-representation}). These relations imply that if~\eqref{eq:wolpert-formula-bow-tie-representations} holds on $\mathcal{B}_c^-$, then it also holds on $\mathcal{B}_c^+$.

We now prove that~\eqref{eq:wolpert-formula-bow-tie-representations} holds on $\mathcal{B}_c^-$. The beginning of the argument is similar to the proof of Theorem~\ref{thm:wolpert-formula-pentagon-representations}, so we will omit certain details. Since the simple closed curves represented by $a_1$, $[a_1,b_1]$, and $a_2$ are disjoint, it follows that $\{\ell_i,\ell_j\}\equiv 0$ and $\{\ell_i,\beta\}\equiv 0$, where $\{\cdot,\cdot\}$ stand for the Poisson bracket associated to $\omega_\mathcal{G}$, see~\eqref{eq:Poisson-bracket}. Furthermore, combining~\eqref{eq:Goldman-formula-one-holed-torus} and Proposition~\ref{prop:gamma-is-dual-to-beta}, as well as the explicit forms of the twist flow along $c$ provided by~\eqref{eq:twist-flows-genus-2}, we conclude that $\{\ell_i,\tau_j\}\equiv-\delta_{ij}$, $\{\ell_i,\gamma\}=\{\beta,\tau_i\}\equiv 0$, and $\{\beta,\gamma\}\equiv -1$. This means that when we express $\omega_{\mathcal{G}}$ in the coordinates provided by Theorem~\ref{thm:coordinates-bow-tie-representations}, we may write
\begin{align}
    \omega_{\mathcal{G}} &= -d\ell_1\wedge d\tau_1 - d\beta\wedge d\gamma - d\ell_2\wedge d\tau_2 \label{eq:symplectic-form-intermediate-form-bow-tie}\\
    & + \alpha_0\, d\ell_1\wedge d\ell_2 + \alpha_1\, d\ell_1\wedge d\beta + \alpha_2\, d\ell_2\wedge d\beta, \nonumber
\end{align}
where the coefficients $\alpha_0,\alpha_1,\alpha_2$ are smooth functions. 

We want to prove that each $\alpha_i$ is identically zero. As in the proof of Theorem~\ref{thm:wolpert-formula-pentagon-representations}, since the coefficients $\alpha_i$ are independent of the twist parameters, it is sufficient to prove that $\alpha_i\equiv 0$ for specific values of the twist parameters. We will thus only consider the points $[\rho]\in\mathcal{B}_c^-$ with $\tau_1=\tau_2=0$ and $\gamma=\beta/2$. Those points have a special symmetry. Consider the involutive automorphism of $\pi_1S$ given by
\[
\imath\colon\begin{cases}
a_1\mapsto a_1^{-1}\\
b_1\mapsto a_1b_1a_1^{-1}\\
a_2\mapsto a_2^{-1}\\
b_2\mapsto a_2b_2a_2^{-1}.
\end{cases}
\]
The outer automorphism of $\pi_1S$ represented by $\imath$ corresponds to the isotopy class of an orientation-reversing homeomorphism of $S$, hence $\imath^\ast\omega_{\mathcal{G}}=-\omega_{\mathcal{G}}$. In general, $\imath$ and $\jmath$ act differently on $\mathcal{B}_c$.
\begin{claim}\label{claim:involution-representation-bow-tie}
However, if $[\rho]\in\mathcal{B}_c^-$ has $\tau_1=\tau_2=0$ and $\gamma=\beta/2$, then
\[
[\rho\circ\imath]=\jmath[\rho].
\]
\end{claim}
\begin{proofclaim}
To prove that $[\rho\circ\imath]=\jmath[\rho]$, we will show that the associated polygonal models (bow-ties) are isometric with matching orientation. We start from the bow-tie associated to $[\rho]$. Recall from the discussion that led to Figure~\ref{fig:bow-tie} that it consists of the two pentagons $X_1Y_1BY_1'X_1'$ and $X_2Y_2BY_2'X_2'$ linked at a common vertex $B$. The two pentagons are the polygonal model of the restrictions of $[\rho]$ to $\pi_1\Sigma_1$ and $\pi_1\Sigma_2$. Since we are assuming that $[\rho]$ satisfies $\tau_1=\tau_2=0$, both pentagons are right-angled pentagons. Since we are also assuming that $\gamma=\beta/2$, the points $Y_1$, $B$, and $Y_2$ lie on a common geodesic line, with $B$
between $Y_1$ and $Y_2$, see Figure~\ref{fig:bow-tie-reflected}. We label $g$ the geodesic line $(Y_1Y_2)$ and $\sigma$ the orientation-reversing reflection across $g$.

\begin{figure}[h]
\centering
\begin{tikzpicture}
\node[anchor = south west, inner sep=0mm] at (0,0) {\includegraphics[width=13cm]{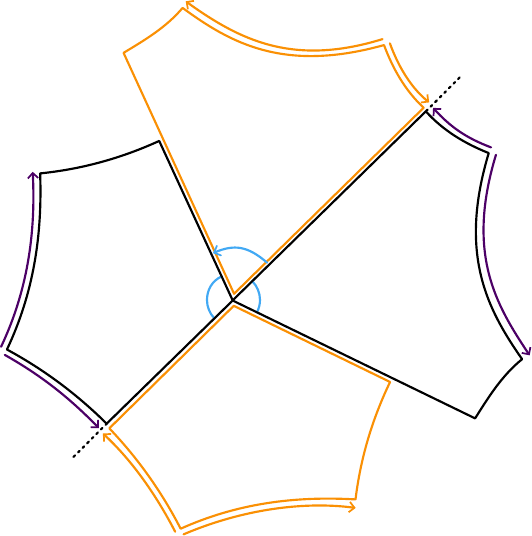}};

\node[sky] at (4.3,5.8) {$\pi-\beta/2$};
\node[sky] at (6.8,5.8) {$\beta/2$};
\node[sky] at (6,7.3) {$\gamma=\beta/2$};

\node[plum] at (.9,3.4) {$\rho(a_1)$};
\node[plum] at (.2,6.9) {$\rho(b_1)$};
\node[plum] at (11.7,10.2) {$\rho(a_2)$};
\node[plum] at (12.4,7) {$\rho(b_2)$};

\node at (-.2,4.5) {$X_1$};
\node at (2.6,3.15) {$Y_1$};
\node at (5.75,5.35) {$B$};
\node at (3.55,9.85) {$Y_1'$};
\node at (.9,9.15) {$X_1'$};
\node at (11.65, 2.6) {$Y_2'$};
\node at (13, 4.15) {$X_2'$};
\node at (12.3, 9.5) {$X_2$};
\node at (10.45, 9.95) {$Y_2$};

\node[apricot] at (4.2,-.2) {$\sigma(X_1)$};
\node[apricot] at (9.3,.9) {$\sigma(X_1')$};
\node[apricot] at (8.9,3.6) {$\sigma(Y_1')$};
\node[apricot] at (2.5,11.8) {$\sigma(Y_2')$};
\node[apricot] at (4,13.15) {$\sigma(X_2')$};
\node[apricot] at (9.95,12.25) {$\sigma(X_2)$};

\node at (1.7,2.1) {$g$};
\end{tikzpicture}
\caption{In black: the bow-tie associated to $[\rho]$, made of the the two right-angled pentagons $X_1Y_1BY_1'X_1'$ and $X_2Y_2BY_2'X_2'$. In orange: its image by the reflection $\sigma$ across the doted geodesic line $g$. For a better visualization, the two orange pentagons have been moved a little bit, but they are really linked at $B$.}
\label{fig:bow-tie-reflected}
\end{figure}

The bow-tie of $\jmath[\rho]$ is obtained by applying $\sigma$ to the bow-tie of $[\rho]$ resulting in the orange bow-tie illustrated on Figure~\ref{fig:bow-tie-reflected}. To prove the claim, we will show that it is also the bow-tie associated to $[\rho\circ\imath]$. Recall from the construction of the bow-tie associated to $\rho$ that $\rho(a_1)$ is the hyperbolic translation along $(X_1Y_1)$ of length $2d(X_1,Y_1)$. Since the geodesic lines $(X_1Y_1)$ and $g$ are perpendicular, it follows that $\rho(a_1)$ sends $X_1$ to $\sigma(X_1)$. By definition of $\imath$, $\rho\circ\imath(a_1)=\rho(a_1)^{-1}$ is the hyperbolic translation along $(\sigma(X_1)Y_1)$ of length $2d(X_1,Y_1)=2d(\sigma(X_1),Y_1)$. Moreover, $\rho\circ\imath(b_2)=\rho(a_1)\rho(b_1)\rho(a_1)^{-1}$ is a hyperbolic translation whose axis is the image by $\rho(a_1)$ of the geodesic line $(X_1X_1')$---the axis of $\rho(b_1)$. Since we are assuming $\tau_1=0$, the geodesic lines $(X_1X_1')$ and $(X_1Y_1)$ are perpendicular, so $\rho(a_1)$ maps $(X_1X_1')$ to the geodesic line $(\sigma(X_1)\sigma(X_1'))$. Therefore, $\rho\circ\imath(b_2)$ is a hyperbolic translation along $(\sigma(X_1),\sigma(X_1'))$ of length $d(X_1,X_1')=d(\sigma(X_1),\sigma(X_1'))$. This shows that the polygonal model for the restriction of $[\rho\circ\imath]$ to $\pi_1\Sigma_1$ is the pentagon $\sigma(X_1)Y_1B\sigma(Y_1')\sigma(X_1')$, which is the image of $X_1Y_1BY_1'X_1'$ under $\sigma$. The analogous conclusion holds for the restriction of $[\rho\circ\imath]$ to $\pi_1\Sigma_2$ and it finishes the proof of the claim.
\end{proofclaim}

Claim~\ref{claim:involution-representation-bow-tie} implies that $(\imath^\ast\omega_{\mathcal{G}})_{[\rho]}=-(\jmath^\ast\omega_{\mathcal{G}})_{[\rho]}$ for every $[\rho]\in\mathcal{B}_c^-$ with $\tau_1=\tau_2=0$ and $\gamma=\beta/2$ by the same computation as in the proof of Theorem~\ref{thm:wolpert-formula-pentagon-representations}. Furthermore, observe that $\imath^\ast d\ell_i=d\ell_i$ because $\rho\circ\imath(a_i)=\rho(a_i)^{-1}$. On the other hand, $\imath^\ast d\beta=-d\beta$, because $\rho\circ\imath([a_1,b_1])=\rho([b_1,a_1])=\rho([a_1,b_1])^{-1}$ and $\beta$ is the clockwise rotation angle of $\rho([a_1,b_1])$. For the same reasons as in the proof of Theorem~\ref{thm:wolpert-formula-pentagon-representations}, it also holds that $\imath^\ast d\tau_i=-d\tau_i$. 
\begin{claim}
Moreover, we have $\imath^\ast d\gamma=d\gamma-d\beta$.  
\end{claim}
\begin{proofclaim}
We will study the action of $\imath$ on the two pentagons associated to the restrictions $[\rho_1]$ and $[\rho_2]$ of $[\rho]\in\mathcal{B}_c^-$. We start by observing that for every $t\in \R$, we have
\[
\imath=\overline{\Phi}_{a_1}^{t}\circ\imath\circ\overline{\Phi}_{a_1}^{t},
\]
where $\overline{\Phi}^t_{a_1}$ is the Hamiltonian twist flow along $a_1$. This identity is a consequence of the explicit formula of $\overline{\Phi}^t_{a_1}$ given in~\eqref{eq:twist-flows-genus-2} and the definition of $\imath$. So, in order to describe the action of $\imath$ on the pentagon associated to $\rho_1$, we perform the following steps.
\begin{enumerate}
    \item Start from the pentagon $X_1Y_1BY_1'X_1'$ of $\rho_1$.
    \item Apply $\overline{\Phi}_{a_1}^{t}$ with $t=-\tau_1([\rho])$ to transform it into a right-angled pentagon. As explained in Proposition~\ref{prop:twist-coordinate-is-geometric-distance}, this step only moves the points $X_1$ and $X_1'$.
    \item Reflect it across the side $(Y_1B)$ as we did in the proof of Claim~\ref{claim:involution-representation-bow-tie}. The resulting pentagon is the one associated to $\imath\circ \overline{\Phi}_{a_1}^{t}(\rho_1)$.
    \item Apply $\overline{\Phi}_{a_1}^{t}$ again. This last step only moves the images of $X_1$ and $X_1'$. 
\end{enumerate}
Following those steps, we conclude that the pentagon associated to $[\rho_1]$ is of the form $\overline{X}_1Y_1B\sigma_1(Y_1')\overline{X}_1'$, where $\sigma_1$ is the reflection across the geodesic line $(Y_1B)$. The points $\overline{X}_1$ and $\overline{X}_1'$ are actually the images of $X_1$ and $X_1'$ under $\rho(a_1)$ since $\rho\circ\imath(b_1)=\rho(a_1)\rho(b_1)\rho(a_1)^{-1}$. A similar argument shows that $\imath$ transforms the pentagon of $\rho_2$ into $\overline{X}_2Y_2B\sigma_2(Y_2')\overline{X}_2'$, where $\sigma_2$ is the reflection across the geodesic line $(Y_2B)$, and $\overline{X}_2=\rho(a_2)X_2$, $\overline{X}_2'=\rho(a_2)X_2'$. 

We have just described the bow-tie associated to $[\rho\circ\imath]$. In particular, we conclude that $\gamma(\rho\circ\imath)$ is the oriented angle from the geodesic ray $[BY_2)$ to the geodesic ray $[B\sigma_1(Y_1'))$ using Definition~\ref{defn:gamma-angle-bow-tie-representation}. Decomposing this angle as the sum of the oriented angle from $[BY_2)$ to $[BY_1')$ plus the angle from $[BY_1')$ to $[B\sigma_1(Y_1'))$, we obtain $\gamma(\rho\circ\imath)=\gamma+2\pi-\beta$. Hence $\imath^\ast d\gamma = d\gamma-d\beta$.
\end{proofclaim}

We can now conclude the proof of Theorem~\ref{thm:wolpert-formula-bow-tie-representations}. We inject the expression~\eqref{eq:symplectic-form-intermediate-form-bow-tie} for $\omega_{\mathcal{G}}$ into the relation $(\imath^\ast\omega_{\mathcal{G}})_{[\rho]}=-(\jmath^\ast\omega_{\mathcal{G}})_{[\rho]}$ and use the pullback relations we computed. This gives
\begin{align*}
(\imath^\ast\omega_{\mathcal{G}})_{[\rho]}&=d\ell_1\wedge d\tau_1+d\beta\wedge d\gamma+d\ell_2\wedge d\tau_2\\
& +\alpha_0([\rho])d\ell_1\wedge d\ell_2 - \alpha_1([\rho]) d\ell_1\wedge d\beta -\alpha_2([\rho]) d\ell_2\wedge d\beta,\\
-(\jmath^\ast\omega_{\mathcal{G}})_{[\rho]}&=d\ell_1\wedge d\tau_1+d\beta\wedge d\gamma+d\ell_2\wedge d\tau_2\\
& -\alpha_0([\rho])d\ell_1\wedge d\ell_2 + \alpha_1([\rho]) d\ell_1\wedge d\beta +\alpha_2([\rho]) d\ell_2\wedge d\beta.
\end{align*}
This calculation implies $\alpha_i([\rho])= 0$ for $i=1,2,3$ as desired.
\end{proof}

\subsection{Geometrization}\label{sec:geometrization-bow-tie-representations}
The goal of this section is to prove that all bow-tie representations are geometrizable in sense of Definition~\ref{defn:geometrizable-representations}. More precisely, for every bow-tie representation $\rho\colon\pi_1S\to\psl$, we will construct a 2-dimensional family of branched hyperbolic structures on $S$ having $\rho$ as holonomy. Our construction is inspired by the samosa assemblies from~\cite{fenyes-maret}.

\subsubsection{Building blocks}
We will construct branched hyperbolic structures using a cut-and-paste procedure on two hyperbolic tori with complementary conical singularities. Recall that a \emph{hyperbolic cone structure} on a closed and oriented topological torus $\Sigma$ with one marked point $p$ is the data of a hyperbolic structure on $\Sigma\setminus \{p\}$ whose metric completion at $p$ is a cone of angle $\beta>0$. We call $(p,\beta)$ the \emph{cone data}. By the Gauss--Bonnet theorem, the angle $\beta$ must be less than $2\pi$. When $\beta$ approaches $0$, the cone singularity turns into a cusp.

\begin{figure}[h]
\centering
\begin{tikzpicture}
\node[anchor = south west, inner sep=0mm] at (0,0) {\includegraphics[width=6cm]{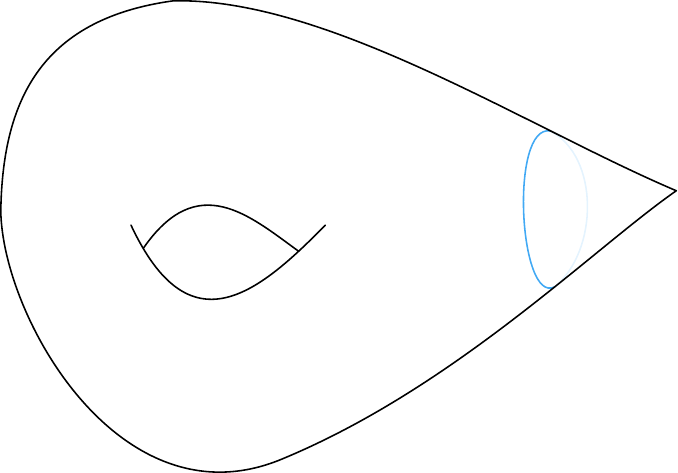}};
\node[sky] at (4.3,2.3) {$\beta$};
\end{tikzpicture}
\caption{A hyperbolic cone torus with cone angle $\beta$.}
\label{fig:hyperbolic-cone-torus}
\end{figure}

\begin{figure}[h]
\centering
\begin{tikzpicture}
\node[anchor = south west, inner sep=0mm] at (0,0) {\includegraphics[width=6cm]{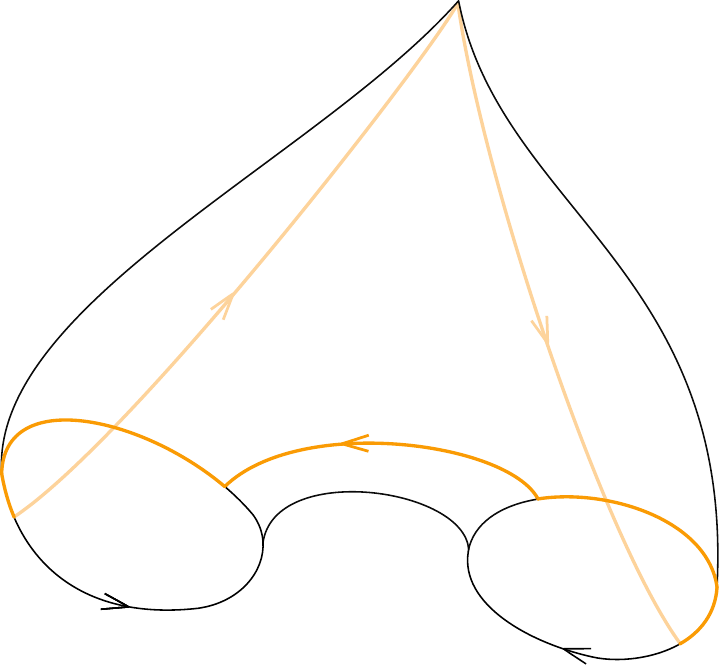}};

\node at (3.9, 5.8) {$p$};
\node at (1,0.2) {$\alpha_j$};
\node at (5.3,-.2) {$\alpha_i$};
\node[apricot] at (3.2,2.1) {$\gamma$};
\node[lightapricot] at (1.2,2.7) {$\xi_j$};
\node[lightapricot] at (5.1,2.2) {$\xi_i$};

\draw (-1.5, 0) edge[out=90, in=225, ->] (.7,1);
\node[anchor=north] at (-1.5,0) {\text{southern hemitorus}};
\draw (-1.5, 4) edge[out=-90, in=135, ->] (2,3.5);
\node[anchor=south] at (-1.5,4) {\text{northern hemitorus}};
\end{tikzpicture}
\caption{In orange, the equator of a hyperbolic cone torus.}
\label{fig:equator}
\end{figure}

Given a simple closed curve $a$ on $\Sigma$ and a hyperbolic cone structure on $\Sigma$ with cone point $p$, there is a unique simple closed geodesic $\alpha$ representing $a$. This geodesic always avoids the cone singularity since the cone angle is less than $2\pi$. Given two oriented simple closed curves $a$ and $b$ on $\Sigma$ that intersect once, we will refer to $(\Sigma,a,b)$ as a \emph{framed hyperbolic cone torus}. We can use the cone point $p$ of a framed hyperbolic torus $(\Sigma,a,b)$ to define the \emph{equator} of $\Sigma$, which we illustrate on Figure~\ref{fig:equator}. To do so, first cut $\Sigma$ along the geodesic representative $\alpha$ of $a$ as we did on Figure~\ref{fig:one-holed-torus-cut-and-uncut}. This produces a hyperbolic pair of pants with one cone singularity and two geodesic boundary $\alpha_1$ and $\alpha_2$ of equal lengths. There are two shortest geodesic arcs $\xi_1$ and $\xi_2$ from $p$ to $\alpha_1$, respectively to $\alpha_2$. Recall that, in the context of Figure~\ref{fig:one-holed-torus-cut-and-uncut}, we denoted by $\gamma$ the shortest geodesic arc between $\alpha_1$ and $\alpha_2$. The geodesic representative of $b$ is oriented and intersects $\alpha$ once. It therefore travels from $\alpha_i$ to $\alpha_j$ for some ordering of $\{i,j\}=\{1,2\}$. We orient $\gamma$ so that it also travels from $\alpha_i$ to $\alpha_j$. 

\begin{defn}\label{def:equator}
The \emph{equator} of $(\Sigma,a,b)$ is the oriented loop on $\Sigma$ illustrated on Figure~\ref{fig:equator} and obtained by concatenating the following arcs:
\begin{itemize}
    \item First, travel positively along $\gamma$;
    \item Then, travel positively along $\alpha_j$ until reaching the endpoint of $\xi_j$;
    \item Then, travel along $\xi_j$ to $p$ and along $\xi_i$ to $\alpha_i$;
    \item Finally, travel negatively along $\alpha_i$ until reaching the starting point of $\gamma$ to close the loop.
\end{itemize}
The equator separates $\Sigma$ into two \emph{hemitori}. Using the orientations of $\Sigma$ and of its equator, we define the \emph{northern hemitorus} as the one to the right of the equator, and the \emph{southern hemitorus} as the one to the left.
\end{defn}

\begin{figure}[h]
\centering
\begin{tikzpicture}
\node[anchor = south west, inner sep=0mm] at (0,0) {\includegraphics[width=6cm]{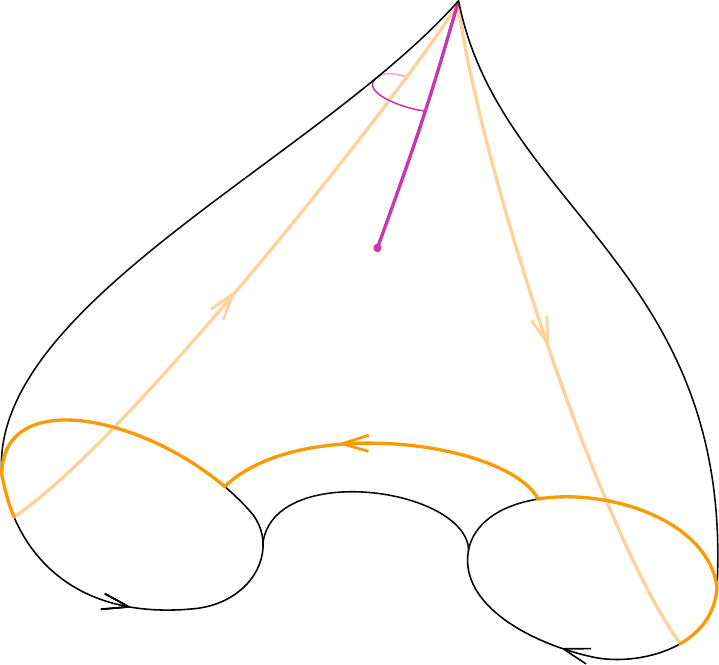}};
\node at (3.9, 5.8) {$p$};
\node[mauve] at (2.9, 5.1) {$\phi$};
\node[mauve] at (3.6, 4.3) {$\ell$};
\end{tikzpicture}
\caption{In mauve, a slit on a hyperbolic cone torus.}
\label{fig:slit}
\end{figure}

\begin{defn}\label{defn:slit}
A \emph{slit} on a framed hyperbolic cone torus $(\Sigma,a,b)$ with cone point $p$ is a simple geodesic arc that starts at $p$, never meets the geodesic representatives of $a$ and $b$, and remains in the closure of one hemitorus. An example of a slit is illustrated on Figure~\ref{fig:slit}. A \emph{slit torus} is a framed hyperbolic cone torus with one slit at its cone singularity.
\end{defn}

A slit is parametrized by its length $\ell>0$ and the angle $\phi\in \R/\beta \Z$ between the equator and the slit. We always measure $\phi$ from the negative side of the equator, so that $\phi=0$ when the slit runs along the equator in the negative direction and $\phi=\beta/2$ when the slit runs along the equator in the positive direction. If we think of the slit angle as a number in $(-\beta/2, \beta/2]$, then $\phi\geq 0$ if the slit lies in the northern hemitorus, and $\phi\leq 0$ if it lies in the southern hemitorus.

\subsubsection{Assembling slit tori}
\begin{figure}[h]
\centering
\begin{tikzpicture}
\node[anchor = south west, inner sep=0mm] at (0,0) {\includegraphics[width=13cm]{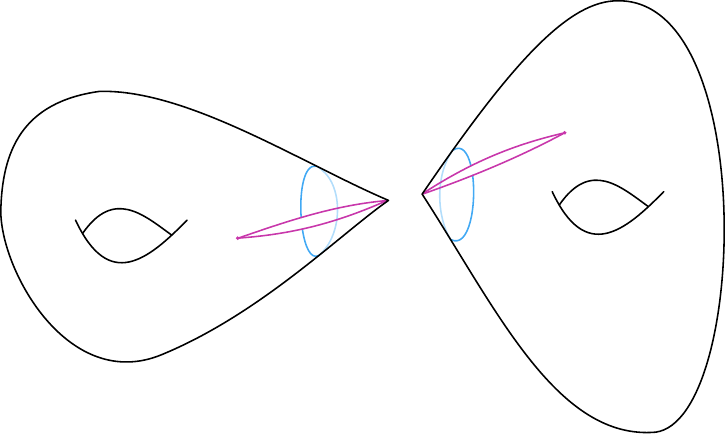}};
\node[sky] at (5.1,4.3) {$\beta_1$};
\node[sky] at (8.8,4.3) {$\beta_2$};
\draw[mauve] (5.8,4.05) edge[out=80, in=130, <->] (9.2,5.18);
\draw[mauve] (5.8,3.7) edge[out=-80, in=-70, <->] (9.3,4.91);
\end{tikzpicture}
\caption{The cut-and-paste procedure to turn two cone singularities with angles $\beta_1$ and $\beta_2$ into a regular point and a branch point.}
\label{fig:cut-and-paste}
\end{figure}

Whenever we are given two slit tori with cone data $(p_1,\beta_1)$ and $(p_2,\beta_2)$ that are \emph{complementary}, meaning that $\beta_1+\beta_2=2\pi$, we may apply a cut-and-paste procedure to turn them into a branched hyperbolic surface of genus two. To do so, we need the two slits to have the same length $\ell >0$. The procedure works as follows. First, cut along the two slits from their starting point to their endpoint. The incised slits have a right and a left \emph{lip}. Next, glue the two tori together by identifying each slit's left lip to the other's slit right lip, as we illustrate on Figure~\ref{fig:cut-and-paste}. The starting points of the two slits are identified in the glued surface to a point with cone angle $\beta_1+\beta_2=2\pi$, hence a regular point. The endpoints of the slits, on the other hand, become a cone singularity of angle $2\pi+2\pi=4\pi$, hence a branch point.

\begin{defn}\label{defn:slit-tori-assemblies}
We denote by $\mathcal{S}$ the space of \emph{slit tori assemblies} which consist of the following data:
\begin{itemize}
    \item Two framed hyperbolic cone tori $(\Sigma_1, a_1, b_1)$ and $(\Sigma_2, a_2, b_2)$ with respective cone data $(p_1,\beta_1)$ and $(p_2, \beta_2)$ satisfying $\beta_2=2\pi-\beta_1$.
    \item Two slits, one on $\Sigma_1$ and one on $\Sigma_2$, with the same length $\ell>0$ and respective slit angles $\phi_1\in \R/\beta_1\Z$ and $\phi_2\in \R/\beta_2\Z$.
\end{itemize}
We will also denote by $\mathcal{S}^\textrm{north}$ the subset of all slit tori assemblies in which both slits lie in the closures of the respective northern hemitori. In terms of our parameters, $\mathcal{S}^\textrm{north}$ is characterized by $\phi_1\in [0, \beta_1/2]$ and $\phi_2\in [0,\beta_2/2]$, as well $\ell>0$ being sufficiently small so that the slits never leave the hemitorus they started in.
\end{defn}

For coherence, we want all the branched hyperbolic structures produced by the cut-and-paste method described above to be attached to the ``same'' topological surface of genus $2$. In other words, we want all branched hyperbolic structures coming from the gluing of a slit tori assembly to be homeomorphic, up to isotopy, in a canonical way. We can achieve this as follows. 

\begin{figure}[h]
\centering
\begin{tikzpicture}
\node[anchor = south west, inner sep=0mm] at (0,0) {\includegraphics[width=8cm]{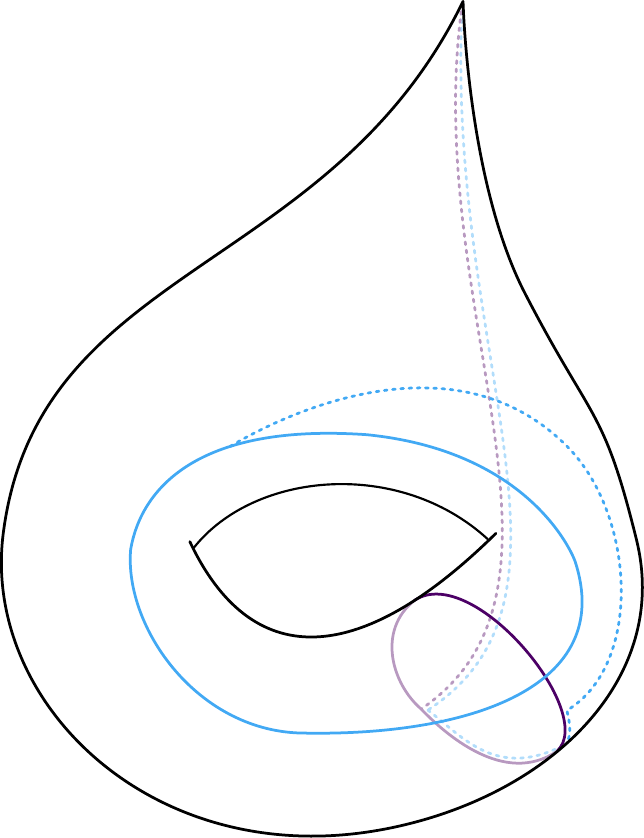}};
\node at (5.78,10.7) {$p$};
\end{tikzpicture}
\caption{The loop $\mathfrak{a}$ in purple and the loop $\mathfrak{b}$ in blue. The dotted lines indicate the equator segments and the unbroken lines indicate the geodesic representatives of $a$ and~$b$.}
\label{fig:curves-a-b}
\end{figure}

Let $(\Sigma,a,b)$ be a framed hyperbolic cone torus with cone data $(p,\beta)$. A slit on $\Sigma$ always entirely lies in the closure of one hemitorus by Definition~\ref{defn:slit}. 
We will consider the case where slits lie in the closure of the northern hemitorus. We can turn the geodesic representatives of $a$ and $b$ into loops $\mathfrak{a}$ and $\mathfrak{b}$ based at $p$ by traveling positively along the equator from $p$ until we first encounter those geodesics, then looping along the geodesics once, and coming back to $p$ by traveling negatively along the equator. We illustrate the resulting loops on Figure~\ref{fig:curves-a-b}. Note that the loops $\mathfrak{a}$ and $\mathfrak{b}$ always avoid meeting any slit in the northern hemitorus (except at $p$) as long as the slit does not run positively along the equator (i.e.~$\phi=\beta/2$). When $\phi=\beta/2$, we push the initial and final segments of $\mathfrak{a}$ and $\mathfrak{b}$ into the southern hemisphere to avoid overlapping the slit. The homotopy classes of $\mathfrak{a}$ and $\mathfrak{b}$ are generators of $\pi_1(\Sigma,p)$, which we abusively denote by $a$ and $b$. Note that the loop $[\mathfrak{a},\mathfrak{b}]$ is homotopic relative to $p$ to the equator of $\Sigma$ with matching orientations by Definition~\ref{def:equator}. When we cut $\Sigma$ along a slit in its norther hemisphere, the loop based at $p$ obtained by following the left lip of the slit and going back to $p$ along the right lip is also homotopic relative to $p$ to both $[\mathfrak{a},\mathfrak{b}]$ and the equator of $\Sigma$.

Now, let $(\Sigma_1,a_1)$ and $(\Sigma_2, a_2)$ be two framed hyperbolic cone tori with conical data $(p_1,\beta_1)$ and $(p_2,\beta_2=2\pi-\beta_1)$, and with slits in the closures of their northern hemitori. Denote by $\mathfrak{a}_1,\mathfrak{b}_1 \subset \Sigma_1$ and $\mathfrak{a}_2,\mathfrak{b}_2\subset \Sigma_2$ the loops based at $p_1$ and $p_2$ we just described. The cut-and-paste procedure applied to $\Sigma_1$ and $\Sigma_2$ produces a branched hyperbolic surface $S$ of genus $2$. The endpoints of the slits turn into the branch point $q$ of $S$, while the cone points $p_1$ and $p_2$ turn into a regular point $p$ of $S$. Since the loops $\mathfrak{a}_i,\mathfrak{b}_i$ avoid meeting the slits (except at the basepoint) by construction, they are also loops on $S$ based at $p$. 
\begin{claim}
The loop $[\mathfrak{a}_1,\mathfrak{b}_1][\mathfrak{a}_2,\mathfrak{b}_2]$ is homotopically trivial on $S$.
\end{claim}
\begin{proofclaim}
By construction, $[\mathfrak{a}_i,\mathfrak{b}_i]$ is homotopic relative to $p_i$ to the concatenation of the slit's left lip with slit's right lip on $\Sigma_i$. Since our gluing procedure identifies a slit's left lip with the other slit's right lip and vice versa, the loop $[\mathfrak{a}_1,\mathfrak{b}_1][\mathfrak{a}_2,\mathfrak{b}_2]$ is indeed homotopically trivial on $S$.
\end{proofclaim}
There is a geodesic arc $\mathfrak{s}$ on $S$ going from $q$ to $p$ that is the gluing of the $\Sigma_1$ slit's left lip with the $\Sigma_2$ slit's right lip. Concatenating $\mathfrak{s}$ and the loops $\mathfrak{a}_i,\mathfrak{b}_i$, we subsequently obtained a preferred geometric presentation of $\pi_1(S,q)$. In conclusion, every slit tori assembly in $\mathcal{S}^\mathrm{north}$ gives a preferred presentation of the fundamental group for the glued branched hyperbolic surface, hence defining an isotopy class of homeomorphisms between any two glued surfaces in this way.

In other words, cutting and gluing slit tori assemblies defines a realization map
\begin{equation}
    R\colon \mathcal{S}^\textrm{north}\to \Hyp_1(S),
\end{equation}
where $\Hyp_1(S)$ denotes the space of isotopy classes of branched hyperbolic structures on a closed surface $S$ of genus $2$ with one branch point, as introduced in Section~\ref{sec:branched-hyperbolic-structures}.

\begin{prop}\label{prop:bow-tie-realization-map-injective}
The realization map $R\colon \mathcal{S}^\textrm{north}\to \Hyp_1(S)$ is injective.
\end{prop}
\begin{proof}
We will construct an inverse map $R^{-1}\colon R(\mathcal{S}^\textrm{north})\to \mathcal{S}^\textrm{north}$ defined on the image of~$R$. 

We start with following observation. Let $\Sigma$ be a slit torus with cone data $(p,\beta)$. We denote by $q'$ the endpoint of the slit and we assume that $\beta\leq \pi$. When we cut $\Sigma$ along its slit, the slit's two lips form a broken geodesic arc that starts at $q'$ and comes back to $q'$. When $\beta< \pi$, the isotopy class of that arc relative to $q'$ contains a unique unbroken geodesic arc representative $\zeta$ based at $q'$, as the one we illustrate on Figure~\ref{fig:kinked-geodesic}. When $\beta=\pi$, the arc $\zeta$ is broken as it is simply made of the two slit's lips. In both cases, the arc $\zeta$ makes a kink angle $\kappa$ at $q'$, with $\kappa=0$ if and only if $\beta=\pi$. 

\begin{figure}[h]
\centering
\begin{tikzpicture}
\node[anchor = south west, inner sep=0mm] at (0,0) {\includegraphics[width=8cm]{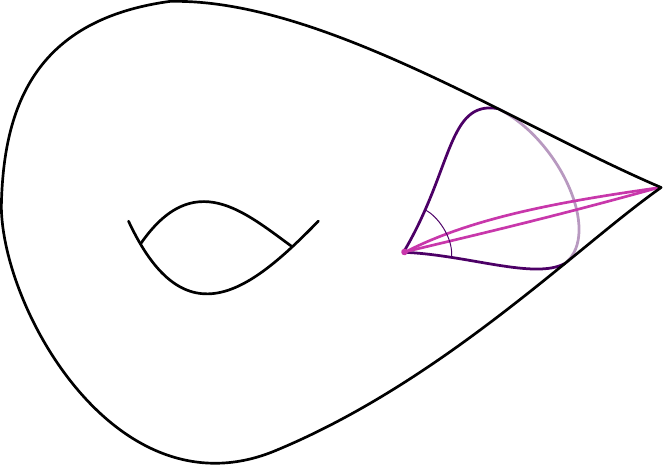}};

\node[plum] at (5.45,3.1) {$\kappa$};
\node[plum] at (5.2,3.75) {$\zeta$};
\node[mauve] at (4.7, 2.65) {$q'$};
\end{tikzpicture}
\caption{The geodesic arc $\zeta$ and the kink angle $\kappa$.}
\label{fig:kinked-geodesic}
\end{figure}

Now, given an isotopy class of branched hyperbolic structures $[h]\in R(\mathcal{S}^\textrm{north})$, recall that it comes with a preferred geometric presentation of $\pi_1(S,q)$ with generators $a_1,b_1,a_2,b_2$ satisfying $[a_1,b_1][a_2,b_2]=1$. In the homotopy class $[a_1,b_1]$ of loops on $S$ based at the branch point $q$, there is a unique geodesic arc representative $\overline{\zeta}$ that starts and ends at $q$. A priori, the arc $\overline{\zeta}$ could pass several times through $q$.
\begin{claim}
The arc $\overline{\zeta}$ meets $q$ only at its endpoints.
\end{claim}
\begin{proofclaim}
By assumption, the surface $S$ is the realization of a slit tori assembly made of two cone tori $\Sigma_1$ and $\Sigma_2$ with cone data $(p_1,\beta_1)$ and $(p_2,\beta_2)$. Since $\beta_1+\beta_2=2\pi$, it holds that $\min\{\beta_1,\beta_2\}\leq \pi$. This means that on either $\Sigma_1$ or $\Sigma_2$, there is a geodesic arc $\zeta$ as we have just described. After gluing, this arc becomes an unbroken geodesic arc on $S$ based at $q$ and lies in the homotopy class $[a_1,b_1]$. By uniqueness, it corresponds to $\overline{\zeta}$ and hence $\overline{\zeta}$ passes through $q$ only once.
\end{proofclaim}

The length $c$ and kink angle $\kappa\in (-\pi,\pi)$ at $q$ of the geodesic arc $\overline{\zeta}$ can be measured using $h$. We define $\kappa$ to be signed positively if the curve turns right at $q$, and negative if it turns left. In other words, $\kappa>0$ if $\beta_1<\pi$, while $\kappa<0$ if $\beta_2<\pi$, and $\kappa=0$ if $\beta_1=\beta_2=\pi$. The quantities $c$ and $\kappa$ are directly related to the angle $\beta=\min\{\beta_1,\beta_2\}$ and the slit length $\ell$. More precisely, when we cut the cone torus illustrated on Figure~\ref{fig:kinked-geodesic} along the slit and along $\zeta$, we obtain the isoceles triangle $T$ of Figure~\ref{fig:isoceles-triangle}. Hyperbolic trigonometry gives
\begin{align*}
    \cos(\beta)&=\sin(\kappa/2)^2\cosh(c)-\cos(\kappa/2)^2,\\
    \cosh(\ell)^2&=\frac{\cos(\kappa/2)^2}{\cos(\kappa/2)^2+\tanh(c/2)^2}.
\end{align*}
These two relations allow us to recover the values of $\beta$ (which is $\leq \pi$) and $\ell$ from those of $\kappa$ and $c$ which we can measure from $[h]$. The sign of $\kappa$ tells us which of $\beta_1$ or $\beta_2$ is equal to $\beta$, hence we also recover the values of $\beta_1$ and $\beta_2$.

\begin{figure}[h]
\centering
\begin{tikzpicture}
\node[anchor = south west, inner sep=0mm] at (0,0) {\includegraphics[width=5cm]{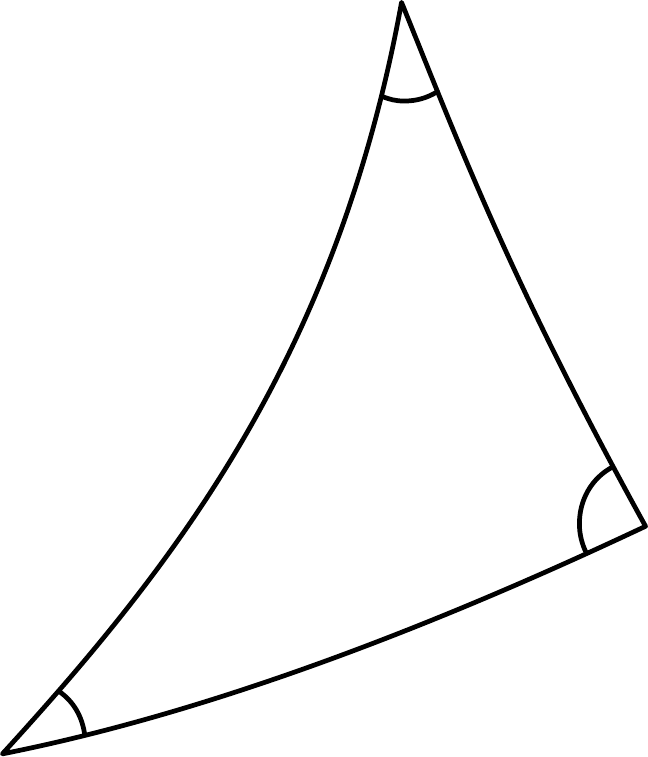}};

\node[plum] at (1.1,.6) {$|\kappa|/2$};
\node[plum] at (3.2,4.3) {$|\kappa|/2$};
\node[plum] at (1.9,3.1) {\Large $c$};
\node[mauve] at (3, .6) {\Large $\ell$};
\node[mauve] at (4.3, 3.7) {\Large $\ell$};
\node[sky] at (4.2, 2) {\Large $\beta$};
\end{tikzpicture}
\caption{The isoceles triangle $T$.}
\label{fig:isoceles-triangle}
\end{figure}

When we cut $S$ along $\overline{\zeta}$, we obtain two hyperbolic tori $\Sigma_+$ and $\Sigma_{-}$ with a boundary curve of length $c$ and corner points on that boundary curve of respective angles $2\pi+\kappa$ and $2\pi-\kappa$, as those illustrated on Figure~\ref{fig:tori-Sigma+-and-Sigma-}. When $\kappa\geq 0$, we may remove a copy of $T$ from $\Sigma_+$ to turn it into a hyperbolic tori with a boundary curve of length $2\ell$ and two corner points of angle $2\pi$ and $2\pi-\beta$. We can then glue together the two segments of length $\ell$ of the boundary curve of $\Sigma_+\setminus T$ that join the two corner points to obtain a hyperbolic cone torus with one cone singularity of angle $2\pi-\beta$. The glued image of the two segments defines a slit on that torus. Similarly, we can glue a copy of $T$ along the boundary of $\Sigma_-$, then glue the analogous two segments of length $\ell$ and obtain a slit torus with cone angle $\beta$. 

\begin{figure}[h]
\centering
\begin{tikzpicture}
\node[anchor = south west, inner sep=0mm] at (0,0) {\includegraphics[width=12cm]{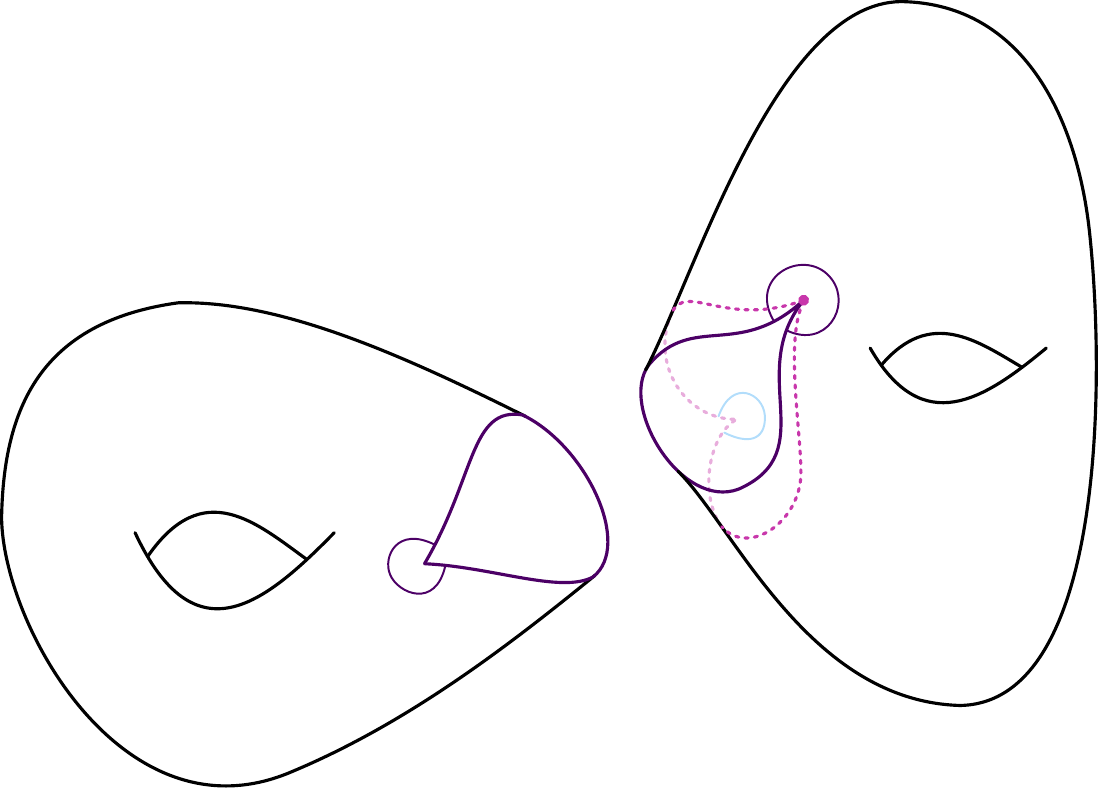}};

\node[plum] at (4.5,1.9) {$2\pi-\kappa$};
\node[plum] at (8.8,6) {$2\pi+\kappa$};
\node[plum] at (4.9,3.75) {$c$};
\node[plum] at (8.2,3.1) {$c$};
\node[mauve] at (9, 3.3) {$\ell$};
\node[mauve] at (7.9, 5.55) {$\ell$};
\node[lightsky] at (7.95,4.5) {\footnotesize $2\pi-\beta$};

\draw (6,5.5) edge[out=-90, in=210, ->] (7.5, 4);
\node[anchor=south] at (6,5.5) {$T$};
\end{tikzpicture}
\caption{On the left: the torus $\Sigma_-$. On the right: the torus $\Sigma_+$ with the copy of $T$ drawn inside it.}
\label{fig:tori-Sigma+-and-Sigma-}
\end{figure}

When $\kappa\leq 0$, we add a copy of $T$ to $\Sigma_+$ and remove it from $\Sigma_-$. Recall that when $\kappa=0$, the triangle $T$ is degenerate as $c=2\ell$ and so we may skip the step of removing/gluing a copy of $T$ to $\Sigma_+$ and $\Sigma_-$. The resulting two slit tori are framed by the generators $a_1,b_1$, respectively $a_2,b_2$, of $\pi_1(S,q)$, and thus define a slit tori assembly.

Finally, we define $R^{-1}([h])$ to be the slit tori assembly we just described and we obtain the desired map $R^{-1}\colon R(\mathcal{S}^\textrm{north})\to \mathcal{S}^\textrm{north}$. The construction guarantees that $R$ and $R^{-1}$ are inverses of each other.
\end{proof}

\subsubsection{Realizing holonomies}
The composition of the realization map from Proposition~\ref{prop:bow-tie-realization-map-injective} with the holonomy map~\eqref{eq:holonomy-map} gives a map
\[
\mathcal{S}^\textrm{north}\stackrel{R}{\longrightarrow} \Hyp_1(S)\stackrel{\hol}{\longrightarrow}\chi(S,\psl).
\]
The image of $R\circ\hol$ lands in the two components of $\chi(S,\psl)$ with Euler number $\pm 1$. More precisely, given a slit tori assembly in $\mathcal{S}^\textrm{north}$, its realization $[h]$ comes with a preferred geometric presentation of $\pi_1S$ with generators $a_1,b_1,a_2,b_2$ satisfying $[a_1,b_1][a_2,b_2]=1$. If $[\rho]=\hol([h])$, then by construction $\rho([a_1,b_1])$ is the peripheral monodromy of the first slit torus in the assembly we are considering and is therefore elliptic. This means that 
\[
\hol\circ R(\mathcal{S}^\textrm{north})\subset \mathcal{B},
\]
where $\mathcal{B}$ is the set of conjugacy classes of bow-tie representations introduced in Definition~\ref{defn:bow-tie-representations}.

\begin{thm}\label{thm:bow-tie-geometrization}
The map $\hol\circ R\colon\mathcal{S}^\mathrm{north}\to\mathcal{B}$ is surjective. In other words, every bow-tie representation is geometrizable (in the sense of Definition~\ref{defn:geometrizable-representations}) by slit tori assemblies.
\end{thm}
\begin{proof}
Let $[\rho]\in \mathcal{B}_c\subset\mathcal{B}$ be the conjugacy class of a bow-tie representation. Recall that here $c$ is a simple closed curve on $S$ and $[\rho]\in\mathcal{B}_c$ means that $\rho$ maps every fundamental group representative of $c$ to an elliptic element of $\psl$. Our goal is to prove that $[\rho]$ belongs to the image of $\hol\circ R$.

Let $a_1,b_1,a_2,b_2$ be a geometric generating set of $\pi_1S$ compatible with $c$, which means that $c$ is represented by $[a_1,b_1]$. Let $\mathfrak{T}([\rho])=(\ell_1,\tau_1,\beta,\gamma,\ell_2,\tau_2)$ be the corresponding coordinates of $[\rho]$, as in Theorem~\ref{thm:coordinates-bow-tie-representations}. 
\begin{claim}\label{claim:optimal-action-angle-coordinates-bow-tie}
Up to changing the geometric generators by applying an automorphism of $\pi_1S$ that maps $[a_1,b_1]$ to a conjugate of itself, we may assume that
\[
\tau_i\in[0,\ell_i)\quad \text{and}\quad \gamma\in [0,\pi).
\]
\end{claim}
\begin{proofclaim}
First, when we replace the generators $b_1$ by $b_1a_1^{m}$ and $b_2$ by $b_2a_2^{n}$ for some integers $m,n$, then only the twist coordinates $\tau_1$ and $\tau_2$ are affected, and they change by $+m\ell_1$, respectively $+n\ell_2$. This is a consequences of~\eqref{eq:Goldman-formula-one-holed-torus} and~\eqref{eq:twist-flow-one-holed-torus}, and the fact that these operations are iterations of the fundamental group automorphisms~\eqref{eq:generators-Mos(S)} from Appendix~\ref{apx:generators}. So, after sufficiently many iterations, we can make sure that $\tau_i\in [0,\ell_i)$ for $i=1,2$. 

Now, if we conjugate $a_2$ and $b_2$ by $[a_1,b_1]^n$, then we only modify the coordinate $\gamma$ by $+n(2\pi-\beta)$. Similarly, when we conjugate $a_1$ and $b_1$ by $[a_2,b_2]^m$, we only modify $\gamma$ by $+m\beta$. This is a consequence of Proposition~\ref{prop:gamma-is-dual-to-beta} and~\eqref{eq:twist-flow-separating-curve-genus-2}. Therefore, assuming that $\gamma\in[\pi,2\pi)$, if $\beta\geq \pi$, then $2\pi-\beta\leq \pi$ and iterating the first operation sufficiently many times will eventually guaranty that $\gamma\in [0,\pi)$. On the other hand, if $\beta\leq \pi$, then we achieve the same goal by iterating sufficiently many times the second operation.
\end{proofclaim}

By Claim~\ref{claim:optimal-action-angle-coordinates-bow-tie}, up to changing the geometric generators of $\pi_1S$, we may assume that the bow-tie of $[\rho]$ looks like the one illustrated on Figure~\ref{fig:bow-tie} in the sense that the two pentagons are convex because $\tau_i\in [0,\ell_i)$, and they do not overlap because $\gamma\in [0,\pi)$, except maybe for the segments $[BY_1']$ and $[BY_2]$ that overlap when $\gamma=0$.

Now, let $\Sigma_1$ and $\Sigma_2$ be the two one-holed tori obtained by cutting $S$ along $c$. The fundamental groups $\pi_1\Sigma_1$ and $\pi_1\Sigma_2$ are the two free subgroups of $\pi_1S$ respectively generated by $a_1,b_1$ and $a_2,b_2$. We consider two hyperbolic cone torus structures on $\Sigma_1$ and $\Sigma_2$ with respective cone angles $2\pi-\beta$ and $\beta$, and respective Fenchel--Nielsen coordinates $(\ell_1,\tau_1)$ and $(\ell_2,\tau_2)$. The holonomies of these two structures are $[\rho_1]\in\chi_t(\Sigma_1,\psl)$ and $[\rho_2]\in\chi_{-t}(\Sigma_2,\psl)$ with $t=-2\cos(\beta/2)$ (compare with~\eqref{eq:relative-character-variety-restrictions-rho_1-rho_2} and the discussion before). It holds that $\Tol(\rho_1)+\Tol(\rho_2)=\pm 1$; so, up to reversing the orientations of the hyperbolic cone torus structures on $\Sigma_1$ and $\Sigma_2$, we may assume that $\Tol(\rho_1)+\Tol(\rho_2)=\eu(\rho)$.

Next, we chose the slit parameters $(\ell,\phi_1)$ and $(\ell,\phi_2)$. Let $p_1$ and $p_2$ be the respective cone points of $\Sigma_1$ and $\Sigma_2$ which are also the starting points of the slits. There is an upper bound $L_{[\rho]}$ such that any slit of length 
\begin{equation}\label{eq:range-slit}
\ell \in (0,L_{[\rho]})
\end{equation}
in the northern hemitori of $\Sigma_1$ and $\Sigma_2$ never meet the geodesic representatives of $b_1$ and $b_2$. It is possible to compute an explicit value of $L_{[\rho]}$ in terms of $\ell_i$, $\tau_i$, and $\beta$ using hyperbolic trigonometry but we will not do this computation. We choose the slit angles $\phi_1\in [0,\pi-\beta/2)$ and $\phi_2\in(0,\beta/2]$ so that they satisfy
\begin{equation}\label{eq:slit-angle-relation}
 \gamma = \phi_1 + \beta/2 - \phi_2.   
\end{equation}
Since we are assuming $\gamma\in [0,\pi)$, we can always pick the slit angles this way. Note that if $\gamma=0$, we have to choose $\phi_1=0$ and $\phi_2=\beta/2$. When $\gamma\in (0,\pi)$, then there is continuum of possible choices for $\phi_1$ and $\phi_2$. We have now constructed a family of slit tori assemblies parametrized by the freedom we have in picking the slit length $\ell$ and the slit angles $\phi_1$ and $\phi_2$ according to~\eqref{eq:range-slit} and~\eqref{eq:slit-angle-relation}.
\begin{claim}
The holonomy of the realization of any slit tori assembly in this family is $[\rho]$.
\end{claim}
\begin{proofclaim}
    We start by considering the two pentagons $X_1Y_1B_1Y_1'X_1'$ and $X_2Y_2,B_2Y_2'X_2'$ associated to $[\rho_1]$ and $[\rho_2]$ by our construction from Section~\ref{sec:polygonal-model-one-holed-tori}. Since we are assuming that $\tau_1\in [0,\ell_1)$ and $\tau_2\in[0,\ell_2)$, those two pentagons are convex. The slits are fully contained in these pentagons because $\ell<L_{[\rho]}$.  Let $M_i$ be the midpoint of $[X_iX_i']$, and $H_i$ and $H_i'$ be two orthogonal projections of $M_i$ on the geodesic line $(X_i,Y_i)$ and $(X_i'Y_i')$, as on Figure~\ref{fig:unfolded-torus}. Now, we reflect the pentagon $H_iY_iB_iY_i'H_i'$ through the geodesic line $(H_iH_i')$ and connect the image of the endpoint of the reflected slit to $X_i$ and $X_i'$ using geodesic segments. These two segments and the reflected image of the slit determine three regions---which we call \emph{flaps}---in the reflected pentagons. We translate two flaps using the hyperbolic transformations $\rho_i(a_i)$ and $\rho(b_ia_ib_i^{-1})$ to obtain the orange octagon on Figure~\ref{fig:unfolded-torus}. The two sides of the octagon that are the images of the reflected slit by the translations are drawn as dotted mauve segments to highlight that they correspond to the slit. These two sides meet at the vertex $B_i$ and make an angle of $2\pi-\beta_i$.
    
    \begin{figure}[h]
    \centering
    \begin{tikzpicture}
    \node[anchor = south west, inner sep=0mm] at (0,0) {\includegraphics[width=12cm]{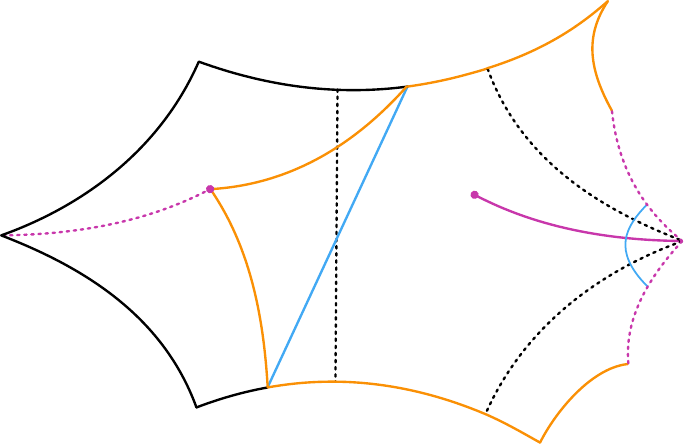}};

    \node at (4.75,.75) {$X_i$};
    \node at (5.9,.85) {$H_i$};
    \node at (8.5,.25) {$Y_i$};
    \node at (7.1,6.55) {$X_i'$};
    \node at (5.9,6.5) {$H_i'$};
    \node at (8.5,6.95) {$Y_i'$};
    \node at (12.25,3.55) {$B_i$};
    \node at (6.25,3.55) {$M_i$};
    \node[sky] at (10.3,3.32) {\small $2\pi-\beta_i$};

    \draw (3,3) edge[out=10, in=160, ->]  (10.3,1.9);
    \node at (8, 3) {\footnotesize $\rho(a_i)$};
    \draw (4,5.5) edge[out=-10, in=200, ->]  (9.8,6);
    \node at (8, 5.1) {\footnotesize $\rho(b_ia_ib_i^{-1})$};
    \end{tikzpicture}
    \caption{In orange, the octagon obtained after translation of the two flaps as indicated by the two arrows.}
    \label{fig:unfolded-torus}
    \end{figure}

    Once we have constructed the two octagons as above (one for $\Sigma_1$ and one for $\Sigma_2$), we can glue them by identifying opposite mauve dotted sides together. Since the octagons have complementary angles at the vertices $B_1$ and $B_2$ ($2\pi-\beta_1+2\pi-\beta_2=2\pi$), we end up with a large polygon $\mathcal{N}$ (that might overlap itself) immersed in the hyperbolic plane. This polygon $\mathcal{N}$ serves as a net for the branched hyperbolic structure on $S$. In other words, using obvious side identifications in $\mathcal{N}$, the flaps can be glued back together to produce the hyperbolic branched surface $S$ that is the realization of our slit tori assembly. We can alternatively develop $\mathcal{N}$ across its sides to obtain the developing map $\widetilde{S}\to \HH$ of the branched structure on $S$. This developing map is equivariant with respect to a representation $\rho'\colon\pi_1S\to\psl$.

    If we remove the six flaps (three for each of the two octagons) from $\mathcal{N}$ and only keep the two original pentagons $X_1Y_1B_1Y_1'X_1'$ and $X_2Y_2B_2Y_2'X_2'$, now attached at the vertex $B_1=B_2$, the resulting shape is a bow-tie. It is not hard to see that $\rho'$ is the representation associated to that bow-tie. By construction, the angle between the geodesic rays $[B_2Y_2)$ and $[B_1Y_1')$ is $\phi_1+\beta/2-\phi_2=\gamma$, so the bow-tie is isometric to the bow-tie of $\rho$. This means that $[\rho]=[\rho']$ as desired.
\end{proofclaim}
This concludes the proof of Theorem~\ref{thm:bow-tie-geometrization}.
\end{proof}

With the same notation as in the proof of Theorem~\ref{thm:geometrization-bow-tie-intro}, we can state the following result on isomodromic deformations.

\begin{cor}\label{cor:isomonodromic-deformations-bow-tie}
When $\gamma\neq 0$, the realization of the family of slit tori assemblies introduced in the proof of Theorem~\ref{thm:geometrization-bow-tie-intro} and parametrized by $\ell\in (0,L_{[\rho]})$ and $\phi_1-\phi_2=\gamma-\beta/2$ (see~\eqref{eq:range-slit} and~\eqref{eq:slit-angle-relation}) is a $2$-dimensional real subset of branched hyperbolic structures on $S$ with holonomy $[\rho]$.
\end{cor}
\begin{proof}
Different values of the slit length $\ell$ and slit angles $\phi_1,\phi_2$ correspond to different branched hyperbolic structures by Proposition~\ref{prop:bow-tie-realization-map-injective}. However, as we explained in the proof of Theorem~\ref{thm:geometrization-bow-tie-intro}, they all have holonomy $[\rho]$.
\end{proof}

\appendix
\setcounter{section}{0}
\renewcommand\sectionname{Appendix}
\addcontentsline{toc}{section}{Appendix}
\renewcommand{\thesection}{\Alph{section}}

\section{Mapping class group generators}\label{apx:generators}
\begin{figure}[h]
\centering
\begin{tikzpicture}
\node[anchor = south west, inner sep=0mm] at (0,0) {\includegraphics[width=8cm]{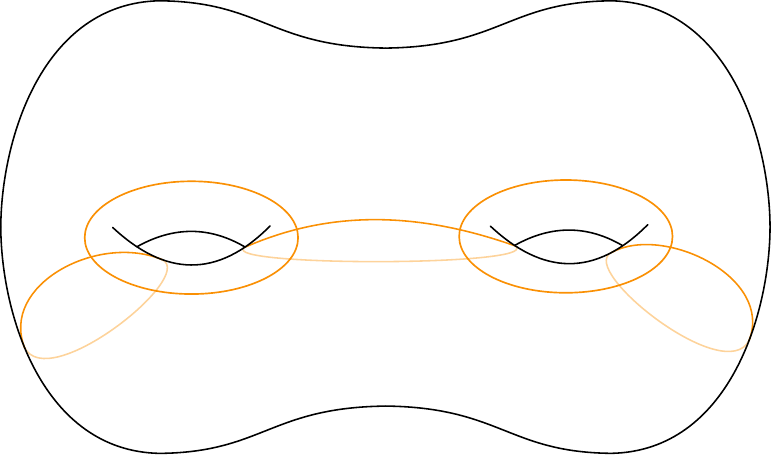}};
\end{tikzpicture}
\caption{Five simple closed curves whose corresponding Dehn twists are Lickorish's generators in genus $2$.}
\label{fig:Lickorish-generators}
\end{figure}

We recall some useful facts about the mapping class group $\Mod(S)$ of a closed and oriented surface $S$ of genus $2$. Mapping class groups are finitely generated by a particular class of elements called \emph{Dehn twists}. In genus $2$, an explicit generating family---called \emph{Lickorish's generating family}---is given by the Dehn twists along the five simple closed curves on $S$ depicted on Figure~\ref{fig:Lickorish-generators} (see~\cite[Theorem~4.13]{mcg-primer} for more details).

\begin{figure}[h]
\centering
\begin{tikzpicture}
\node[anchor = south west, inner sep=0mm] at (0,0) {\includegraphics[width=8cm]{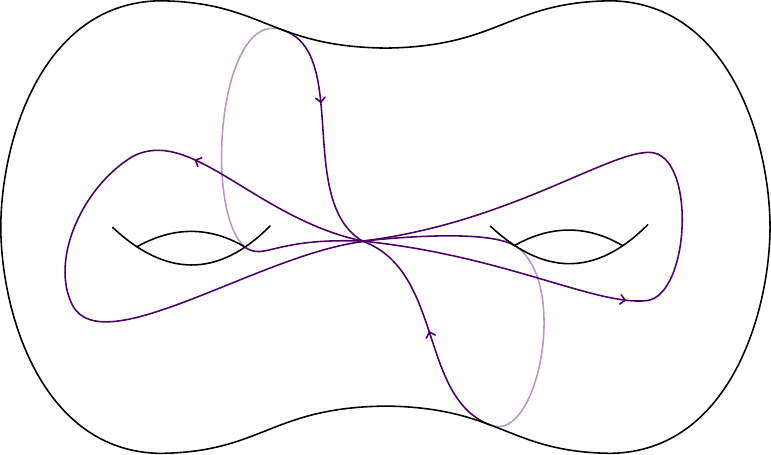}};

\node[plum] at (4.2,1) {$a_1$};
\node[plum] at (7.35,2.1) {$b_1$};
\node[plum] at (3.7,3.5) {$a_2$};
\node[plum] at (.4,2.1) {$b_2$};
\end{tikzpicture}
\caption{Fundamental group generators for a surface of genus $2$.}
\label{fig:generators-pi1S}
\end{figure}

The Dehn--Nielsen--Baer Theorem (see~\cite[Theorem~8.1]{mcg-primer}) identifies $\Mod(S)$ with an index $2$ subgroup of $\Out(\pi_1S)$---the group of outer automorphisms of $\pi_1S$. It is useful for our computations to work with explicit automorphisms of $\pi_1S$ representing Lickorish's generators. They can be expressed using a presentation of $\pi_1S$ of the form
\begin{equation}\label{eq:presentation-pi1S-appendix}
\pi_1S=\langle a_1,b_1,a_2,b_2 \, : \, [a_1,b_1][a_2,b_2]=1\rangle.
\end{equation}
The generators $a_1,b_1,a_2,b_2$ can be taken to be the homotopy classes of the loops on $S$ illustrated on Figure~\ref{fig:generators-pi1S}.

Our convention for fundamental group multiplication is the following: if $x,y\in\pi_1S$, then $xy$ stands for the homotopy class of the based loop obtained by going along a representative of $y$ first, and then along a representative of $x$. With this convention, the loops on Figure~\ref{fig:generators-pi1S} satisfy the fundamental group relation $[a_1,b_1][a_2,b_2]=1$, where $[x,y]=xyx^{-1}y^{-1}$

Working with the presentation~\eqref{eq:presentation-pi1S-appendix}, we can define five automorphisms of $\pi_1S$ which represent Lickorish's generators of $\Mod(S)$. Four of them are labeled $\tau_{a_1}, \tau_{b_1}, \tau_{a_2}, \tau_{b_2}$ because they are Dehn twists along the curves $a_1,b_1,a_2,b_2$. The last one is the Dehn twist along the curve $c=a_1a_2$ and is accordingly labeled $\tau_c$.
\begin{equation*}
    \tau_{a_1}\colon\begin{cases} a_1\mapsto a_1\\
    b_1\mapsto b_1a_1^{-1}\\
    a_2\mapsto a_2\\
    b_2\mapsto b_2
    \end{cases}
    \quad
    \tau_{b_1}\colon \begin{cases}
    a_1\mapsto a_1b_1\\
    b_1\mapsto b_1\\
    a_2\mapsto a_2\\
    b_2\mapsto b_2
    \end{cases}
    \quad
    \tau_{a_2} \colon\begin{cases}
    a_1\mapsto a_1\\
    b_1\mapsto b_1\\
    a_2\mapsto a_2\\
    b_2\mapsto b_2a_2^{-1}
    \end{cases}
    \quad 
    \tau_{b_2}\colon\begin{cases}
    a_1\mapsto a_1\\
    b_1\mapsto b_1\\
    a_2\mapsto a_2b_2\\
    b_2\mapsto b_2
    \end{cases}
\end{equation*}
\begin{equation}\label{eq:generators-Mos(S)}
    \tau_c\colon\begin{cases}
    a_1\mapsto a_1\\
    b_1\mapsto a_2a_1b_1\\
    a_2\mapsto a_2\\
    b_2\mapsto a_1a_2b_2
    \end{cases}
\end{equation}

\bibliographystyle{amsalpha}
\bibliography{bibliography.bib}
\end{document}